\documentclass[reqno, 10pt]{article}
\usepackage{amssymb, amsmath, amsthm, amsfonts, amscd, epsfig}
\usepackage{color}
\usepackage{subfigure}
\usepackage{bbm}
\usepackage{graphicx, epsfig}
\usepackage{geometry,graphicx,pgfplots,caption}
\usepackage{verbatim}
\usepackage{booktabs}

\newtheorem{theorem}{Theorem}[section]
\newtheorem{corollary}{Corollary}[section]
\newtheorem{lemma}{Lemma}[section]
\newtheorem{ass}{Assumption}[section]
\newtheorem{proposition}{Proposition}[section]
\newtheorem{definition}{Definition}[section]

\newtheorem{example}{Example}[section]
\newtheorem{remark}{Remark}[section]
{\bf}{\it}

\newcommand{\Bd}{\mathbf{d}}

\newcommand{\Om}{\Omega}

\newcommand{\Bx}{\mathbf{x}}

\newcommand{\hBx}{\hat{\mathbf{x}}}
\newcommand{\By}{\mathbf{y}}
\newcommand{\tBy}{\tilde{\mathbf{y}}}
\newcommand{\Bz}{\mathbf{z}}

\newcommand{\pa}{\partial}

\newcommand{\inc}{\textrm{i}}
\newcommand{\scat}{\textrm{s}}

\newcommand{\rmi}{\mathrm{i}}

\newcommand{\p}{\partial}

\newcommand{\ds}{\displaystyle}

\newcommand{\beq}{\begin{equation}}
	\newcommand{\eeq}{\end{equation}}
\newcommand{\RN}[1]{%
\textup{\uppercase\expandafter{\romannumeral#1}}%
}
\renewcommand{\d}{\mathrm{d}}
\numberwithin{equation}{section}
\numberwithin{figure}{section}

\begin{document}
	
\title{Solving Inverse Acoustic Obstacle Scattering Problem with Phaseless Far-Field Measurement Using Deep Neural Network Surrogates\thanks{The work of B. Jin is supported by Hong Kong RGC General Research Fund (Project 14306824), Hong Kong RGC  ANR / RGC Joint
Research Scheme (A-CUHK402/24) and a sart-up fund from The Chinese University of Hong Kong.}}
	
\author{Yuxin Fan\thanks{Department of Mathematics, The Chinese University of Hong Kong, Shatin, N.T., Hong Kong (\texttt{1155205162@link.cuhk.edu.hk}, \texttt{jihohong@cuhk.edu.hk}, \texttt{b.jin@cuhk.edu.hk})} \and Jiho Hong\footnotemark[2] \and Bangti Jin\footnotemark[2] \thanks{Author to whom any correspondence should be addressed}}
	
\maketitle
	
\begin{abstract}
In this work, we investigate the use of deep neural networks (DNNs) as surrogates for solving the inverse acoustic scattering problem of recovering a sound-soft obstacle from phaseless far-field measurements. We approximate the forward maps from the obstacle to the far-field data using DNNs, and for star-shaped domains in two and three dimensions, we establish the expression rates for fully connected feedforward neural networks with the ReLU activation for approximating the forward maps.  The analysis is based on the weak formulation of the direct problem, and can handle variable coefficients. Numerically we validate the accuracy of the DNN surrogates of the forward maps, and  demonstrate the use of DNN surrogates in the Bayesian treatment of the inverse obstacle scattering problem. Numerical experiments indicate that the surrogates are effective in both two- and three-dimensional cases, and can significantly speed up the exploration of the posterior distribution of the shape parameters using Markov chain Monte Carlo. \\
\textbf{Key words}:  inverse obstacle scattering, phaseless data, deep neural network, expression rate, Bayesian inversion
\end{abstract}
	
\section{Introduction}

Inverse acoustic scattering problems are concerned with determining the nature of an unknown scattering phenomenon, e.g., shape, locations, size, and medium properties, from the measurement of the scattered acoustic field, and have found a wide range of real-world applications, e.g., nondestructive evaluation, biomedical imaging, and microwave remote
sensing \cite{ODonnel:1981,ChandraZhou:2015}.
In practice, accurate phased data is usually difficult to acquire \cite[Chapter 8]{Chen:2018:book}: the accuracy of the phased measurement cannot be guaranteed especially for operating frequencies approaching the millimeter-waveband and beyond, so acquiring phase information is sophisticated and expensive, and moreover, the phase information is more susceptible to noise pollution than the amplitude information. In contrast,  collecting high-accuracy phaseless data is easier and cheaper. Phaseless inverse scattering rises also in quantum inverse scattering \cite{Klibanov:2014}. However, due to a lack of  phase information, phaseless reconstruction is more ill-posed and nonlinear than the phased counterpart. Therefore, it is of great importance to develop efficient and accurate computational techniques for inverse scattering problems with phaseless data.

In the literature, several numerical methods have been proposed for inverse acoustic obstacle scattering using phaseless data, which roughly can be categorized into two groups: regularized reconstruction and direct methods. The methods (see, e.g., \cite{Kress:1997:IOS,BaoLi:2013,AgaltsovHohage:2019,Zhang:2017:RSO}) in the former class are based on variational regularization \cite{ItoJin:2015}, by minimizing a discrepancy function that measures the mismatch between the predicted and measured phaseless far-field data plus a suitable penalty term, and can generally provide accurate reconstructions. However these methods utilize suitable optimization algorithms, which often take many iterations to reach convergence, require good initial guesses and can be expensive to deploy, especially in the three-dimensional case. In contrast, direct methods employ suitable indicator functions to indicate the presence of obstacles and have also been proposed for phaseless data, including linearized reconstruction via scattering coefficients \cite{HabibChowZou:2016}, reverse time migration \cite{ChenHuang:2017}, (approximate) factorization method \cite{ZhangZhang:2020siap}, and direct sampling method \cite{NingHanZhou:2025}, etc. Direct methods are computationally much more efficient but the reconstruction resolution is often modest.

In this work we investigate the use of deep neural networks (DNNs) to assist the task of recovering a sound-soft obstacle from phaseless far-field data. The approach consists of two steps: first construct DNN surrogates for the forward map(s), and then use the surrogates in exploring the posterior distribution of the shape parameters using Markov chain Monte Carlo in the Bayesian treatment of phaseless inverse obstacle scattering.  The main contributions of this work are as follows. First, we establish DNN expression rates for the shape-to-solution map  with a reference ball including the case of inhomogeneous coefficients. The key of the analysis is to establish the $(\boldsymbol{\beta},p,\varepsilon)$ shape holomorphy \cite{Chkifa:2015:BCD} of the forward maps. The analysis strategy employs the weak formulation of the problems (involving nonlocal operators \cite{Chandler:2008:WNEBTIS}) inside a region in which the shape deformation is conducted, and the approach is flexible and can handle the case of inhomogeneous coefficients. The result extends the work \cite{Dolz:2024:PSH}, in which D\"olz and Henr\'iquez proved the $(\boldsymbol{\beta},p,\varepsilon)$-holomorphy of the forward maps under an affine-parametric boundary transformation of the obstacle using a boundary integral formulation. Second, we conduct several numerical experiments in two- and three-dimensional cases to illustrate its potential. The numerical results in Section \ref{section:numerics} show that the approach can achieve significant speedup in Bayesian computation of inverse obstacle scattering, while the accuracy of the resulting approximate posterior distribution is comparable with the exact one based on the boundary element method. The comparative study with the generalized polynomial chaos expansion indicates that the DNN approach is superior in terms of both reconstruction accuracy and computational efficiency.

The last few years have witnessed significant progress on using DNNs to solve inverse scattering problems (see, e.g., \cite{wei2018deep,khoo2019switchnet,NingHanZou:2023,ChenJinLiu:2024,Zhou:2025:RTFVC} for an incomplete list). These techniques are constructed in various different ways, e.g., postprocessing \cite{wei2018deep,xu2020deep,NingHanZou:2023}, low-rank structure \cite{khoo2019switchnet}, Born approximation \cite{Zhou:2025:RTFVC}, unrolled optimization \cite{Deshmukh:2022,Zumbo:2024}, and learned regularizer \cite{ChenJinLiu:2024}, and have shown impressive empirical performance. Note that postprocessing type methods \cite{wei2018deep,xu2020deep,NingHanZou:2023} often require many paired training data, and the generalization property may suffer when tested on out-of-distribution data. In contrast, inversion methods that build on physical constraints, e.g., algorithmic unrolling \cite{Deshmukh:2022,Zumbo:2024}, often exploit the forward model and its adjoint operator and require less paired training data. Inspired by the well-established decomposition method in inverse scattering, Yin and Yan \cite{YinYan:2025sisc} proposed a novel physics-aware deep decomposition method for 2D acoustic obstacle scattering from limited aperture data, and it consists of data completion, Herglotz kernel  network and boundary recovery network, closely leveraging the scattering information. This approach was extended in \cite{YinYan:2025jcp} to the 3D case using transfer learning. See the reviews \cite{ChenWei:2020,GuoHuangEldar:2023} for in-depth discussions of deep learning-based techniques for inverse scattering. Several recent studies also explored the use of DNNs for inverse scattering with phaseless data \cite{xu2020deep,YinYangLiu:2021jcp,LimShin:2021,LuoWang:2021,Deshmukh:2022,ChenJinLiu:2024,NingHanZhou:2025}. The two-stage strategy is very popular: Xu et al \cite{xu2020deep} obtain initial estimates by, e.g., direct imaging / contrast source inversion which are then postporcessed using trained DNNs, whereas Ning et al \cite{NingHanZhou:2025} obtain rough estimates via a direct sampling method. Luo et al \cite{LuoWang:2021} first recover the phase information via phase retrieval net which is followed by a reconstruction net. Yin et al \cite{YinYangLiu:2021jcp} propose a two-layer sequence-to-sequence neural network (with parameters representing the boundary curve of obstacle) for recovering obstacle from the limited phaseless data. See also the work \cite{LimShin:2021} on using feedforward fully connected neural networks to predict discrete Fourier coefficients of a radially symmetric function from phaseless data. Deshmukh et al \cite{Deshmukh:2022} developed a reconstruction method based on unrolling gradient descent for a regularized objective function for inverse scattering with phaseless data. Chen et al \cite{ChenJinLiu:2024} propose a learned regularizer approach via latent representation for recovering the obstacle, including phaseless far field data. The present work complements these existing works on using DNNs as surrogates to solve phaseless inverse obstacle scattering problems in the Bayesian framework.

The rest of the paper is organized as follows. In Section \ref{section:prelim}, we describe the mathematical formulation of inverse obstacle scattering with phaseless data and the admissible set of shape parameters.
In Section \ref{section:holo:ext}, we establish shape holomorphy of the forward maps and derive the expression rates for DNN approximations.
Finally, in Section \ref{section:numerics}, we present numerical illustrations about approximating the forward maps and Bayesian treatment of inverse obstacle scattering with phaseless data. Throughout we denote by $\cdot$ the Euclidean inner product on $\mathbb{R}^d$, and by $|\cdot|$ the Euclidean norm of vectors. The notation $C$ denotes a generic constant which may change at each occurrence.
	
\section{Preliminaries}
In this section we describe phaseless inverse acoustic obstacle scattering and  admissible set of shape parameters for describing the obstacle.
\label{section:prelim}
\subsection{Inverse acoustic obstacle scattering from phaseless data}

Let $\Omega\subset \mathbb{R}^d$ $(d=2,3)$ be an open bounded Lipschitz domain with a connected complement $\mathbb{R}^d\setminus \overline{\Omega}$ and a boundary $\partial\Omega$. Physically $\Omega$ represents an impenetrable sound soft obstacle. Let
\begin{equation}\label{u:inc}
u^{\rm i}(\mathbf{x}\Bd) = e^{{\rm i}k\mathbf{d}\cdot\mathbf{x}}
\end{equation}
be a time harmonic incident plane wave, with $\rm i$, $k>0$ and $\mathbf{d}\in \mathbb{S}^{d-1}:=\{\mathbf{x}\in\mathbb{R}^d: |\mathbf{x}|=1\}$ being the imaginary unit, the fixed wave number and the incident direction, respectively. The interaction of the incident field $u^\inc$ and the obstacle $\Omega$ generates the scattered field $u^{\scat}$. The total field $u$ is given by $u=u^{\rm i}+u^{\scat}$, which  satisfies the following Helmholtz equation
\begin{equation}\label{HelmholtzEquation}
\left\{\begin{aligned}
	(\nabla\cdot\sigma(\Bx)\nabla + k^2 \tau(\Bx)) u &= 0,\quad\mbox{in }\mathbb{R}^d\backslash\overline{\Om},\\
	u &= 0,\quad\mbox{on }\partial\Omega.
\end{aligned}\right.
\end{equation}
The conductivity  $\sigma$ and the refractive index $\tau$ represent inhomogeneous density and compressibility, respectively, of the medium for pressure waves. We assume that $\sigma$ and $\tau$ are continuous, piecewise analytic, $\sigma\ge\sigma_0$ and $\tau\ge\tau_0$ for constants $\sigma_0,\tau_0>0$, and the functions $\sigma-1$ and $\tau-1$ are compactly supported. The scattered field $u^{\scat}$ satisfies the Sommerfeld radiation condition:
\begin{equation}\label{SmmerfeldRadiationCondition}
	\lim_{r\to\infty} r^{(d-1)/2}\left(\frac{\pa u^{\scat}(\mathbf{x})}{\partial r} - \rmi k u^{\scat}(\mathbf{x})\right) =0, \quad \text{with } r:= |\mathbf{x}|.
\end{equation}
Problem \eqref{HelmholtzEquation}--\eqref{SmmerfeldRadiationCondition} is well-posed for $u\in H_{\rm loc}^1(\mathbb{R}^d\backslash\overline{\Om})$ for the case $\sigma\equiv1$ \cite{Qu:2017:AFMIMS}.
The well-posedness when $\sigma$ is not a constant will be proved in Section \ref{sec:reformulation}. It is well known that the scattered field $u^{\rm s}$ satisfies the following asymptotic expansion (see, e.g., \cite[Lemma 2.5]{Chandler:2012:NABIM})
\begin{equation}\label{FarFieldPattern}
	u^{\scat}(\mathbf{x}) = r^{(1-d)/2}\mathrm{e}^{\rmi k r}\left( u^{\infty}(\hat{\mathbf{x}},\mathbf{d}) + O( r^{-1})\right),\quad \mbox{as } r\rightarrow\infty,
\end{equation}
which holds uniformly in all observation directions $\hat{\mathbf{x}}={\mathbf{x}}/{|\mathbf{x}|}$.
The function $u^{\infty}(\hat{\mathbf{x}},\mathbf{d}):\mathbb{S}^{d-1}\times\mathbb{S}^{d-1} \to \mathbb{C}$ is known as the far-field pattern, and denoted by $u^{\infty}[\Omega](\hat{\mathbf{x}},\mathbf{d})$ below to explicitly indicate its dependence on  $\Omega$. The concerned inverse problem is to recover the obstacle $\Omega$ from a knowledge of the phaseless far-field pattern $|u^\infty[\Omega](\hat{\mathbf{x}},\mathbf{d})|$.
When the background is homogeneous (i.e., $\sigma\equiv 1$ and $\tau\equiv1$), the phaseless far-field pattern is translation invariant, i.e., $|u^{\infty}[\Omega]|=|u^{\infty}[\mathbf{v}+\Om]|$ for any $\mathbf{v}\in\mathbb{R}^d$  \cite[Section 2]{Kress:1997:IOS}, which represents a natural obstruction to unique recovery. Since the seminal works \cite{Kress:1997:IOS,Klibanov:2014}, inverse scattering with phaseless  data has received much attention; see, e.g., \cite{HabibChowZou:2016,LiLiuWang:2017,Zhang:2017:RSO,ZhangGuo:2018,XuZhangZhang:2018,Sun:2019:UPI}. Several approaches have been proposed to break translation invariance and to ensure unique determination of $\Omega$ (see, e.g., \cite{Zhang:2017:RSO,ZhangGuo:2018,XuZhangZhang:2018}), by the superposition of incident plane waves and introduction of a reference ball into the scattering system.

We follow the approach of Zhang and Guo \cite{ZhangGuo:2018}. Specifically, let $\Gamma$ be the fundamental solution to the Helmholtz equation $(\Delta+k^2)u=0$ in $\mathbb{R}^d$ for $d=2,3$: for $\mathbf{x}\in\mathbb{R}^d\backslash\{0\}$,
$$\Gamma(\mathbf{x})=\begin{cases}
	\ds -\frac{\rm i}{4} H_0^{(1)}(k|\mathbf{x}|),&\mbox{if }d=2,\\[3mm]
	\ds -\frac{e^{{\rm i}k|\mathbf{x}|}}{4\pi|\mathbf{x}|},&\mbox{if }d=3,
\end{cases}
$$
where $H_0^{(1)}$ is the Hankel function of the first kind of order $0$. Now consider the incident field of the point source $v^{\inc}(\Bx;\Bz)=\Gamma(\Bx-\Bz)$, with
$\Bz$ located on the boundary $\partial P$ of a convex polygon $P\subset\mathbb{R}^d$ such that $\overline{P}\cap\overline{\Om}=0$.
Let $v^\infty[\Om]$ be the far field pattern for the point source. Then the inverse problem reads: given a reference ball $B$ and a convex polygon $P$ such that $\overline{B}\subset\mathbb{R}^d\backslash(\overline{P}\cup \overline{\Om})$, determine $\Omega$ from the following phaseless far-field data: $\bigl\{|u^\infty[\Omega\cup B](\widehat{\mathbf{x}},\mathbf{d})| : \widehat{\mathbf{x}}\in\mathbb{S}^{d-1}\bigr\}$ for a fixed $\mathbf{d}\in\mathbb{S}^{d-1}$, $\bigl\{|v^\infty[\Omega\cup B](\widehat{\mathbf{x}},\mathbf{z})| : \widehat{\mathbf{x}}\in\mathbb{S}^{d-1} \text{ and } \mathbf{z}\in\partial P\bigr\}$, and $\bigl\{|u^\infty[\Omega\cup B](\widehat{\mathbf{x}},\mathbf{d}) + v^\infty[\Omega\cup B](\widehat{\mathbf{x}},\mathbf{z})| : \widehat{\mathbf{x}}\in\mathbb{S}^{d-1} \text{ and } \mathbf{z}\in\partial P\bigr\}$.
Then uniqueness holds for the specific setting \cite[Theorem 3.1]{ZhangGuo:2018}, which is the focus of the present work.

\begin{theorem}[{\cite[Theorem 3.1]{ZhangGuo:2018}}]
\label{theorem:ZhangGuo}
Let $D_1$ and $D_2$ be open, simply connected, bounded domains with $C^2$ boundaries.
Let $B$ and $P$ be a ball and a convex polygon, respectively, such that $B$, $P$ and $D_j$ are pairwise disjoint for
 each $j=1,2$. Fix $\sigma\equiv 1$, and $\mathbf{d}\in\mathbb{S}^{d-1}$.
If the following relations hold
\begin{align*}
|u^\infty[D_1\cup B](\widehat{\Bx},\mathbf{d})|&=|u^\infty[D_2\cup B](\widehat{\Bx},\mathbf{d})|,\\
|v^\infty[D_1\cup B](\widehat{\Bx},\Bz)|&=|v^\infty[D_2\cup B](\widehat{\Bx},\Bz)|,\\
|u^\infty[D_1\cup B](\widehat{\Bx},\mathbf{d})+v^\infty[D_1\cup B](\widehat{\Bx},\Bz)|&=|u^\infty[D_2\cup B](\widehat{\Bx},\mathbf{d})+v^\infty[D_2\cup B](\widehat{\Bx},\Bz)|,\quad
\end{align*}
for all $\widehat{\Bx}\in\mathbb{S}^{d-1}$ and $\Bz\in\p P$, then there holds $D_1=D_2$.
\end{theorem}

\subsection{Admissible shape parameters}\label{sect:starshape}

Let $D$, $B$ and $P$ be a bounded Lipschitz domain, a reference ball and a convex polygon in $\mathbb{R}^d$, respectively.
We adopt the following assumption on $D$, $B$ and $P$, under which the far field patterns $u^{\infty}[D\cup B](\cdot,\Bd)$ and $v^\infty[D\cup B](\cdot,\Bz)$ are well defined for all $\Bd\in\mathbb{S}^{d-1}$ and $\Bz\in\p P$.
\begin{ass}\label{ass:separation}
There exist two balls $B_1$ and $B_2$ in $\mathbb{R}^d$ such that $\overline{B}\subset B_1$, $\overline{P}\subset B_2$, and 
the sets $\overline{B_1}$, $\overline{B_2}$ and $\overline{D}$ are pairwise disjoint.
\end{ass}

Now we define real-valued shape parameters that describe the domain $\Omega$ (so that $\Omega\cup B$ is a sound-soft obstacle and Assumption \ref{ass:separation} holds for $D=\Omega$). In this study, we focus on star-shaped obstacles, which have been extensively investigated in the numerical studies of phaseless inverse scattering \cite{Zhang:2017:RSO,YinYangLiu:2021jcp}. Note that every bounded, star-shaped domain with a $C^1$ boundary has the boundary with the following representation
\begin{equation*}
\left\{\Bx\in\mathbb{R}^d\,:\,|\Bx-\Bx_0|=e^{\rho(\frac{\Bx-\Bx_0}{|\Bx-\Bx_0|})}\right\},
\end{equation*}
where $\Bx_0\in\mathbb{R}^d$ is an interior point, and $\rho$ is a $C^1$ real-valued function on $\mathbb{S}^{d-1}$.
To parameterize $\rho$, we employ the natural orthonormal basis of $L^2(\mathbb{R}^{d-1})$, i.e., Fourier basis on $\mathbb{S}^{1}$ and real-valued spherical harmonics on $\mathbb{S}^{2}$. More precisely, we define $X_0:=(2\pi)^{-1/2}$, and for $m\in\mathbb{N}$,
\begin{align*} X_m(\varphi):=\pi^{-1/2}\cos (m\varphi)\quad\mbox{and}\quad X_{-m}(\varphi):=
\pi^{{-1/2}}\sin (m\varphi).
\end{align*}
We also define, for all $\ell\in\mathbb{N}\cup\{0\}$,
	$$Y_{\ell,m}(\theta,\varphi):=\begin{cases}
		\sqrt{2} a_{\ell,m} P_\ell^m(\cos\theta)\cos (m\varphi),&m=1,\dots,\ell,\\
		\sqrt{2} a_{\ell,m} P_\ell^{|m|}(\cos\theta)\sin (|m|\varphi),&m=-1,\dots,-\ell,\\
		a_{\ell,0} P_{\ell}^0(\cos\theta),&m=0,
	\end{cases}$$
where $a_{\ell,m}:=\left(\frac{2\ell+1}{4\pi}\frac{(\ell-|m|)!}{(\ell+|m|)!}\right)^{1/2}$, $P_\ell^m$ is the Legendre polynomial of degree $\ell$ and order $m$ (cf. Definition \ref{def:Legendre}), and $\mathbb{S}^{d-1}$ is parametrized by $(\cos\varphi,\sin\varphi)$ for $d=2$ and $(\sin\theta\cos\varphi,\sin\theta\sin\varphi,\cos\theta)$ for $d=3$ with $0\le \theta<\pi$ and $0\le \varphi<2\pi$.
Then the family of functions
$\{X_m\,:\,m\in\mathbb{Z}\}$ ($d=2$) and $\{Y_{\ell,m}\,:\,\ell\in\mathbb{N}\cup\{0\},\ m\in\mathbb{Z}\mbox{ and }|m|\le\ell\}$ ($d=3$)
is an orthonormal basis of $L^2(\mathbb{S}^{d-1})$.
Let $\mathbf{a}\in\mathbb{R}^d$, $N\in\mathbb{N}$, $\mathbf{b}\equiv (b_{m})_{|m|\le N}\in\mathbb{R}^{2N+1}$ for $d=2$, and $\mathbf{b}\equiv(b_{\ell,m})_{0\le\ell\le N,|m|\le\ell}\in\mathbb{R}^{(N+1)^2}$ for $d=3$.
Let $\Om(\mathbf{a},\mathbf{b})$ be the domain with the boundary
\begin{equation}\label{eq:define:parameters}
\p\Om(\mathbf{a},\mathbf{b})=\left\{\Bx\in\mathbb{R}^{n}\,:\,|\Bx-\mathbf{a}|=
e^{\rho_{\mathbf{b}}(\frac{\Bx-\mathbf{a}}{|\Bx-\mathbf{a}|})}\right\},
\end{equation}
with
\begin{equation}\label{eq:define:parameters-2}
\rho_{\mathbf{b}}=\begin{cases}
\sum_{m=-N}^N b_{m} X_{m}&\mbox{if }d=2,\\
\sum_{\ell=0}^N\sum_{m=-\ell}^\ell b_{\ell,m} Y_{\ell,m}&\mbox{if }d=3.
\end{cases}
\end{equation}

Also consider the limit case $N=\infty$: Let $\mathbf{a}\in\mathbb{R}^d$, $\mathbf{b}\equiv (b_{m})_{m\in\mathbb{Z}}\in\mathbb{R}^{\infty}$ for $d=2$, and $\mathbf{b}\equiv(b_{\ell,m})_{\ell\in\mathbb{N}\cup\{0\},|m|\le\ell}\in\mathbb{R}^{\infty}$ for $d=3$ satisfy the following assumption.
\begin{ass}\label{ass:regularity}
    There exist $C>0$ and $q>0$ such that $|b_m|\le C(1+|m|)^{-2-q}$ for all $m\in\mathbb{Z}$ in 2D, and $|b_{\ell,m}|\le C\ell^{-3/2-q}{\ell+|m|\choose 2|m|}^{-1/2}\frac{(|m|+1)^{1/4}}{\ell^2 + |m|^2-\ell|m|}$ for all $\ell\in\mathbb{N}$ and $-\ell\le m\le\ell$ in 3D.
\end{ass}

Assumption \ref{ass:regularity} imposes suitable decay property of the expansion coefficients in order to ensure sufficient regularity of the boundary $\partial\Omega$; see the following proposition for the precise regularity. Several existing theoretical studies require $C^2$ regularity of the boundary $\partial\Omega$ (see, e.g., \cite{XuZhangZhang:2018,ZhangGuo:2018}). Hence, the imposed $C^1$ regularity is not very restrictive.
It always holds when $N<\infty$. Assumption \ref{ass:regularity} in 3D for the cases $m=0$ and $m=\ell$ are $|b_{\ell,0}|\le C\ell^{-7/2-q}$ and $|b_{\ell,\ell}|\le C\ell^{-13/4-q}$ for all $\ell\in\mathbb{N}$. The proof of the next result is given in Section \ref{appendix:Smoothness}.
\begin{proposition}\label{prop:regularity}
    Under Assumption \ref{ass:regularity}, the function $\rho_{\mathbf{b}}$ defined by the limit of \eqref{eq:define:parameters-2} as $N\to\infty$ is $C^1$.
\end{proposition}

\begin{definition}
The pair $(\mathbf{a},\mathbf{b})$ is said to be admissible if it satisfies Assumption \ref{ass:regularity}, and the set $\overline{\Om(\mathbf{a}, \mathbf{b})}$ is disjoint from $\overline{B_1}$ and $\overline{B_2}$.     Let $\mathcal{A}_N(B_1,B_2)$ be the set of all admissible pairs.
\end{definition}

The far-field responses $u^{\infty}[\Om(\mathbf{a},\mathbf{b})\cup B](\cdot,\Bd)$ and $v^{\infty}[\Om(\mathbf{a},\mathbf{b})\cup B](\cdot,\Bz)$ are well defined for all admissible tuples $(\mathbf{a},\mathbf{b})\in\mathcal{A}_N(B_1,B_2)$, $\Bd\in\mathbb{S}^{d-1}$ and $\Bz\in\p P$.

\subsection{Proof of Proposition \ref{prop:regularity} }
\label{appendix:Smoothness}
We only give the proof in the 3D case, since the 2D case is similar and simpler. The proof is based on the facts that $X_m$ and $Y_{\ell,m}$ are $C^1$ functions and that the series defining $\rho_{\mathbf{b}}$ have uniformly convergent term-by-term derivatives under Assumption \ref{ass:regularity}. The proof uses the Legendre and Jacobi polynomials, defined using Rodrigues' formula.
\begin{definition}\label{def:Legendre}
The Legendre polynomials $P_\ell^m$ are defined by
\begin{equation*}
P_\ell^m(x)=\frac{(-1)^m}{2^\ell \Gamma(\ell+1)}(1-x^2)^{{m}/{2}}\frac{\d^{\ell+m}}{\d x^{\ell+m}}(x^2-1)^\ell,\quad  \forall \ell\in\mathbb{N}\cup\{0\}, -\ell\le m\le \ell.
\end{equation*}
The Jacobi polynomials $P_\gamma^{(\alpha,\beta)}$ are defined by
\begin{equation*}
 P_\gamma^{(\alpha,\beta)}(x)=\frac{(-1)^\gamma}{2^\gamma \Gamma(\gamma+1)}(1-x)^{-\alpha}(1+x)^{-\beta}\frac{\d^\gamma}{\d x^\gamma}\{(1-x)^\alpha (1+x)^\beta (1-x^2)^\gamma\}.
\end{equation*}
\end{definition}

The proof of Proposition \ref{prop:regularity} uses crucially the following lemma.
	\begin{lemma}[{\cite[Theorem 7.32.1]{Szego:1975:OP}}]
		\label{lemma:max}
		Let $\alpha>-1$, $\beta>-1$ and $n\in\mathbb{N}\cup\{0\}$. If $\max(\alpha,\beta)\ge-1/2$,
		the Jacobi polynomial $P_\gamma^{(\alpha,\beta)}$ satisfies
		$$\max\{|P_\gamma^{(\alpha,\beta)}(x)|\,:\,-1\le x\le 1\}={\gamma+\max(\alpha,\beta) \choose \gamma}=\frac{\Gamma(\gamma+\max(\alpha,\beta)+1)}{\Gamma(\gamma+1)\Gamma(\max(\alpha,\beta)+1)}.$$
	\end{lemma}

Now we can give the proof of Proposition \ref{prop:regularity} in the 3D case.
\begin{proof}
	Let $\ell\in\mathbb{N}$ and $m\in\{0,1,\dots,\ell\}$.
Using the relations
\begin{align*} P_\ell^0(x)&=P_\ell^{(0,0)}(x),\quad \frac{\d^m}{\d x^m}P_\ell^{(0,0)}(x)=\frac{\Gamma(\ell+m+1)}{2^m\Gamma(\ell+1)}P_{\ell-m}^{(m,m)}(x),\\
P_\ell^m(x)&=(-1)^m(1-x^2)^{{m}/{2}}\frac{\d^m}{\d x^m}P_\ell^0(x),
\end{align*}
we obtain
$$|P_\ell^m(x)|=(1-x^2)^{{m}/{2}}\frac{\Gamma(\ell+m+1)}{2^m\Gamma(\ell+1)}|P_{\ell-m}^{(m,m)}(x)|,\quad \forall x\in[-1,1].$$
This and Lemma \ref{lemma:max} with $\alpha=\beta=m$ and $\gamma=\ell-m$ imply \begin{align*}
\max_{x\in[-1,1]}|P_\ell^m(x)|&\le\frac{\Gamma(\ell+m+1)}{2^m\Gamma(\ell+1)}\max_{x\in[-1,1]}| P_{\ell-m}^{(m,m)}(x)|\\
&=\frac{\Gamma(\ell+m+1)}{2^m\Gamma(m+1)\Gamma(\ell-m+1)}.
\end{align*}
It follows from the inequality
\begin{align*}
\left|\frac{\d}{\d x}P_\ell^m(x)\right|&=\left|\frac{\d}{\d x}((1-x^2)^{{m}/{2}}\frac{\d^m}{\d x^m}P_\ell^0(x))\right|\\
&\le m\left|\frac{\d^m}{\d x^m}P_\ell^0(x)\right|+\left|\frac{\d^{m+1}}{\d x^{m+1}}P_\ell^0(x)\right|
\end{align*}
that	
\begin{align*}
\max_{x\in[-1,1]}\left|{\frac{\d}{\d x}}P_\ell^m(x)\right|\le& m\frac{\Gamma(\ell+m+1)}{2^m\Gamma(\ell+1)}\max_{x\in[-1,1]}| P_{\ell-m}^{(m,m)}(x)|\\
&+\frac{\Gamma(\ell+m+2)}{2^{m+1}\Gamma(\ell+1)}\max_{x\in[-1,1]}|P_{\ell-m-1}^{(m+1,m+1)}(x)|\\
=& m\frac{\Gamma(\ell+m+1)}{2^m\Gamma(\ell-m+1)\Gamma(m+1)}\\
&+\frac{\Gamma(\ell+m+2)}{2^{m+1}\Gamma(\ell-m)\Gamma(m+2)},\quad\mbox{if }m<\ell,
\end{align*}
and
$$\max_{x\in[-1,1]}|{\textstyle\frac{\d}{\d x}}P_\ell^\ell(x)|\le \frac{\ell\Gamma(2\ell+1)}{2^\ell\Gamma(\ell+1)}.
$$
This and the estimate $${2m\choose m}\le2^{2m}(2m+1)^{-\frac12}$$
imply that for all $0\le m\le\ell-1$,
	\begin{align*}
        \max_{\mathbf{y}\in\mathbb{S}^{n-1}}|\nabla Y_{\ell,m}(\mathbf{y})|&\le C|a_{\ell,m}|\left(m\max_{x\in[-1,1]}|P_\ell^m(x)|+\max_{x\in[-1,1]}|{\frac{\d}{\d x}}P_\ell^m(x)|\right)\\
	&\le C \frac{1}{2^{m}}\left(\ell{\ell+m\choose 2m}{2m\choose m}\right)^{\frac12}\left(m+\frac{(\ell+m+1)(\ell-m)}{2(m+1)}\right)\\
    &\le C \ell^{1/2}{\ell+m\choose 2m}^{\frac12}(2m+1)^{-\frac14}\left(m+\frac{\ell(\ell-m)}{m+1}\right)
	\end{align*}
	 and $$\max_{\mathbf{y}\in\mathbb{S}^{d-1}}|\nabla Y_{\ell,\ell}(\mathbf{y})|\le C\ell^{1/2}{2\ell\choose2\ell}^{1/2}\ell^{1-1/4}= C\ell^{{5}/{4}}.$$
Therefore, under Assumption \ref{ass:regularity},
we have
\begin{align*}
\sum_{\ell=0}^\infty\sum_{m=-\ell}^\ell |b_{\ell,m} \nabla Y_{\ell,m}|&\le C\sum_{\ell=0}^\infty(\ell+1)^{-1-q}\sum_{m=-\ell}^{\ell}\frac{1}{|m|+1}\\
&\le C\sum_{\ell=0}^\infty(\ell+1)^{-1-q}(1+\log(\ell+1))<\infty.
\end{align*}	
Since $Y_{\ell,m}\in C^1(\mathbb{S}^{2})$ for all $\ell,m$, the series of functions $\sum_{\ell=0}^\infty\sum_{m=-\ell}^\ell b_{\ell,m} \nabla Y_{\ell,m}$ converges absolutely and uniformly on $\mathbb{S}^2$. Similarly, one can prove that $\sum_{\ell=0}^\infty\sum_{m=-\ell}^\ell b_{\ell,m} Y_{\ell,m}$ converges absolutely and uniformly on $\mathbb{S}^2$. Thus we deduce  $\rho_{\mathbf{b}}\in C^1(\mathbb{S}^{2})$.
\end{proof}

\section{Expression rates of DNN approximations} \label{section:holo:ext}
In this section, we establish expression rates of DNN approximations of the forward maps, which represent the main technical novelty of the study. The key in the analysis is the shape holomorphy of the forward maps. The discussion focuses on the plane wave excitation. The case of point source excitation can be analyzed similarly and is given in the appendix.
The analysis is lengthy and technical, and thus we provide a brief summary  of  the overall analysis strategy of establishing the expression rates.
\begin{itemize}
\item[Step 1.] First we reformulate the forward problem on an unbounded domain into an equivalent problem on a bounded domain which involves the Dirichlet-to-Neumann map and prove its well-posedness (cf. Lemma \ref{lemma:iso:pw}).
\item[Step 2.] We define the forward map on a class of domain transformations, and then extend the map to a complex Banach space based on the reformulation in Step 1 (see Lemma \ref{lemma:welldefined:complex} for its well-definedness).
\item[Step 3.] We prove that the extended forward map is complex Fr\'{e}chet differentiable (cf. Theorem 
\ref{theorem:holomorphy:planewave}), and then construct the parametric forward map that has $(\boldsymbol{\beta},p,
\varepsilon)$-holomorphic property (cf. Lemma \ref{lemma:paratransform:holomorphy}).
\item[Step 4.] We establish the expression rate using the $(\boldsymbol{\beta},p,\varepsilon)$-holomorphic property in the existing literature.
\end{itemize}
In particular, by carefully choosing the parameters $\boldsymbol{\beta}$ and $p$ in the analysis, we can quantify the impact of the boundary regularity of the obstacle on the expression rate; see Remark \ref{remark:convergencrate} for more details.

\subsection{Reformulation of the forward problems}
\label{sec:reformulation}
To deal with the unbounded domain $\mathbb{R}^d\setminus\Omega$, following  \cite{Chandler:2008:WNEBTIS,Kirby:2021:nonlocalBC}, we replace the Sommerfeld radiation condition \eqref{SmmerfeldRadiationCondition}  by an equivalent nonlocal boundary condition on a bounded domain. 	
Let $B_1$ and $B_2$ satisfy Assumption \ref{ass:separation}. Since $\sigma-1$ and $\tau-1$ are compactly supported, there exist two balls $B_{\rm i},B_{\rm o}\subset \mathbb{R}^d$ satisfying $\sigma=\tau=1$ in $\mathbb{R}^d\backslash\overline{B_{\rm i}}$, $\overline{B_1}\cup\overline{B_2}\cup\overline{D}\subset B_{\rm i}$ and $\overline{B_{\rm i}}\subset B_{\rm o}$.
Fix any such $B_{\rm i}$ and $B_{\rm o}$ below.
Let $u$ be the solution of problem \eqref{u:inc}-\eqref{SmmerfeldRadiationCondition}.
Let $\Lambda :H^{1/2}(\p B_{\rm o})\to H^{-1/2}(\p B_{\rm o})$ be the (exterior) Dirichlet-to-Neumann map defined by $\Lambda u = -\partial_\nu \widetilde{u}$ for all $u\in H^{1/2}(\p B_{\rm o})$, where $\nu$ is the unit outward normal vector to $B_{\rm o}$ and $\widetilde{u}$ is the solution of
$$
\left\{\begin{aligned}	(\Delta+k^2)\widetilde{u}&=0,\quad\mbox{in }\mathbb{R}^d\backslash\overline{B_{\rm o}},\\
\widetilde{u}&=u,\quad\mbox{on }\p B_{\rm o},
\end{aligned}\right.
$$
with the Sommerfeld radiation condition \eqref{SmmerfeldRadiationCondition}.
For the well-definedness of the nonlocal operator $\Lambda$, see, e.g., \cite[Theorem 2.31]{Chandler:2012:NABIM}.
Then the restriction of the solution $u$ to problem \eqref{HelmholtzEquation}-\eqref{SmmerfeldRadiationCondition} with $\Om=B\cup D$ and the incident field \eqref{u:inc} to the region $E:=B_{\rm o}\backslash(\overline{D}\cup\overline{B})$ satisfies
\begin{equation}\label{eq:bddomain:reformulation:u}
\left\{
\begin{aligned}
(\nabla\cdot\sigma(\Bx)\nabla + k^2\tau(\Bx)) u &= 0,\quad\mbox{in }E,\\
	u &= 0,\quad\mbox{on }\partial D \cup\p B,\\		
    \partial_\nu(u-u^{\inc}) &= -\Lambda (u-u^{\inc}),\quad\mbox{on }\partial B_{\rm o}.
\end{aligned}\right.
\end{equation}
Let
$
\mathcal{H}=\{f\in H^{1}(E)\,:\, f|_{\p B\cup \p D}=0\}.
$
The weak formulation of  \eqref{eq:bddomain:reformulation:u} reads: find $u\in H^1(E)$ such that
\begin{equation}\label{eq:reform:u}
a(u,w) = b(w),\quad\forall w\in \mathcal{H},	
\end{equation}
with the sesquilinear form $a$ and linear form $b$ given respectively by
\begin{align*}
a(u,w)&:=\int_{E}\sigma\nabla u\cdot\nabla \overline{w}- k^2\tau u \overline{w}~\d\mathbf{x} +\int_{\p B_{\rm o}}\sigma(\Lambda u)\overline{w}~\d s\\
b(w)&:=\int_{\p B_{\rm o}}(\Lambda u^{\rm i} + \partial_\nu u^{\rm i})\overline{w}~\d s.
\end{align*}
Note that $a$ is bounded on $\mathcal{H}$. We define a map $A:\mathcal{H}\to \mathcal{H}'$ by
$\langle A f,g\rangle=a(f,g)$, for all $f,g\in \mathcal{H}.$
Then problem \eqref{eq:reform:u} is well-posed. Below the notation $(\cdot,\cdot)_{L^2(D)}$ denotes either the $L^2(D)$ inner product or duality product.
\begin{lemma}\label{lemma:iso:pw}
$A:\mathcal{H}\to \mathcal{H}'$ is an isomorphism.
\end{lemma}
\begin{proof}
For all $u\in\mathcal{H}$, there holds
$$\langle A u, u \rangle = \|\sqrt{\sigma}\nabla u\|_{L^2(E)}^2 - k^2(\tau u,u)_{L^2(E)}+(\Lambda u,u)_{L^2(\partial B_{\rm o})}.$$
By \cite[Corollary 3.1]{Chandler:2008:WNEBTIS}, we have $\Re\left[(\Lambda u,{u})_{L^2(\p B_{\rm o})}\right]\ge0$. Thus Garding's inequality holds
$$\Re\left[\langle Au, u \rangle\right]\ge \sigma_0\|u\|_{H^1(E)}^2 - (k^2\|\tau\|_\infty+\sigma_0)\|u\|_{L^2(E)}^2,\quad\forall u\in\mathcal{H}.$$
Since $H^1(E)$ is compactly embedded in $L^2(E)$, $A$ is an isomorphism by the Fredholm alternative  \cite[Theorem 5.4.5]{Nedelec:2001:AEE}.
\end{proof}

\subsection{Domain transformation}
Fix a bounded Lipschitz domain $\widehat{D}$, a reference ball $B$ and a polygon $P$ satisfying Assumption \ref{ass:separation} for $D=\widehat{D}$. Fix also $B_1$, $B_2$, $B_{\rm i}$ and $B_{\rm o}$ as in Section \ref{sec:reformulation}.
Consider a family $\mathcal{T}$ of domain transformations ${T}:B_{\rm o}\to B_{\rm o}$ satisfying the following assumption:
\begin{ass}\label{ass:T}
$T$ is bijective on $B_{\rm o}$, $T$ and $T^{-1}$ are Lipschitz continuous, and $T$ coincides with the identity map on $B_1\cup B_2\cup(B_{\rm o}\backslash\overline{B_{\rm i}})$.
\end{ass}

Moreover, we also assume the following condition. Let $J_T\in L^\infty(B_{\rm o},\mathbb{C}^{d\times d})$ be the Jacobian of $T$, whose entries are weak derivatives of $T\in W^{1,\infty}(B_{\rm o},\mathbb{C}^d)$.
\begin{ass}\label{ass:T:compact}
    $\mathcal{T}$ is a compact subset of $W^{1,\infty}(B_{\rm o},\mathbb{R}^d)$.
\end{ass}

Let $D_T=T(\widehat{D})$ for all $T\in\mathcal{T}$. We denote by $u_T$ the solution to problem \eqref{eq:reform:u} in the domain $D=D_T$. For all $T\in\mathcal{T}$, $\widehat{u}_T:=u_T\circ T$ and $E_T:=B_{\rm o}\backslash(\overline{D_T}\cup\overline{B})$ (with the shorthand $E=E_{\rm id}$) so that
${u}_T\in \mathcal{H}_T:=\{f\in H^{1}(E_T)\,:\, f|_{\p B\cup \p D_T}=0\}$ and $\widehat{u}_T\in\mathcal{H}\equiv \mathcal{H}_{\rm id}$ for any $ T\in\mathcal{T}$.
Then $\widehat{u}_T$ satisfies the following variational problem
\begin{align}
\label{eq:varform:extended:u}	a_T(\widehat{u}_T,w) = b(w), \quad\forall w\in \mathcal{H},
\end{align}
with the sesquilinear form $a $ given by
\begin{align*}
a_T(u,w):=&\int_{\p B_{\rm o}}(\Lambda u)\overline{w}~\d s+\int_{E} (\sigma\circ T)|J_T| J_T^{-1}J_T^{-\top}\nabla u\cdot\nabla \overline{w}~\d\mathbf{x} \\
&-k^2\int_{E} (\tau\circ T)|J_T| u \overline{w} ~\d\mathbf{x},
\end{align*}
and $b(w)$ given in \eqref{eq:reform:u} with $D=\widehat{D}$.
Let $A_T:\mathcal{H}\to \mathcal{H}'$ be the operator induced by $a_T$.
\begin{lemma}\label{lem:isomorphism}
If the operator $A_{\rm id}:\mathcal{H}\to \mathcal{H}'$ is an isomorphism, then so is $A_T$ for every $T\in\mathcal{T}$.
\end{lemma}

Next, we establish the holomorphic extension of the map $T\mapsto\widehat{u}_T$ in the following sense.
\begin{definition}[{\cite[Definition 13.1, Theorem 14.7]{Mujica:1986:CAB}}]\label{definition:holomorphic}
Let $X$ and $Y$ be two complex Banach spaces and $X_0$ be any nonempty open subspace of $X$. A map $f:X_0\to Y$ is called \emph{complex Fr\'{e}chet-differentiable} or \emph{holomorphic} if, for every $x\in X_0$, there exists a bounded linear operator $f'(x):X\to Y$ satisfying
\begin{equation*} \lim_{\varepsilon\to0,\,\varepsilon\in\mathbb{C}}\frac{\|f(x+\varepsilon x')-f(x)-\varepsilon f'(x)[x']\|_Y}{\varepsilon}=0,\quad\forall x'\in X.
\end{equation*}
\end{definition}

In order to extend the forward map $T\mapsto \widehat{u}_T$, we extend the piecewise linear functions $\sigma$ and $\tau$  defined 
in $\mathbb{R}^d$ to complex variables by: for all $\mathbf{z}\in\mathbb{C}^d$,
\begin{equation}\label{eq:extensions:sigmatau}
\sigma(\mathbf{z}):=\sigma(\Re\mathbf{z})+{\rm i}\Im\mathbf{z}\cdot\nabla\sigma(\Re\mathbf{z})\quad\mbox{and}\quad\tau(\mathbf{z}):=\tau(\Re\mathbf{z})+{\rm i}\Im\mathbf{z}\cdot\nabla\tau(\Re\mathbf{z}).
\end{equation}
The treatment of piecewise analytic functions is given in Remark \ref{remark:piecewiseanalytic} below.
We will prove the existence of a holomorphic extension of the domain-to-solution map $\mathcal{T}\to \mathcal{H}$ defined by $T\mapsto \widehat{u}_T$ to the domain
$$\mathcal{T}_\delta:=\{T\in W^{1,\infty}(B_{\rm o},\mathbb{C}^d)\,:\,\|T-T_0\|_{W^{1,\infty}(B_{\rm o},\mathbb{C}^d)}<\delta\mbox{ for some }T_0\in\mathcal{T}\}$$
for some $\delta>0$. In the next lemma, we extend the definition of $\widehat{u}_T$ for $T\in\mathcal{T}$ to $T\in\mathcal{T}_\delta$.
\begin{lemma}\label{lemma:welldefined:complex}
There exists some $\delta>0$ such that for all $T\in\mathcal{T}_\delta$, the variational problem \eqref{eq:varform:extended:u} with \eqref{eq:extensions:sigmatau} has a unique solution.
\end{lemma}
\begin{proof}
The maps $T''\mapsto J_{T''}\in L^\infty(B_{\rm o},\mathbb{C}^{d\times d})$ and $T''\mapsto \sigma\circ T''\in L^\infty(B_{\rm o},\mathbb{C})$ are continuous in $T''\in W^{1,\infty}(B_{\rm o},\mathbb{C}^d)$ and we have $|J_T|\neq0$ for all $T\in\mathcal{T}$.
Thus, for each $T\in\mathcal{T}$, there exists a $\delta_T>0$ such that for all $T'\in W^{1,\infty}(B_{\rm o},\mathbb{C}^d)$ satisfying $\|T-T'\|_{W^{1,\infty}(B_{\rm o},\mathbb{C}^d)}<2\delta_T$, there holds
\begin{equation}\label{eq:smallperturb:complexval}
\||J_{T'\circ T^{-1}}|(\sigma\circ T'\circ T^{-1})J_{T'\circ T^{-1}}^{-1}J_{T'\circ T^{-1}}^{-\top}-\sigma{\rm Id}\|_{L^\infty(B_{\rm o},\mathbb{C}^{d\times d})}<\tfrac{1}{2}\min\{\sigma(\Bx)\,:\,\Bx\in B_{\rm o}\}.
\end{equation}
We define an open cover $\{B_{\delta_T}(T)\}_{T\in\mathcal{T}}$ of $\mathcal{T}$ by
$$B_{\delta_T}(T):=\{T'\in W^{1,\infty}(B_{\rm o},\mathbb{C}^d)\,:\,\|T-T'\|_{W^{1,\infty}(B_{\rm o},\mathbb{C}^d)}<\delta_T\},\quad T\in\mathcal{T}.$$
By Assumption \ref{ass:T:compact}, $\mathcal{T}$ is compact and hence it has a finite subcover, namely, $\{B_{\delta_{T_i}}(T_i)\}_{i=1}^n$.
Let $\delta := \min\{\delta_{T_i}\,:\,1\le i\le n\}$.
Then, for every $T'\in \mathcal{T}_\delta$, there exists a $T\in\mathcal{T}$ and $i\in\{1,\dots,n\}$ satisfying $\|T'-T\|_{W^{1,\infty}(B_{\rm o},\mathbb{C}^d)}<\delta$ and $\|T-T_i\|_{W^{1,\infty}(B_{\rm o},\mathbb{C}^d)}<\delta_{T_i}$, which implies $\|T'-T_i\|_{W^{1,\infty}(B_{\rm o},\mathbb{C}^d)}<2\delta_{T_i}$, so that the estimate \eqref{eq:smallperturb:complexval} holds for $T=T_i$.
Meanwhile, for all $u\in \mathcal{H}$, we have
\begin{align*}
    a_{T'} (u, u) &= \int_{E}\left[ (\sigma\circ {T'})J_{T'}^{-1}J_{T'}^{-\top}\nabla u\cdot\nabla \overline{u} - k^2(\tau\circ {T'})|u|^2\right]|J_{T'}|\d \mathbf{x}+ (\Lambda u,{u})_{L^2(\p B_{\rm o})}.
\end{align*}
We change the variable $y=T_i x$ in the integral and then apply \eqref{eq:smallperturb:complexval} with $T=T_i$ and $\Re\left[(\Lambda u,{u})_{L^2(\p B_{\rm o})}\right]\ge0$ from \cite[Corollary 3.1]{Chandler:2008:WNEBTIS}.
Then, there exist positive $C$ and $C'$ depending only on $i$ satisfying
$$
\Re [a_{T'}(u, u)]\ge C\|u\circ T_i\|^2_{H^1(E_{T_i})}-C'\|u\circ T_i\|^2_{L^2(E_{T_i})},\quad\forall T'\in B_{\delta_{T_i}}(T_i).
$$
Since $H^1(E_{T_i})$ is compactly embedded in $L^2(E_{T_i})$, $A_{T'}:\mathcal{H}_{T_i}\to\mathcal{H}_{T_i}'$ is an isomorphism for all $T'\in B_{\delta_{T_i}}(T_i)$ by the Fredholm alternative \cite[Theorem 5.4.5]{Nedelec:2001:AEE}.
\end{proof}

We also use the following lemma \cite[Lemma 4.1]{JSZ:2017:EWS}.
\begin{lemma}\label{lemma:differentials:T}
Let $D'$ be a bounded domain. The map $W^{1,\infty}(D',\mathbb{C}^d)\to L^\infty(D',\mathbb{C}^{d\times d})$ defined by $T\mapsto J_T$ is holomorphic.
If $T\in W^{1,\infty}(D',\mathbb{C}^d)$  satisfies $T^{-1}\in W^{1,\infty}(D',\mathbb{C}^d)$, then the map $T\mapsto J_T^{-1}$ is holomorphic at $T$.
\end{lemma}
\begin{theorem}\label{theorem:holomorphy:planewave}
There exist a constant $\delta>0$ and a holomorphic extension $F:\mathcal{T}_\delta\to H^{1}(E,\mathbb{C})$ of the domain-to-solution map $T\mapsto\widehat{u}_T$.
\end{theorem}
\begin{remark}\label{remark:FrechetDerivative}
The proof of Theorem \ref{theorem:holomorphy:planewave} shows that the Fr{\'e}chet derivative $F'(T)(H)$ of $F$ at $T\in\mathcal{T}_\delta$ is given by the unique solution $U_{T}(H)\in H^{1}(E,\mathbb{C})$ of the following problem:
$$
{a}_T(U_{T}(H),w)=g_{T,H}(\widehat{u}_T,w),\quad\forall w\in\mathcal{H},
$$
with the linear form $g_{T,H}(\widehat u_T,w)$ given by
\begin{align}
&	g_{T,H}(\widehat{u}_T,w)\\
    :=&k^2\int_{E}\left[(\tau\circ T)\operatorname{tr}(J_T^*J_H) + (H\cdot\nabla\tau\circ T)|J_ T|\right]\widehat{u}_{T} \overline{w}~\d \mathbf{x}\nonumber\\
&-\int_{E}(\sigma\circ T)\left[\operatorname{tr}(J_T^*J_H)J_T^{-1}J_T^{-\top}-|J_T|J_T^{-1}\left(J_H J_T^{-1}+J_T^{-\top}J_H^{\top}\right)J_T^{-\top}\right]\nabla\widehat{u}_{T}\cdot\nabla \overline{w}~\d \mathbf{x}\nonumber\\
&+\int_{E}(H\cdot\nabla\sigma\circ T)|J_T| J_T^{-1}J_T^{-\top} \nabla \widehat{u}_{T}\cdot\nabla \overline{w}~\d\mathbf{x}.\label{eqn:g-TH}	
\end{align}
\end{remark}
\begin{proof}
Fix $T\in\mathcal{T}_\delta$, $H\in W^{1,\infty}(B_{\rm o},\mathbb{C}^d)$ such that $T+H\in\mathcal{T}_\delta$ for $\delta$ satisfying Lemma \ref{lemma:welldefined:complex}. Then for any   $w\in\mathcal{H}$, using the relation
$$a_{T+H}(\widehat{u}_{T+H},w)=b(w)=a_{T}(\widehat{u}_{T},w),$$
we obtain
\begin{align}
&a_T(\widehat{u}_{T+H}-\widehat{u}_T,w)=a_T(\widehat{u}_{T+H},w)-a_{T+H}(\widehat{u}_{T+H},w)\label{eqn:bilin-diff}\\
=&-k^2\int_{E}\left((\tau\circ T)|J_T| - (\tau\circ (T+H))|J_{T+H}\right)\widehat{u}_{T+H} \overline{w}~\d \mathbf{x}\nonumber\\
& + \int_{E}\left((\sigma\circ T)|J_T| J_T^{-1}J_T^{-\top} - (\sigma\circ (T+H))|J_{T+H}| J_{T+H}^{-1}J_{T+H}^{-\top}\right)\nabla \widehat{u}_{T+H}\cdot\nabla \overline{w}~\d\mathbf{x}.\nonumber
\end{align}
Note that  by Lemma \ref{lemma:differentials:T}, the following identities hold
\begin{align*}	
|J_{T+H}| =& |J_T| + \operatorname{tr}(J_T^*J_H) + R_1(H),\\
|J_{T+H}| J_{T+H}^{-1}J_{T+H}^{-\top} =& |J_T| J_T^{-1}J_T^{-\top} + \operatorname{tr}(J_T^*J_H)J_T^{-1}J_T^{-\top}\\
				&-|J_T| \left(J_T^{-1} J_H J_T^{-1}J_T^{-\top} + J_T^{-1}J_T^{-\top}J_H^{\top}J_T^{-\top}\right) + R_2(H),
\end{align*}
with $\|R_1(H)\|_{L^\infty(E,\mathbb{C})}$ and $\|R_2(H)\|_{L^\infty(E,\mathbb{C}^{d\times d})}$ being of the order $o(\|H\|_{W^{1,\infty}(E,\mathbb{C}^d)})$ as $\|H\|_{W^{1,\infty}(E,\mathbb{C}^d)}\to0$.
Further, for all $T\in\mathcal{T}_\delta$ and $H\in W^{1,\infty}(B_{\rm o},\mathbb{C}^d)$ such that $T+H\in\mathcal{T}_\delta$, there holds the identity
\begin{align}
\sigma\circ(T+H)&=\sigma\circ \Re (T+H) + {\rm i}\Im (T+H)\cdot\nabla\sigma\circ\Re (T+H)\nonumber\\
&=\sigma\circ \Re T + \Re H\cdot\nabla\sigma\circ \Re T + {\rm i}\Im (T+H)\cdot\nabla\sigma\circ\Re T+o(\|H\|_{L^\infty(B_{\rm o},\mathbb{C}^d)})\nonumber\\
&=\sigma\circ T + H\cdot\nabla\sigma\circ T+o(\|H\|_{L^\infty(B_{\rm o},\mathbb{C}^d)}).\nonumber
\end{align}
Similarly, for piecewise linear $\tau$, we have
\begin{equation}\label{eq:tau:Frechetderivative}
    \|\tau\circ(T+H)-\tau\circ T - H\cdot\nabla\tau\circ T\|_{L^\infty(B_{\rm o},\mathbb{C})}=o(\|H\|_{L^\infty(B_{\rm o},\mathbb{C}^d)}).
\end{equation}
Substituting the preceding identities into the identity \eqref{eqn:bilin-diff} gives
$$a_T(u_{T+H}-u_T,w)=g_{T,H}(\hat{u}_{T},w)+R_3(w,H),\quad\forall w\in\mathcal{H},$$
where ${\|R_3(w,H)\|_{L^\infty(E,\mathbb{C})}}=o(\|w\|_{\mathcal{H}}{\|H\|_{W^{1,\infty}(E,\mathbb{C}^d)}})$ as $\|H\|_{W^{1,\infty}(E,\mathbb{C}^d)}\to0$.
Thus, for $U_T(H)$ defined in Remark \ref{remark:FrechetDerivative}, we have
$$a_T(u_{T+H}-u_T-U_T(H),w)=R_3(w,H),\quad\forall w\in\mathcal{H}.$$
Since $a_T$ induces an isomorphism, we conclude that $U_T(H)$ is the complex Fr{\'e}chet derivative of the map $T \mapsto u_T$ for all $T\in\mathcal{T}_\delta$ as in Definition \ref{definition:holomorphic}.
\end{proof}

\begin{remark}\label{remark:piecewiseanalytic}
If $\sigma$ and $\tau$ are piecewise analytic, we may use an extension to complex variables different from \eqref{eq:extensions:sigmatau}, by taking higher order Taylor expansions.
Let $\sigma$ be piecewise analytic in the following sense: There exists a finite set $\{\omega_j\}_{j}$ of open, pairwise disjoint subsets of $\mathbb{R}^d$ with Lipschitz boundaries such that $\cup_j\overline{\omega_j}=\mathbb{R}^d$ and, for every $\omega\in\{\omega_j\}_{j}$, there exists an analytic function $\widetilde{\sigma}:\mathbb{R}^d\to\mathbb{R}$ such that $\widetilde{\sigma}|_{\omega}={\sigma}|_{\omega}$.
Let $\mathbf{n}=(n_1,\dots,n_d)\in \mathbb{N}_0^d$, $\p^{\mathbf{n}}=\p_{x_1}^{n_1}\cdots\p_{x_d}^{n_d}$, $\mathbf{n}!=\prod_{i=1}^d 
(n_i!)$ and $\mathbf{z}^{\mathbf{n}}=\prod_{i=1}^d z_i^{n_i}$.
Then by assumption, for each $\omega\in\{\omega_j\}_{j}$, $\widetilde{\sigma}$ admits the following Taylor expansion:
$$\widetilde{\sigma}(\Bx)=\sigma(\Bx_0)+\sum_{|\mathbf{n}|=1}^\infty \frac{1}{\mathbf{n}!}(\p^{\mathbf{n}}\sigma)(\Bx_0)(\Bx-\Bx_0)^{\mathbf{n}},\quad\forall \Bx\in \mathbb{R}^d,\,\Bx_0\in \omega.$$
Then we define an extension of $\sigma$ to $\mathbb{C}^d$ almost everywhere by
\begin{equation}\label{eq:sigma:extension}
\sigma(\Bz)=\sigma(\Re\Bz)+\sum_{|\mathbf{n}|=1}^\infty \frac{1}{\mathbf{n}!}(\p^{\mathbf{n}}\sigma)(\Re\Bz)({\rm i}\Im\Bz)^{\mathbf{n}},\quad\forall \Bz\in \mathbb{C}^d\backslash({\rm i}\mathbb{R}^d+\cup_j\p \omega_j).
\end{equation}
By following the proof of Theorem \ref{theorem:holomorphy:planewave}, we prove that $F$ is holomorphic: for all $T\in \mathcal{T}_\delta$ and $H\in W^{1,\infty}(B_{\rm o},\mathbb{C})$, the definition \eqref{eq:sigma:extension} of the extended function $\sigma$ gives
\begin{align*}
    &\sigma(T+H)-\sigma(T)=\sum_{|\mathbf{n}|=0}^\infty \frac{1}{\mathbf{n}!}\left[(\p^{\mathbf{n}}\sigma\circ\Re(T+H))({\rm i}\Im(T+H))^{\mathbf{n}}-(\p^{\mathbf{n}}\sigma\circ\Re T) ({\rm i}\Im T)^{\mathbf{n}}\right]\\
    =&\sum_{|\mathbf{n}|=0}^\infty \frac{1}{\mathbf{n}!}\left[\Re H\cdot(\nabla \p^{\mathbf{n}}\sigma\circ\Re T)({\rm i}\Im T)^{\mathbf{n}} + (\p^\mathbf{n}\sigma\circ\Re T) \nabla(\Bz^{\mathbf{n}})|_{\Bz={\rm i}\Im T}\cdot({\rm i}\Im H) \right]+ o(\|H\|_{L^\infty(B_{\rm o},\mathbb{C}^d)})\\
    =&\Re H\cdot\nabla\sigma\circ T+{\rm i}\Im H\cdot\nabla\sigma\circ T+o(\|H\|_{L^\infty(B_{\rm o},\mathbb{C}^d)})
\end{align*}
as $\|H\|_{L^{\infty}(B_{\rm o},\mathbb{C}^d)}\to0$.
Similarly, we can extend the relation \eqref{eq:tau:Frechetderivative} for $\tau$ as \eqref{eq:sigma:extension}.
Thus, the proof of Theorem \ref{theorem:holomorphy:planewave} is valid under the weaker assumption that $\sigma$ and $\tau$ are piecewise analytic, and the complex Fr{\'e}chet derivative of $F$ is of the form as in Remark \ref{remark:FrechetDerivative}.
\end{remark}

By Assumption \ref{ass:T}, for all $T\in\mathcal{T}$, there holds $\widehat{u}_T|_{\p B_{\rm o}}=u_T|_{\p B_{\rm o}}$.
Finally, fix any $m\in\mathbb{N}$ and define the near-to-far field operator $G:H^{1/2}(\p B_{\rm o})\to C^m(\mathbb{S}^{d-1})$ by $u|_{\p B_{\rm o}}\mapsto u^{\infty}(\cdot,\mathbf{d})$, where $u\in H_{\rm loc}^{1}(\mathbb{R}^d\backslash\overline{B_{\rm o}})$ is the solution to $\Delta u+ k^2u=0$ in $\mathbb{R}^d\backslash\overline{B_{\rm o}}$ subject to the Sommerfeld radiation condition \eqref{SmmerfeldRadiationCondition} for $u^{\rm s}=u-u^{\rm i}$ and $u^{\infty}(\cdot,\mathbf{d}) \in C^\infty(\mathbb{S}^{d-1})$ is the far-field pattern in \eqref{FarFieldPattern} \cite[Lemma 2.5]{Chandler:2008:WNEBTIS}.
By Theorem \ref{theorem:holomorphy:planewave}, since the trace operator $\operatorname{tr}_{\p B_{\rm o}}:H^1(E,\mathbb{C})\to 
H^{1/2}(\p B_{\rm o})$ and $G:H^{1/2}(\p B_{\rm o})\to C^{m}(\p B_{\rm o})$ are $T$-independent continuous operators, the composition 
$G\circ \operatorname{tr}_{\p B_{\rm o}}\circ F$
is a holomorphic extension of $T\mapsto u_T^\infty$.
In view of Theorem \ref{theorem:holomorphy:pointsource}, the assertion holds also for the map $T\mapsto v_T^\infty$.
\begin{corollary}\label{cor:JtoH:holomorphy}
Let $m\in\mathbb{N}$. There exist a $\delta>0$ and one holomorphic map $\mathcal{T}_\delta\to C^m(\mathbb{S}^{d-1})$ whose restriction 
to $\mathcal{T}$ is $T\mapsto u_T^\infty$.
\end{corollary}

\subsection{Parametric shape holomorphy}
Now we establish the holomorphy of the forward maps with respect to the shape parameters.
Fix a convex $C^1$ open set $B_3$ satisfying $\overline{B_3}\subset B_{\rm i}\backslash(\overline{B_1}\cup\overline{B_2})$, fix any $r>0$ and restrict the discussions to domains of the class
\begin{align*}
	\mathcal{C}_r:=\Big\{ &\Om \,:\, \Om\mbox{ is star shaped, }\overline{\Om}\subset B_3\mbox{ and }\operatorname{dist}(\Om,\p B_3)>r\Big\}.
\end{align*}
Suppose that $\widehat{D}$ is a ball  centered at $\mathbf{a}_0$ of radius $r_0$ in $\mathcal{C}_r$.
Then there is a natural one-to-one correspondence between $\p\widehat{D}$ and $\p\Om(\mathbf{a},\mathbf{b})$ defined by
$$\p\widehat{D}\to \p\Om(\mathbf{a},\mathbf{b}),\quad \mathbf{x}\mapsto\mathbf{a}+\frac{|\mathbf{x}-\mathbf{a}_0|}{r_0}\exp\rho_b\left(\frac{\mathbf{x}-\mathbf{a}_0}{|\mathbf{x}-\mathbf{a}_0|}\right), \quad\forall \mathbf{x}\in\p\widehat{D},$$
where $(\mathbf{a},\mathbf{b})$ are the shape parameters introduced in Section \ref{sect:starshape}.
By the convexity of $B_3$, we have the following natural extension of the boundary correspondence to  transformations in $\mathcal{T}$:
Given a ball $\widehat{D}\in\mathcal{C}_r$  centered at $\mathbf{a}_0$ of radius $r_0$, for every $\Om=\Om(\mathbf{a},\mathbf{b})\in\mathcal{C}_r$, define $T(\cdot;\mathbf{a},\mathbf{b})\in\mathcal{T}$ by $T(\cdot;\mathbf{a},\mathbf{b})={\rm id}$ in $B_{\rm o}\backslash{B_3}$, $T(\mathbf{a}_0;\mathbf{a},\mathbf{b})=\mathbf{a}$ and
\begin{align*}
T(\mathbf{x};\mathbf{a},\mathbf{b})=&\begin{cases}
		\displaystyle	 \mathbf{a}+\frac{|\mathbf{x}-\mathbf{a}_0|}{r_0}{\exp\rho_{\mathbf{b}}(\mathbf{c}_{\mathbf{a}}(\Bx))}\mathbf{c}_{\mathbf{a}}(\Bx),&\mbox{if }\mathbf{x}\in\widehat{D}\backslash\{\mathbf{a}_0\},\\[1mm]
	\displaystyle		\mathbf{a}
			 +\frac{|\mathbf{c}(\Bx)-\mathbf{a}_0|-|\mathbf{x}-\mathbf{a}_0|}{|\mathbf{c}(\Bx)-\mathbf{a}_0|-r_0}\exp\rho_{\mathbf{b}}(\mathbf{c}_{\mathbf{a}} (\Bx))\mathbf{c}_{\mathbf{a}}(\Bx)\\
			\quad+\frac{|\mathbf{x}-\mathbf{a}_0|-r_0}{|\mathbf{c}(\Bx)-\mathbf{a}_0|-r_0}\left(\mathbf{c}(\Bx)-\mathbf{a}\right),&\mbox{if }\mathbf{x}\in B_3\backslash\widehat{D},
		\end{cases}
	\end{align*}
	with $\{\mathbf{c}(\Bx)\}=\{\mathbf{a}_0+y(\mathbf{x}-\mathbf{a}_0)\,:\,y>0\}\cap\p B_3$ and $\mathbf{c}_{\mathbf{a}}(\Bx)=\frac{\mathbf{c}(\Bx)-\mathbf{a}}{|\mathbf{c}(\Bx)-\mathbf{a}|}$.
Note that there hold $T(\cdot;\mathbf{a},\mathbf{b})\in\mathcal{T}$, $T(\widehat{D};\mathbf{a},\mathbf{b})=\Om(\mathbf{a},\mathbf{b})$ and $T(B_3;\mathbf{a},\mathbf{b})=B_3$. Next we recall the concept of the $(\boldsymbol{\beta},p,\varepsilon)$-holomorphy.

\begin{definition}[{\cite[Definition 15.3.3]{OSZ:2022:DLHD}}]
Let $\mathcal{X}$ be a complex Banach space equipped with the norm $\|\cdot\|_{\mathcal{X}}$. For $\varepsilon >0$ and $\boldsymbol{\beta}\in\ell^p(\mathbb{N})$ with some $p\in (0,1)$, the map $u:[-1,1]^{\mathbb{N}}\to \mathcal{X}$
is said to be $(\boldsymbol{\beta},p,\varepsilon)$-holomorphic if and only if the following three conditions hold:
\begin{itemize}
	\item[\rm(i)] The map $u:[-1,1]^{\mathbb{N}}\to \mathcal{X}$ is continuous.
	\item[\rm(ii)] There exists a sequence $\boldsymbol{\beta} := (\beta_j)_{j\geq 1} \in \ell^p(\mathbb{N})$ of positive numbers
	such that for any sequence $\boldsymbol{\rho} := (\rho_j)_{j\geq 1}\subset (1,\infty)^{\mathbb{N}}$ that is $(\boldsymbol{\beta} ,\varepsilon)$-admissible, i.e., satisfying
	$$\sum_{j\geq 1}(\rho_j - 1)\beta_j \leq \varepsilon,$$
	the map $u:[-1,1]^{\mathbb{N}}\to \mathcal{X}$ admits a complex extension $u:\mathcal{E}_{\boldsymbol{\rho}}\to \mathcal{X}$ that is holomorphic with respect to each component $z_j\in\mathbb{C}$ of $(z_j)_{j\in\mathbb{N}}$ in the set $\mathcal{E}_{\boldsymbol{\rho}}$ defined as		$$\mathcal{E}_{\boldsymbol{\rho}}:=\bigotimes_{j\geq1}\mathcal{E}_{\rho_j},\quad\mbox{with } \mathcal{E}_{\rho_j}=\{(z+z^{-1})/2\,:\,1\le |z|\le \rho_j\}.$$
	\item[\rm(iii)] Any extension $u$ in (ii) is uniformly bounded. In other words, there exists a constant $M>0$ independent of $\boldsymbol{\rho}$ satisfying $\sup\{\|u(\boldsymbol{z})\|_{\mathcal{X}}\,:\,\boldsymbol{z}\in\mathcal{E}_{\boldsymbol{\rho}}\}\le M<\infty$.
\end{itemize}
\end{definition}

For all $m\in\mathbb{Z}$ and $|m|\leq \ell$, let  \begin{equation*}
\widetilde{w}_m:=(1+|m|)^{-1}\quad \mbox{and}\quad \widetilde{w}_{\ell,m}:= (\ell+1)^{-1/2}{\ell+|m|\choose 2|m|}^{-1/2}\frac{(|m|+1)^{1/4}}{\ell^2 + |m|^2-\ell|m|+1}.
\end{equation*}

\begin{lemma}\label{lemma:paratransform:holomorphy}
	Let $w>0$, $\{w_m\}_{m\in\mathbb{Z}}$ and $\{w_{\ell,m}\}_{m\in\mathbb{Z},\,\ell\ge|m|}$ satisfy, for some $C>0$ and $q>0$ such that
    $$0\le w_m\le C(|m|+1)^{-1-q}\widetilde{w}_m\quad\mbox{and}\quad 0\le w_{\ell,m}\le C(\ell+1)^{-1-q}\widetilde{w}_{\ell,m},\quad \forall m\in\mathbb{Z},\ell\ge|m|.$$
	Suppose also that for all $(\widehat{\mathbf{a}},\widehat{\mathbf{b}})\in[-1,1]^{\mathbb{N}}$, the pair $(\mathbf{a},\mathbf{b})$ satisfies $\Om(\mathbf{a},\mathbf{b})\in \mathcal{C}_r$, where $\mathbf{a}=\mathbf{a}_0+w\widehat{\mathbf{a}}$, $b_m=w_{m}\widehat{b}_m$ in 2D and $b_{\ell,m}=w_{\ell,m}\widehat{b}_{\ell,m}$ in 3D for all $m\in\mathbb{Z}$ and $\ell\ge|m|$.
	Let $\delta>0$ be as in Corollary \ref{cor:JtoH:holomorphy}.
The map $[-1,1]^{\mathbb{N}}\ni(\widehat{\mathbf{a}},\widehat{\mathbf{b}})\mapsto T(\cdot;\mathbf{a},\mathbf{b})\in\mathcal{T}_\delta$ is $(\boldsymbol{\beta},p,\varepsilon)$-holomorphic for some $\varepsilon>0$, with
\begin{equation}\label{eq:beta:p}\begin{cases}
	\boldsymbol{\beta}=\{1\}^2\times(\widetilde{w}_{m}^{-1} w_{m})_{m\in\mathbb{Z}}\quad\mbox{and all}~ p\in((q+1)^{-1},1)&\mbox{in }2D,\\
	\boldsymbol{\beta}=\{1\}^3\times(\widetilde{w}_{\ell,m}^{-1} w_{\ell,m})_{\ell\ge0,|m|\le\ell}~\mbox{and all}~ p\in((q+1)^{-1/2},1)&\mbox{in }3D.
\end{cases}
\end{equation}
\end{lemma}
\begin{proof}
Since $\Om(\mathbf{a},\mathbf{b})\in\mathcal{C}_r$, we have $|\mathbf{c}(\Bx)-\mathbf{a}|\ge r$ for all $\Bx\neq\mathbf{a}_0$.
Hence, there exists a holomorphic extension of $(\widehat{\mathbf{a}},\widehat{\mathbf{b}})\mapsto T(\cdot;\mathbf{a},\mathbf{b})$ in the variable $\mathbf{a}$ to an open subset of $\mathbb{C}^d$.
The holomorphic extension in the variable $\mathbf{b}$ is obvious in the case $N<\infty$ because the dependence of $T(\cdot;\mathbf{a},\mathbf{b})$ on $\mathbf{b}$ is a composition of linear combinations and exponential. In the case $N=\infty$, we will show that there exists an $\varepsilon>0$ such that $\|\mathbf{c}_{\mathbf{a}}\|_{W^{1,\infty}(B_{3},\mathbb{C})}$, $\|\rho_{\mathbf{b}}(\mathbf{c}_{\mathbf{a}})\|_{W^{1,\infty}(B_{3},\mathbb{C})}$, and thus, $\|T(\cdot;\mathbf{a},\mathbf{b})\|_{W^{1,\infty}(B_{\rm o},\mathbb{C}^d)}$ are bounded uniformly on all $(\boldsymbol{\beta},\varepsilon)$-admissible sequences $\boldsymbol{\rho}$.
For $(\widehat{\mathbf{a}}, \widehat{\mathbf{b}})\in\mathcal{E}_{\boldsymbol{\rho}}$, there holds
\begin{align*}
	\|\rho_{\mathbf{b}}(\mathbf{c}_{\mathbf{a}})\|_{W^{1,\infty}(B_{3})}&\le C\|\rho_{\mathbf{b}}\|_{W^{1,\infty}(\mathbb{S}^{d-1})}\\
&\le \begin{cases}
C\sum_{m\in\mathbb{Z}}\widetilde{w}_{m}^{-1}w_{m}\rho_m, &\mbox{in 2D},\\
			 C\sum_{\ell=0}^\infty\sum_{|m|\le\ell}\widetilde{w}_{\ell,m}^{-1} w_{\ell,m}\rho_{\ell,m},&\mbox{in 3D}.
		\end{cases}
	\end{align*}
    The right-hand sides are uniformly bounded for all $(\boldsymbol{\beta},\varepsilon)$-admissible sequences $\boldsymbol{\rho}$, with $\boldsymbol{\beta}$ in \eqref{eq:beta:p}.
	Redefining indices to be one single indice in $\mathbb{N}$ are $m\mapsto 2|m|+1$ (linear in $|m|$) in 2D and $\ell^2+2|m|+1$ (quadratic in $\ell$) in 3D, which gives $p$ in \eqref{eq:beta:p}.
    We choose $\varepsilon>0$ so small that the range of the map $\mathcal{E}_{\boldsymbol{\rho}}\ni(\widehat{\mathbf{a}},\widehat{\mathbf{b}})\mapsto T(\cdot;\mathbf{a},\mathbf{b})$ is included in $\mathcal{T}_\delta$ for all $(\boldsymbol{\beta},\varepsilon)$-admissible $\boldsymbol{\rho}$, with $\delta$ in Corollary \ref{cor:JtoH:holomorphy}.
\end{proof}

\begin{remark}
The parametric $(\boldsymbol{\beta},p,\varepsilon)$-holomorphy with respect to usual affine shape parameters follows directly from the simple assumption on the abstract basis for transformation \cite[Proposition 5.1]{JSZ:2017:EWS}.
We instead choose non-affine shape parameters defined by the extensions of the natural Fourier-type parametrization for star-shaped domains, which requires a more delicate analysis of the $(\boldsymbol{\beta},p,\varepsilon)$-holomorphy.
Lemma \ref{lemma:paratransform:holomorphy} shows that the decay rate $q$ of the shape parameters determines $p$, which in turn determines the approximation rate by neural networks, cf. Remark \ref{remark:convergencrate}.
\end{remark}

\subsection{Expression rates}
\label{section:exp:rate}
Now we derive the expression rate of fully connected feedforward neural networks (FNNs) with the rectified linear unit (ReLU) activation function for approximating the forward maps.
The forward maps are the parametric shape-to-solution maps $U,V:[-1,1]^{\mathbb{N}}\to \mathcal{X}_m$, with $m\in\mathbb{N}$ and $\mathcal{X}_m=C^m(\mathbb{S}^{d-1})$, defined, respectively, by
\begin{align}\label{def:forwardmap} U(\widehat{\mathbf{a}},\widehat{\mathbf{b}})=u^{\infty}[\Om(\mathbf{a},\mathbf{b})\cup B](\cdot,\mathbf{d})\quad \mbox{and}\quad V(\widehat{\mathbf{a}},\widehat{\mathbf{b}})=v^{\infty}[\Om(\mathbf{a},\mathbf{b})\cup B](\cdot,\mathbf{z}).
\end{align}
Note that for all $(\widehat{\mathbf{a}},\widehat{\mathbf{b}})\in[-1,1]^{\mathbb{N}}$, we have $\Om(\mathbf{a},\mathbf{b})\in\mathcal{A}_{N}(B_1,B_2)$ under the assumptions in Lemma \ref{lemma:paratransform:holomorphy}.

\begin{lemma}\label{lem:forward:reg}
$U$ and $V$ are $(\boldsymbol{\beta},p,\varepsilon)$-holomorphic.
\end{lemma}
\begin{proof}
From Lemma \ref{lemma:paratransform:holomorphy}, the map $[-1,1]^{\mathbb{N}}\to\mathcal{T}_\delta$ defined by $(\widehat{\mathbf{a}},\widehat{\mathbf{b}})\mapsto T(\cdot;{\mathbf{a}},{\mathbf{b}})$ is $(\boldsymbol{\beta},p,\varepsilon)$-holomorphic.
From Corollary \ref{cor:JtoH:holomorphy}, the maps $\mathcal{T}_\delta\to C^m(\mathbb{S}^{d-1})$ defined by $T\mapsto u_T^\infty$ and $T\mapsto v_T^\infty$ are both holomorphic.
Thus $U$ and $V$ are $(\boldsymbol{\beta},p,\varepsilon)$-holomorphic.
\end{proof}

The next result gives the expression rate for the DNN approximation of the forward map $U$. We can also prove the rate for the forward map $V$.
\begin{theorem}
\label{theorem:nn:approx}
Fix any $\mathbf{x}_0\in\mathbb{S}^{d-1}$.
There exist a constant $C>0$ and a sequence $\{U_n\}_{n=1}^\infty$ of ReLU FNNs such that for every $n>2$, $U_n$ has input $(\eta_k)_{k=1}^n$ in $[-1,1]^n$ and output in $\mathbb{R}^2$, and satisfies
$\operatorname{size}(U_n) \le C(1+n\log n \log\log n)$, $
\operatorname{depth}(U_n) \le C(1+\log n\log\log n)$,
and the uniform error bound
$$\sup_{\boldsymbol{\eta}\in [-1,1]^{\mathbb{N}}} |U(\boldsymbol{\eta})(\mathbf{x}_0)-U_n(\eta_1,\dots,\eta_n)|\le C n^{1-1/p},$$
where $p$ satisfies \eqref{eq:beta:p}.
Next, fix any $m\in\mathbb{N}$.
There exist a constant $C'>0$ and a sequence $\{\widetilde{U}_n\}_{n=1}^\infty$ of ReLU FNNs such that for every $n>2$, $\widetilde{U}_n$ has input in $[-1,1]^n\times\mathbb{S}^{d-1}$, output in $\mathbb{R}^2$, and satisfies
$\operatorname{size}(\widetilde{U}_n) \le C'(1+n\log n \log\log n)$, $
\operatorname{depth}(\widetilde{U}_n) \le C'(1+\log n\log\log n)$,
and the uniform error bound
$$\sup_{\boldsymbol{\eta}\in [-1,1]^{\mathbb{N}}} \|U(\boldsymbol{\eta})(\cdot)-\widetilde{U}_n((\eta_k)_{k=1}^n,\cdot)\|_{W^{1,\infty}(\mathbb{S}^{d-1})}\le C' n^{1-1/p}.$$
\end{theorem}
\begin{proof}
The proof is mainly based on \cite[Theorem 15.4.9]{OSZ:2022:DLHD} for ${U}_n$ and \cite[Theorem 15.5.2]{OSZ:2022:DLHD} for $\widetilde{U}_n$.
First, we approximate $U$ by a linear combination of Legendre polynomials.
Let
$$\mathcal{F}=\{\boldsymbol{\nu}\in\mathbb{N}_0^{\mathbb{N}}\,:\,\|\boldsymbol{\nu}\|_{1}<\infty\}.
$$
Also, for $\boldsymbol{\nu}\in\mathcal{F}$, $\boldsymbol{\eta}\in[-1,1]^{\mathbb{N}}$ and $j\in\mathbb{N}$, let
\begin{align*}
&\boldsymbol{\nu}=(\boldsymbol{\nu}_{\le j},\boldsymbol{\nu}_{> j})\quad\mbox{with}\quad \boldsymbol{\nu}_{\le j} = (\nu_1,\dots,\nu_j)\quad\mbox{and}\quad \boldsymbol{\nu}_{> j} = (\nu_{j+1},\nu_{j+2},\dots),\\
&\boldsymbol{\eta}=(\boldsymbol{\eta}_{\le j},\boldsymbol{\eta}_{> j})\quad\mbox{with}\quad \boldsymbol{\eta}_{\le j} = (\eta_1,\dots,\eta_j)\quad\mbox{and}\quad \boldsymbol{\eta}_{> j} = (\eta_{j+1},\eta_{j+2},\dots).
\end{align*}
Let $P_\ell^0$, with $\ell\in\mathbb{N}\cup\{0\}$, be the Legendre polynomial of degree $\ell$. We use the following notation for multi indices:
$$P_{\boldsymbol{\nu}}(\boldsymbol{\eta})=\textstyle\prod_{i=1}^{\dim\boldsymbol{\nu}} P_{\nu_i}^0(\eta_i)\quad \mbox{and}\quad  \boldsymbol{\eta}^{\boldsymbol{\nu}}=\textstyle\prod_{i=1}^{\dim\boldsymbol{\nu}} \eta_i^{\nu_i}.
$$
From \cite[Theorem 15.3.7]{OSZ:2022:DLHD}, for $p$ in \eqref{eq:beta:p}, there exists a positive integer $J$ such that there hold:
\begin{itemize}
\item[\rm(i)] For each $\boldsymbol{\nu}\in\mathcal{F}$,
$$c_{\boldsymbol{\nu}}:=\int_{[-1,1]^J}P_{\boldsymbol{\nu}_{\le J}}(\boldsymbol{\eta}_{\le J})\frac{\p_{\boldsymbol{\eta}_{>J}}^{\boldsymbol{\nu}_{> J}} U(\boldsymbol{\eta}_{\le J},\mathbf{0})}{\boldsymbol{\nu}_{> J}!}\,\d\boldsymbol{\eta}_{\le J}\in C^m(\mathbb{S}^{d-1})
$$
is well-defined and satisfies
$$(\|P_{\boldsymbol{\nu}_{\le J}}\|_{L^\infty([-1,1]^J)}\|c_{\boldsymbol{\nu}}\|_{C^m(\mathbb{S}^{d-1})})_{\boldsymbol{\nu}\in\mathcal{F}}\in\ell^p(\mathcal{F});
$$
\item[\rm(ii)] The series expression
$U(\boldsymbol{\eta})=\sum_{\boldsymbol{\nu}\in\mathcal{F}}c_{\boldsymbol{\nu}}P_{\boldsymbol{\nu}_{\le J}}(\boldsymbol{\eta}_{\le J})\boldsymbol{\eta}_{> J}^{\boldsymbol{\nu}_{> J}}$ for $ \boldsymbol{\eta}\in [-1,1]^{\mathbb{N}}$
converges absolutely and uniformly in the norm of $C^m(\mathbb{S}^{d-1})$;
\item[\rm(iii)] There exist constants $C_1,C_2>0$ and a monotonely increasing sequence $\boldsymbol{\delta}=(\delta_i)_{i\in\mathbb{N}}\subset(1,\infty)$ such that $(\delta_i^{-1})_{i\in\mathbb{N}}\in\ell^{p/(1-p)}(\mathbb{N})$, $\delta_i\le C_1 i^{2/p}$ for all $i\in\mathbb{N}$,
$(\boldsymbol{\delta}^{\boldsymbol{\nu}}\|P_{\boldsymbol{\nu}_{\le J}}\|_{L^\infty([-1,1]^J)}\|c_{\boldsymbol{\nu}}\|_{C^m(\mathbb{S}^{d-1})})_{\boldsymbol{\nu}\in\mathcal{F}}\in\ell^1(\mathcal{F})$, and for $\Lambda_\tau:=\{\boldsymbol{\nu}\in\mathcal{F}\,:\,\boldsymbol{\delta}^{-\boldsymbol{\nu}}\ge\tau\}$,
$$\sup_{\boldsymbol{\eta}\in [-1,1]^{\mathbb{N}}}\left\| U(\boldsymbol{\eta}) - \sum_{\boldsymbol{\nu}\in\Lambda_\tau} c_{\boldsymbol{\nu}} P_{\boldsymbol{\nu}_{\le J}}(\boldsymbol{\eta}_{\le J})\boldsymbol{\eta}_{>J}^{\boldsymbol{\nu}_{>J}}
\right\|_{C^m(\mathbb{S}^{d-1})}\le C_2|\Lambda_\tau|^{1-1/p},\quad\forall\tau\in(0,1).$$
\end{itemize}
Then by \cite[Theorem 15.4.9]{OSZ:2022:DLHD}, there exist $c>0$ and a family of ReLU FNNs $\{U_\tau\}_{\tau\in(0,1)}$ with the input variables indexed by ${i\in D_\tau:= \cup_{\boldsymbol{\nu}\in\Lambda_\tau}\operatorname{supp}\boldsymbol{\nu}}$ satisfying
$\operatorname{size}(U_\tau) \le c(1+|\Lambda_\tau|\log |\Lambda_\tau|\cdot \log\log |\Lambda_\tau|)$ and $
\operatorname{depth}(U_\tau) \le c(1+\log |\Lambda_\tau|\cdot\log\log |\Lambda_\tau|)$,
and the uniform error bound
$$\sup_{\boldsymbol{\eta}\in [-1,1]^{\mathbb{N}}} |U(\boldsymbol{\eta})(\mathbf{x}_0)-U_\tau((\eta_i)_{i\in D_\tau})|\le C_p |\Lambda_\tau|^{1-1/p}.$$
By \cite[Proposition 15.3.8]{OSZ:2022:DLHD}, there exist some $\tau\in(0,1)$ such that $D_\tau=\{j\in\mathbb{N}\,:\,1\le j\le n\}$ and $|\Lambda_\tau|\ge|\{e_j\,:\, j\in D_\tau\}|=n$, which completes the proof for $U_n$.
For the proof for $\widetilde{U}_n$, we additionally use the expression rates by FNNs for the functions in $C^m(\mathbb{S}^{d-1})\hookrightarrow W^{m,\infty}(\mathbb{S}^{d-1})$ \cite[Theorem 1]{Yarotsky:2017:EBA} for arbitrarily large $m\in\mathbb{N}$, which implies that the hypothesis for \cite[Theorem 15.5.2]{OSZ:2022:DLHD} is satisfied for arbitrarily large $\gamma>0$.
\end{proof}

\begin{remark}[Convergence rate]
\label{remark:convergencrate}
We briefly comment on the relation between the boundary regularity of the obstacle and the expression rate of $U$.
The argument of Proposition \ref{prop:regularity} indicates that the boundary regularity of the obstacle increases to $C^\infty$ as the decay rate $q$ of the tail of \eqref{eq:define:parameters-2} increases to infinity. The rate of convergence of $U_n$ and $\widetilde{U}_n$ to $U$ in Theorem \ref{theorem:nn:approx} depends on the regularity index $q$ of the target boundary via \eqref{eq:beta:p}, which gives $n^{1-1/p}\approx n^{-q}$ in 2D and $n^{1-1/p}\approx n^{1-\sqrt{1+q}}$ in 3D.
For $q\to\infty$, we can choose $p$ arbitrarily close to $0$, leading to an algebraic convergence rate of arbitrarily high degree.
\end{remark}

\section{Numerical experiments and discussions}
\label{section:numerics}

Now we present numerical experiments in two- and three-dimensions to illustrate the performance of neural network surrogates and the phaseless inverse obstacle scattering in the Bayesian setting.

\subsection{Neural network surrogates of the forward maps}
\label{subsection:NN training}
First we validate the accuracy of neural network surrogates. To generate training data, we employ the boundary element method (BEM). Let $\mathcal{S}_D[\varphi](\mathbf{x})=\int_{\p D}{\Gamma(\mathbf{x}-\mathbf{y})}\varphi(\mathbf{y})\,\d\sigma(\mathbf{y})$, $
\mathcal{D}_D[\varphi](\mathbf{x})=\int_{\partial D}\partial_{\nu_\mathbf{y}}\Gamma(\mathbf{x}-\mathbf{y})\varphi(\mathbf{y})\,\d s_{\mathbf{y}}$ for $ \mathbf{x}\in\mathbb{R}^d\backslash\p D$ and $\mathcal{K}_D[\varphi](\mathbf{x})= \int_{\partial D}\partial_{\nu_{\mathbf{y}}}\Gamma(\mathbf{x}-\mathbf{y})\varphi(\mathbf{y})\d s_{\mathbf{y}}$ for $ \mathbf{x}\in\p D$.
Then we solve the scattering problem using the Brackage-Werner formulation:
\begin{align*}
    u-u^{\inc}=\left({\rm i}k\beta\mathcal{S}_{\Om\cup B} + \mathcal{D}_{\Om\cup B}\right) \mu, \quad\mbox{in }\mathbb{R}^d\backslash\overline{\Om\cup B},
\end{align*}
where the density $\mu\in L^2(\p\Om\cup\p B)$ solves
\begin{align*}
    &\left({\rm i}k\beta\mathcal{S}_{\Omega\cup B} - \tfrac{1}{2}{\rm id}+\mathcal{K}_{\Omega\cup B}\right) \mu = -u^{\inc},\quad\mbox{on }\partial\Omega\cup \partial B.
\end{align*}
It has a unique solution if  $k\beta$ satisfies $\beta\in\mathbb{R}\backslash\{0\}$ \cite[Section 4.3.2]{Dolz:2024:PSH} ($\beta$ is fixed at ${\rm i}k/2$ below). Then, the far-field pattern $u^\infty(\Omega\cup B)$ is given by
$$
u^\infty[\Om\cup B](\hat{\mathbf{x}},\mathbf{d}) =
\frac{1}{4\pi}\int_{\p\Om\cup\p B}\left(-{\rm i}k\beta + {\rm i}k\hBx\cdot\nu(\tBy) \right) e^{-{\rm i}k\hBx\cdot \tBy}\mu(\tBy;\mathbf{d})\d s_{\tBy}.
$$
We discretize the integral equation using piecewise linear functions on the triangulation of the boundary $\p B\cup\p \Om$, and employ GMRES (with the threshold $10^{-3}$ for the residual) to solve the resulting linear system, as in  ``Gypsilab'' (version 0.61), an open toolbox in MATLAB \cite{Alouges:2018:FBS}.

Next we train FNNs $u_{\rm nn}^\infty$ and $v_{\rm nn}^\infty$ that approximate the maps $U$ and $V$, respectively. To this end, we fix a set of nodes $\{\hBx_q\}_{q=1}^{Q}$ on $\mathbb{S}^{d-1}$: for $q=1,\dots,Q$,
\begin{equation}
\label{eq:nodes:on:boundaries}	\hBx_q:=\begin{cases}
		(\cos (2\pi q/Q),\sin(2\pi q/Q)),&\mbox{if }d=2,\\(\sin\theta_q\cos\varphi_q,\sin\theta_q\sin\varphi_q,\cos\theta_q),&\mbox{if }d=3,
	\end{cases}
\end{equation}
with $\theta_q=\cos^{-1}\left(-1+2(q-1)/(Q-1))\right)$ and $\varphi_q=4\pi(q-1)/(1+\sqrt{5})$.
Let $U_Q$ and $V_Q$ be the evaluations at $\{\hBx_q\}_{q=1}^{Q}$ of the maps $U$ and $V$.
We approximate $U_Q$ and $V_Q$ by FNNs with the ReLU activation function, with five hidden layers: The numbers of nodes in the hidden layers are $3N_{\rm in}$, $6N_{\rm in}$, $12N_{\rm in}$, $24N_{\rm in}$ and $12N_{\rm in}+N_{\rm out}/2$ from the input layer to the output layer, with the numbers of nodes of the input and output layers being $N_{\rm in}$ and $N_{\rm out}$, respectively. We employ the standard mean squared error as the loss.  For the training, we employ the Adam optimizer \cite{Kingma2014AdamAM}, using Glorot uniform initialization scheme for the weights and the learning rate $10^{-3}$, with mini-batches of size $128$.

To measure the accuracy of the learned FNNs $u_{\rm nn}^\infty$ and $v^\infty_{\rm nn}$, we use the relative error
$E_{\rm pw}:={\sum_{i}\|(u^\infty-u_{\rm nn}^\infty)[\Omega_i\cup B]\|_{L^2(\mathbb{S}^{d-1})}^2}/{\sum_{i}\|u^\infty
[\Omega_i\cup B]\|_{L^2(\mathbb{S}^{d-1})}^2}$ for the plane wave excitation, and similarly $E_{\rm ps}$ for the point 
source excitation. The reference solutions $u^{\infty}[\Omega_i\cup B]$ and $v^{\infty}[\Omega_i\cup B]$ are computed using the BEM with $200$ and $2000$ nodes for the triangulation in 2D and 3D, respectively. Throughout we fix the wavenumber $k$ at $2$ in 2D and $5.8509$ in 3D.

\begin{figure}[hbt!]
\centering\setlength{\tabcolsep}{0pt}
\begin{tabular}{ccc}	
\includegraphics[height=5cm]{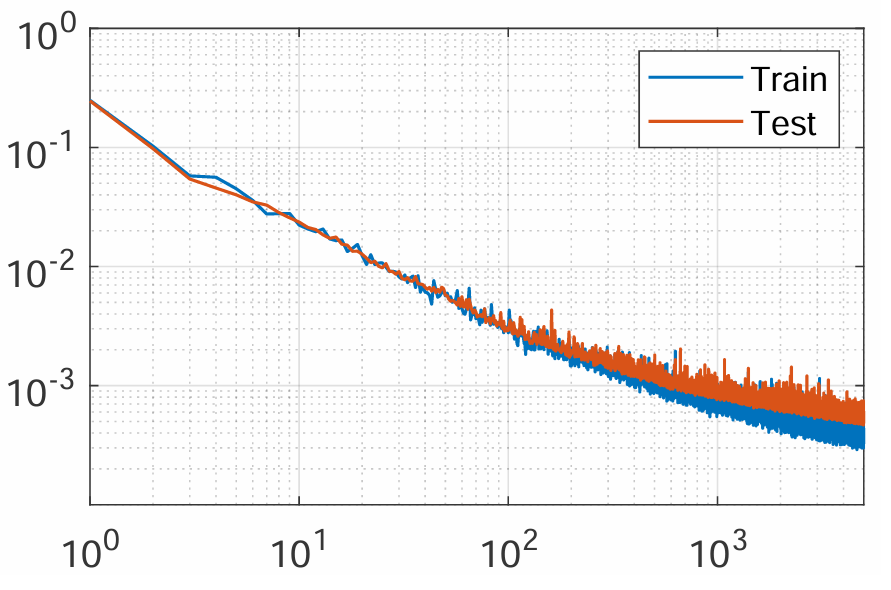}
&\includegraphics[height=5cm]{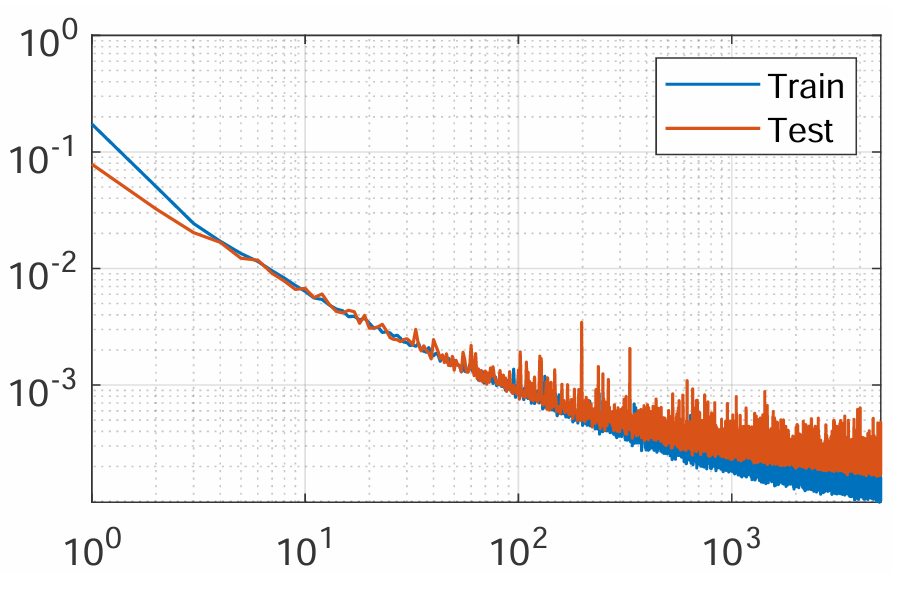}\\
(a) $E_{\rm pw}$ & (b) $E_{\rm ps}$
\end{tabular}
\caption{The evolution of the training and testing errors in 2D for the plane wave excitation (left) and point source excitation (right) in log scale. }
\label{fig:NN:2D:error}
\end{figure}

In the 2D case, we define $(\widehat{\mathbf{a}},\widehat{\mathbf{b}})\mapsto(\mathbf{a},\mathbf{b})$ by $\mathbf{a}=(4,1)+\widehat{\mathbf{a}}$ and $\rho_{\mathbf{b}}(\cos\theta,\sin\theta)=\frac{\widehat{b}_0-1}{4}+\sum_{m=1}^4 \frac{\widehat{b}_m \cos m\theta+\widehat{b}_{-m} \sin m\theta}{|m|+2}$ for $ \theta\in[0,2\pi]$
for all $(\widehat{\mathbf{a}},\widehat{\mathbf{b}})\in[-1,1]^{2+9}$. To avoid excessive deformation, we confine the input to the parameter configuration satisfying
\begin{equation}\label{eq:2D:NN:data:constraint}
    \sum_{m=1}^4 \frac{|\widehat{b}_m| +|\widehat{b}_{-m}|}{m+2}<\log 2.
\end{equation}
The numbers of nodes in the input and output layers are $N_{\rm in}=11$ and $N_{\rm out}=200$, respectively.
We employ the BEM to generate training data for $20000$ different obstacles  by sampling the input $(\widehat{\mathbf{a}},\widehat{\mathbf{b}})$ from the uniform distribution on $[-1,1]^{11}$ satisfying \eqref{eq:2D:NN:data:constraint}, among which $80\%$ and $20\%$ are used for training and testing, respectively. The output layers consist of evaluations at $100$ grid points. The dynamics of training and test errors are in Fig. \ref{fig:NN:2D:error}, which show a steady decrease of the training loss. The accuracy of the trained DNNs is $E_{\rm pw}=4.1\times10^{-4}$ and $E_{\rm ps}=1.0\times10^{-4}$ on the training set and $E_{\rm pw}=5.1\times10^{-4}$ and $E_{\rm ps}=1.9\times10^{-4}$ on the test set, showing decent accuracy of the neural network surrogates.

\begin{figure}[hbt!]
\centering\setlength{\tabcolsep}{0pt}
\begin{tabular}{ccc}
\includegraphics[width=0.32\linewidth]{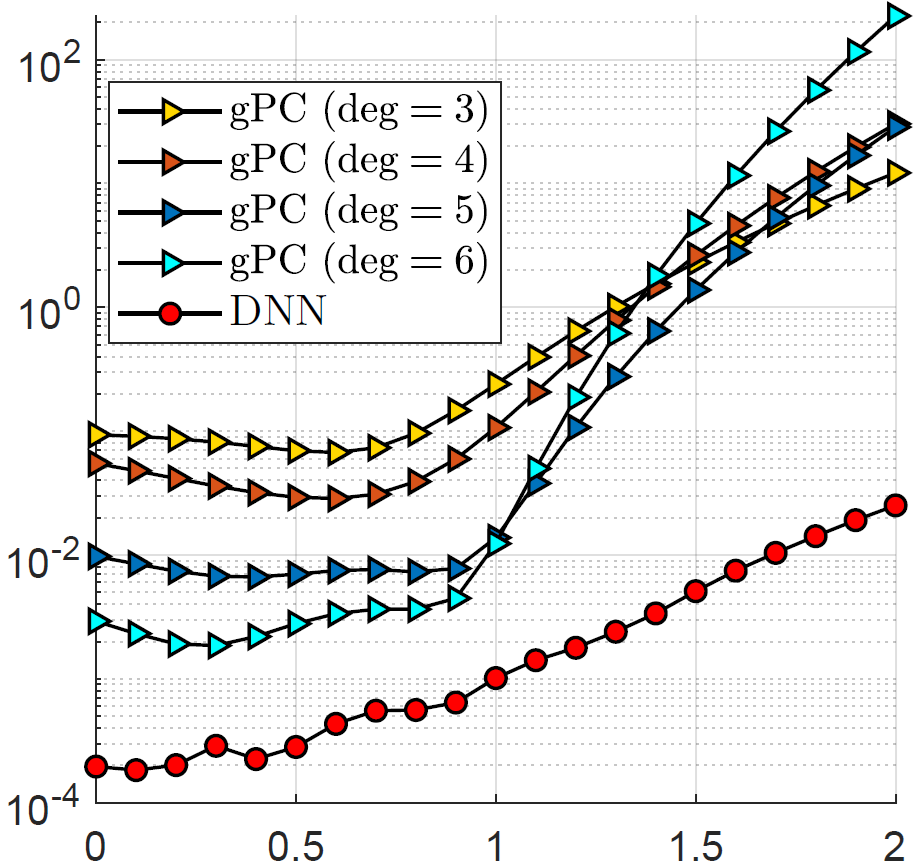}&
\includegraphics[width=0.32\linewidth]{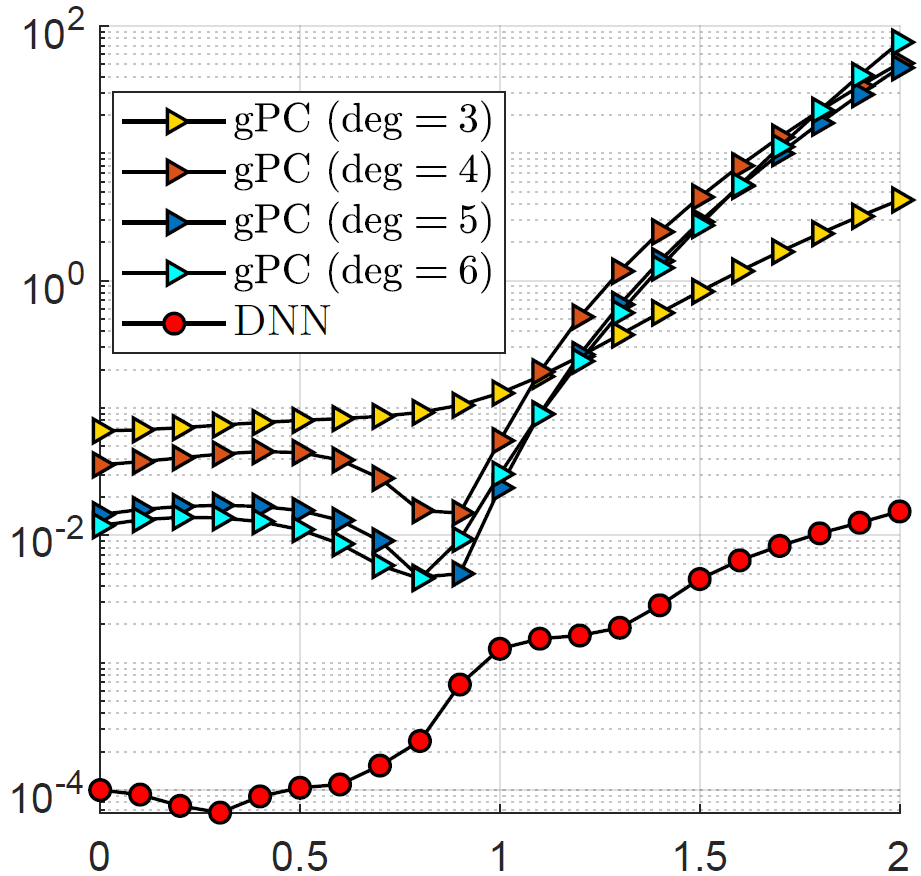}&
\includegraphics[width=0.26\linewidth]{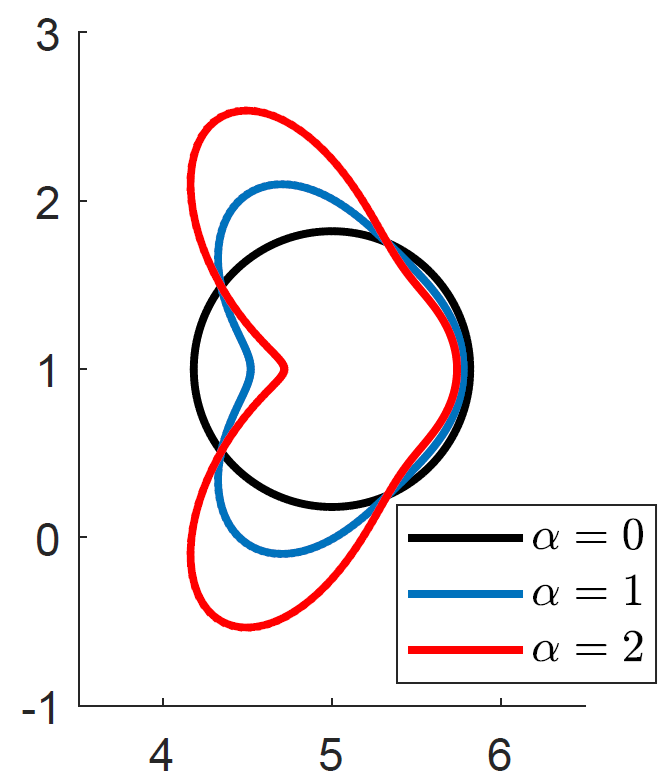}\\
(a) $U$ & (b) $V$ & (c) $\partial\Omega$
\end{tabular}
\caption{The accuracy of the surrogate models for (a) $U$ and (b) $V$, based on DNN and gPC expansion of various degrees, for variations of the kite  Example \ref{ex:2D:inversion}, with ${b}_{m}$ replaced by $\alpha b_{m}$ ($0\le\alpha\le2$) for all $m\in\mathbb{Z}\backslash\{0\}$. The training data were drawn from $\alpha\in[0,1]$.
In (a) and (b), the horizontal axis denotes $\alpha$, and the vertical axis denotes the relative mean squared error (at 100 uniform nodes) on $\mathbb{S}^1$.}\label{fig:shape:evolution:Surrogates}
\end{figure}

For the purpose of comparison, we report also the numerical results by one popular method for constructing surrogates, i.e., the generalized polynomial chaos (gPC) expansion \cite{Marzouk:2009}. For gPC expansion, one can adopt several well established techniques, including Christoffel least squares (CLS) method (see, e.g., \cite{NarayanJakemanZhou2017} or \cite[Section 4.3]{YanZhang:2017}). Note that the standard gPC expansion suffers from the curse of dimensionality: a gPC expansion of degree $d_p$ requires optimizing $M':=\binom{11+d_p}{d_p}$ coefficients for each node, and the CLS method requires the data size $Q'$ larger than $M'$. When implementing the algorithm, we draw $Q'=2M'$ random samples from the CLS sampling density over the parameter domain $[-1,1]^{11}$ for approximating the forward maps $U$ and $V$. The numerical results for the DNN and gPC surrogates are shown in Fig. \ref{fig:shape:evolution:Surrogates}, which includes test cases either inside ($0\le \alpha\le 1$) or outside ($\alpha>1$) of the training distribution, which represent the in-distribution and out-of-distribution cases, respectively.
In both cases, it is observed that the errors are smaller for the neural network surrogates than the gPC surrogate, and the error of the gPC surrogates decreases as the polynomial degree increases on the in-distribution test. Note that the gPC surrogate with degree $6$ requires data of size $Q'=24752$ and the degree of freedom $QM'=1237600$, which are both much larger than that for the DNN surrogate. This observation aligns well with the fact that the expression rate by the neural networks in Section \ref{section:exp:rate} is based on the sparse gPC expansion that can alleviate the curse of dimensionality \cite{Chkifa:2015:BCD}. Note that for both DNN and gPCE surrogates, the test error deteriorates with the increase of the parameter $\alpha$ over the range $(1,2)$ (i.e., the out-of-distribution test data), and high-order gPC surrogates are less robust with respect to the polynomial degree. Thus it is necessary to adapt the surrogates on out-of-distribution test cases, for which pre-training and test-time-adaption seem very promising.

In the 3D case, the input layer has $N_{\rm in}=12$ neurons, since $(\widehat{\mathbf{a}},\widehat{\mathbf{b}})\in[-1,1]^{3+9}$. We use the training data generated with $2000$ nodal points for the triangulation of the boundary.
The BEM dataset contains $4\times3^8$ different obstacles, among which we use $80\%$ ($20995$) for training and the rest ($6249$) for testing.
The accuracy of the trained DNNs is $E_{\rm pw}=0.8\times 10^{-3}$ and $E_{\rm ps}=1.1\times 10^{-3}$ on the training set and $E_{\rm pw}=1.2\times 10^{-3}$ and $E_{\rm ps}=0.9\times 10^{-3}$ for the test set. This again shows the good accuracy of the DNN surrogates for approximating the forward maps.

\subsection{Bayesian reconstruction}
Now we reconstruct the obstacle from the phaseless far-field data using the Bayesian approach \cite{Stuart:2010}. This approach can provide not only point estimators but also associated uncertainties, and thus is very attractive. We employ the following parameterization for the obstacle boundary $\partial\Omega$:
\begin{equation*} \partial\Omega(\mathbf{a},\mathbf{b})=\{\Bx+ \mathbf{a}\,:\,\Bx\in\mathbb{R}^d,\,|\Bx|=\exp{\rho_{\mathbf{b}}(\hBx)}\},\quad \mathbf{a}\in\mathbb{R}^d,\,\mathbf{b}\in\mathbb{R}^{\mathbb{N}}.
\end{equation*}
Then using the far-field values at the nodes $\{\hBx_q\}_{q=1}^Q$ (cf.  \eqref{eq:nodes:on:boundaries}), we form the losses, for $q=1,\dots,Q$,
\begin{align*}
\ds F_{1,q}(\Omega)&:=\big||u_{\rm m}^{\infty}(\hBx_q;\Bd)+v_{\rm m}^{\infty}(\hBx_q;\Bz)| - |{u}^{\infty}[\Omega\cup B](\hBx_q;\Bd)+{v}^{\infty}[\Omega\cup B](\hBx_q;\Bz)|\big|^2,\\
\ds F_{2,q}(\Omega)&:=\big||u_{\rm m}^{\infty}(\hBx_q,\Bd)| - |u^{\infty}[\Omega\cup B](\hBx_q,\Bd)|\big|^2,\\
F_{3,q}(\Omega)&:=\big||v_{\rm m}^{\infty}(\hBx_q,\Bz)| - |v^{\infty}[\Omega\cup B](\hBx_q,\Bz)|\big|^2,
\end{align*}
where the subscript $\mathrm{m}$ indicates the measured phaseless far-field data.
Given the prior distribution $\exp(-\lambda R(\widetilde{\Om}))$ (with the penalty $R$ given below), the posterior distribution $\pi$ is given by
$$
\pi:=\exp\bigg(-\sum_{j=1}^3\frac{1}{2\sigma_j^2}\sum_{q=1}^{Q} F_{j,q}(\widetilde{\Om}) - \lambda R(\widetilde{\Om})\bigg),
$$
where $\sigma_j^2$ are the variances of the additive noise to the phaseless data.
Note that evaluating the posterior distribution $\pi$ requires evaluations of the forward maps $U$ and $V$, which is 
computationally demanding, especially in the 3D case. To reduce the computational expense, we replace $u^\infty$ and $v^\infty$ 
with neural network surrogates $u^\infty_{\rm nn}$ and $v^\infty_{\rm nn}$ and obtain an approximate posterior distribution 
$\tilde \pi$. To explore the posterior distributions $\pi$ and $\tilde \pi$ of  $(\mathbf{a},\mathbf{b})$, we employ Markov 
chain Monte Carlo (MCMC) (specifically, a blockwise Metropolis-Hastings algorithm). The detail of the algorithm is given in 
Appendix \ref{appendix:MCMC}.
The approximate posterior distribution $\widetilde{\pi}$ and the exact one $\pi$ are close to each other with respect to the Kullback-Leibler distance, since the latter can be controlled by the $L^\infty$ norm of the classical BEM solutions and DNN surrogate. See Appendix \ref{appendix:Posteriori approximation} for an error bound in terms of Kullback-Leibler divergence on the posterior approximation based on the surrogate model.

\begin{example}\label{ex:2D:inversion}
The ground truth $(\mathbf{a}^*,\mathbf{b}^*)$ is given by: $\mathbf{a}^*=(5,1)$, and for all $\theta\in[0,2\pi]$,
\begin{align*}					
\rho_{\mathbf{b}_{\rm T}}(\cos\theta,\sin\theta)=\begin{cases}
-0.1-0.2\cos(2\theta)-0.2\sin(2\theta),\quad\mathrm{(ellipse)},\\
-0.1+0.2\cos(3\theta)+0.2\sin(3\theta),\quad\mathrm{(pear)},\\
-0.2+0.04\cos\theta-0.25\cos(2\theta)+0.2\cos(3\theta)-0.04\cos(4\theta),\quad\mathrm{(kite)}.
\end{cases}
\end{align*}
Consider also the case with the rectangular obstacle $\{(x,y)\in\mathbb{R}^2\,:\,|x-5|<0.5\mbox{ and }|y-1|<0.75\}$.
\end{example}

\begin{figure}[hbt!]
\centering\setlength{\tabcolsep}{0pt}
\begin{tabular}{c|c|c}
\toprule
\includegraphics[width=0.32\linewidth]{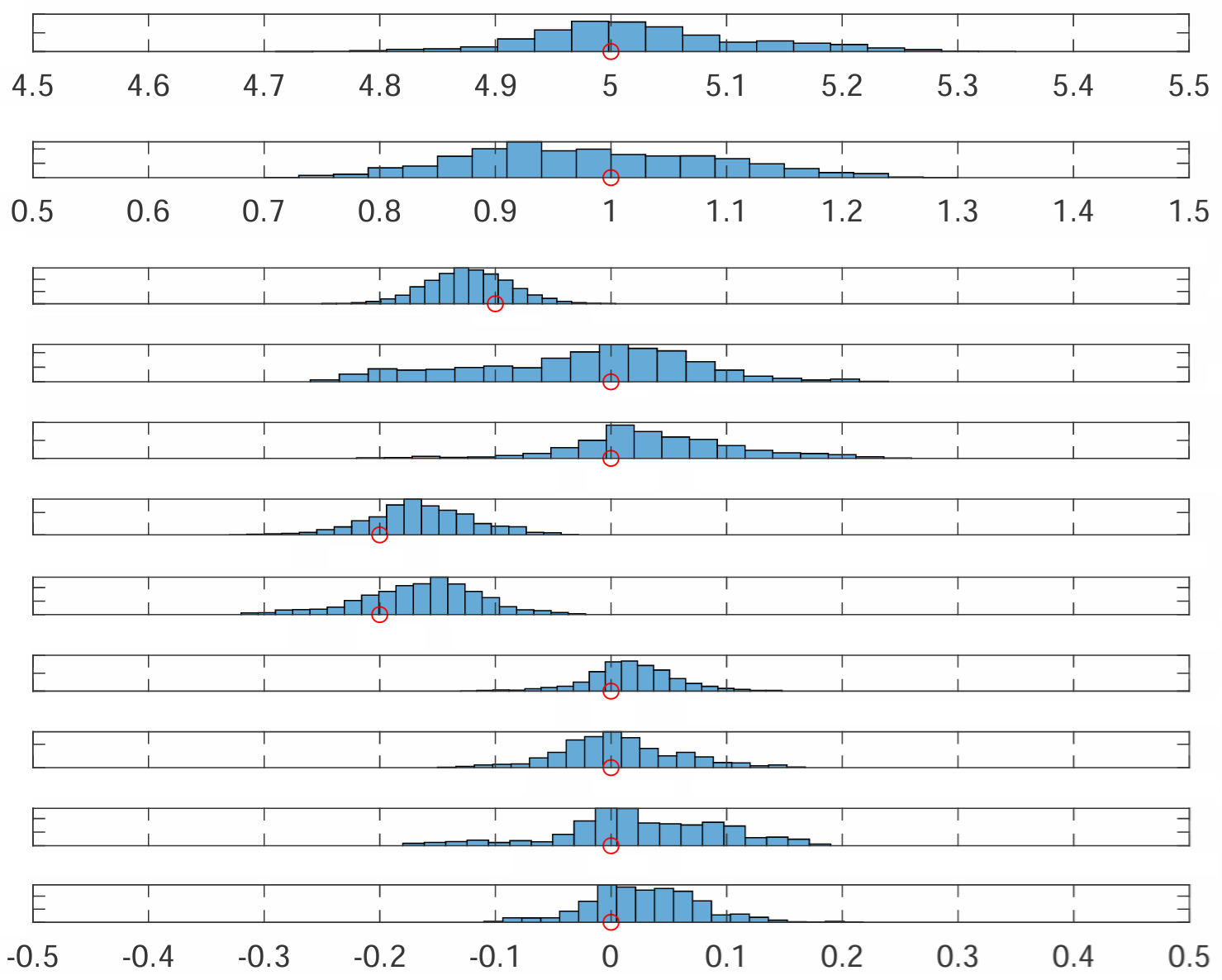}
&\includegraphics[width=0.32\linewidth]{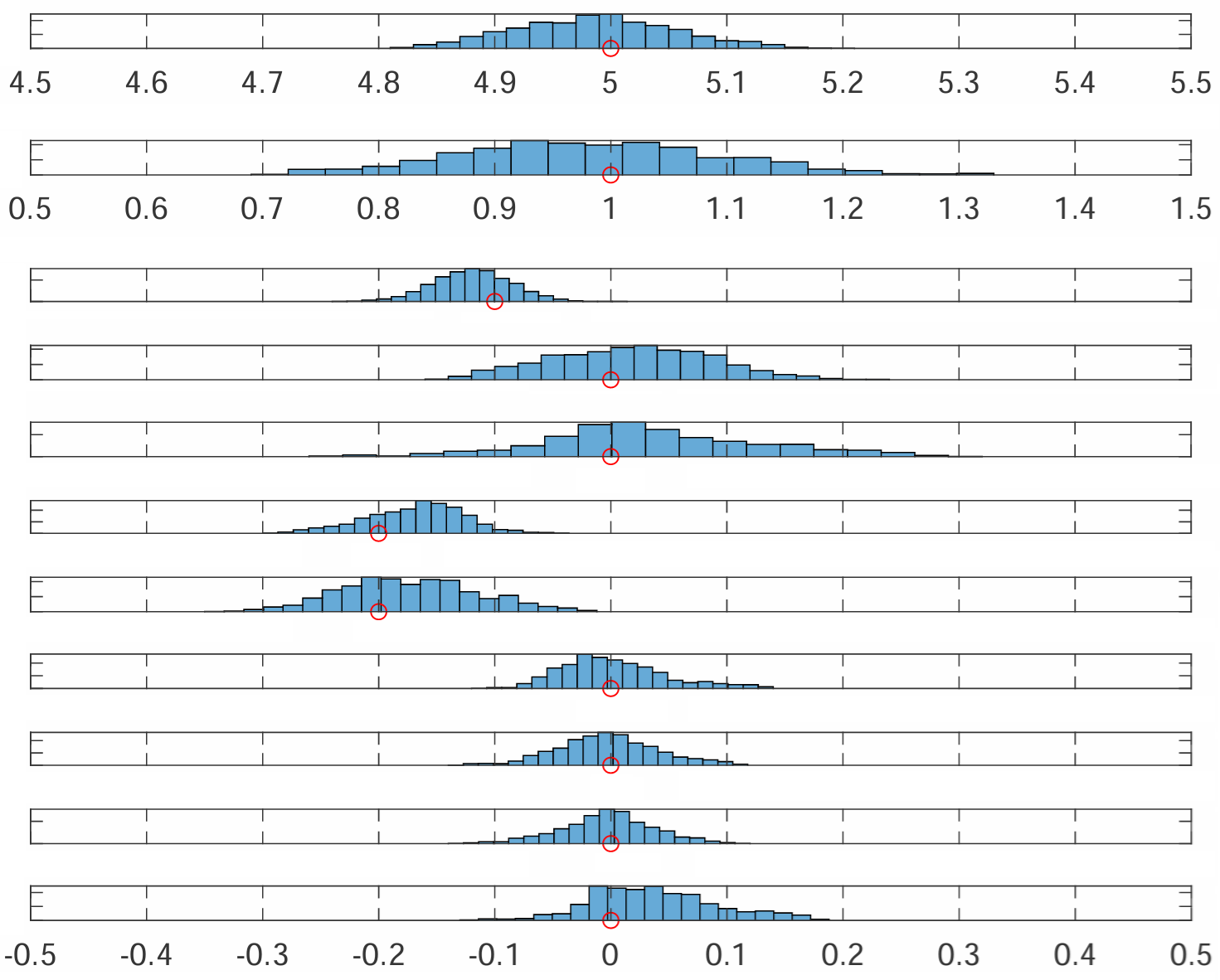}
&\includegraphics[width=0.32\linewidth]{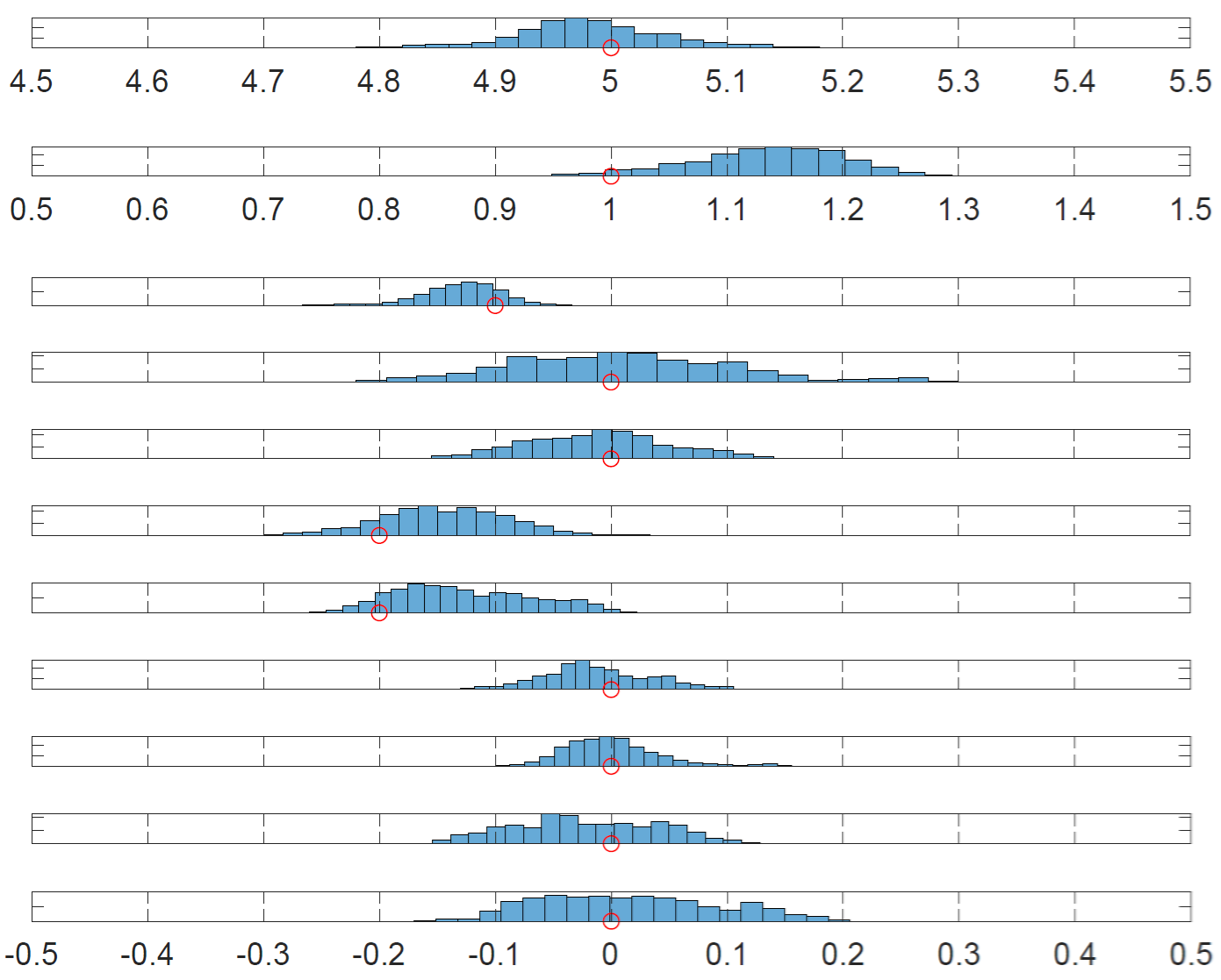}\\
\hline
\includegraphics[width=0.32\linewidth]{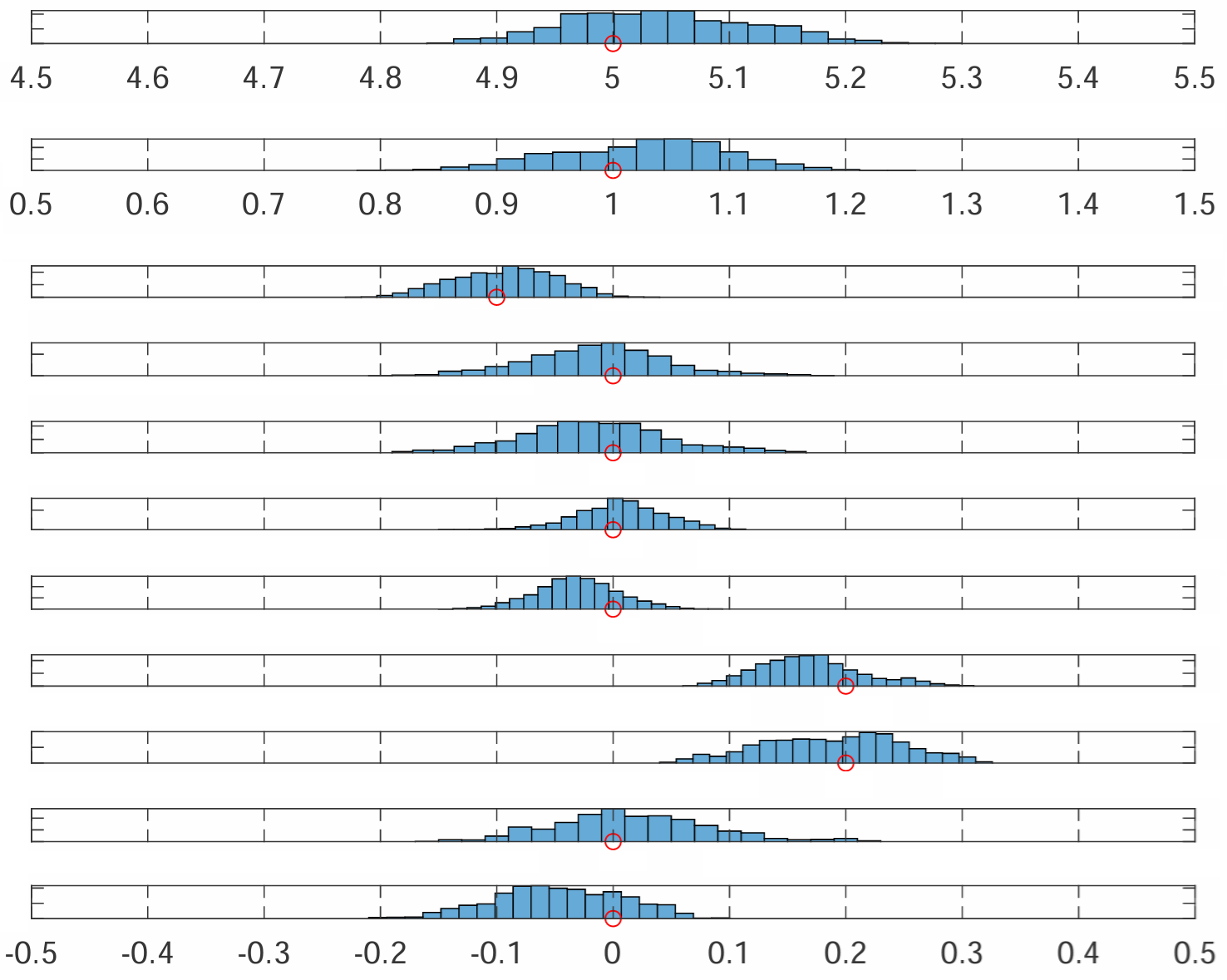}
&\includegraphics[width=0.32\linewidth]{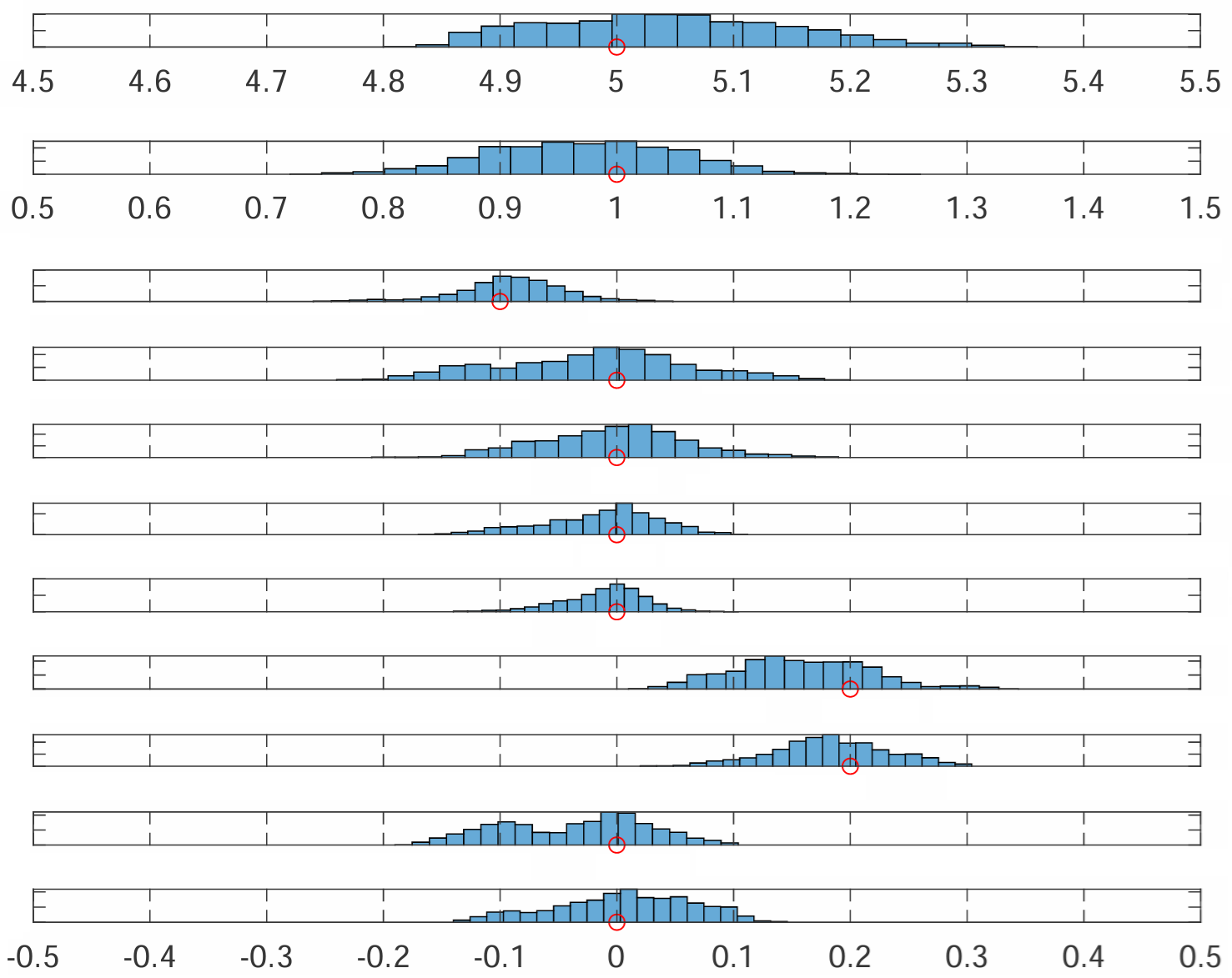}
&\includegraphics[width=0.32\linewidth]{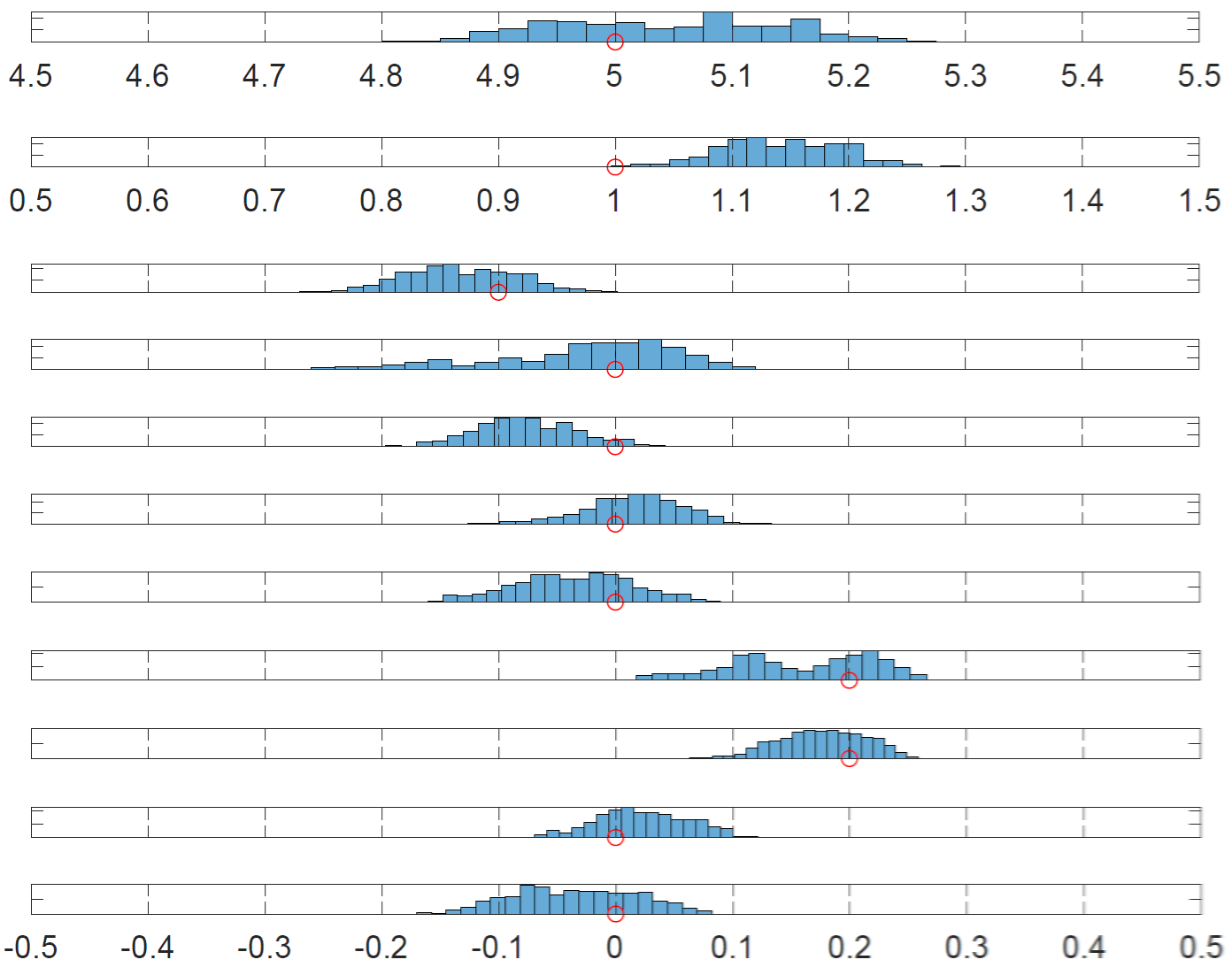}\\
\hline
\includegraphics[width=0.32\linewidth]{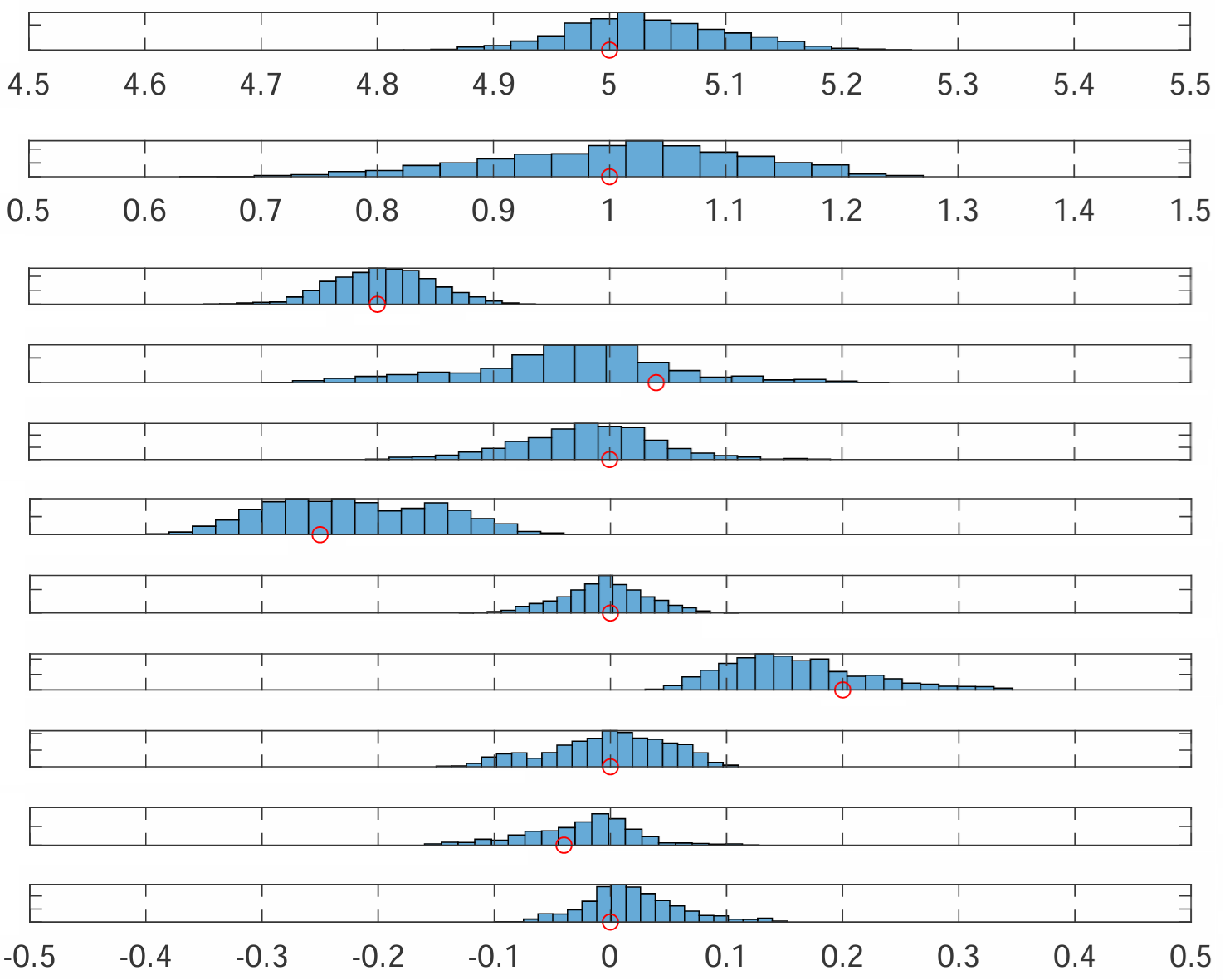}
&\includegraphics[width=0.32\linewidth]{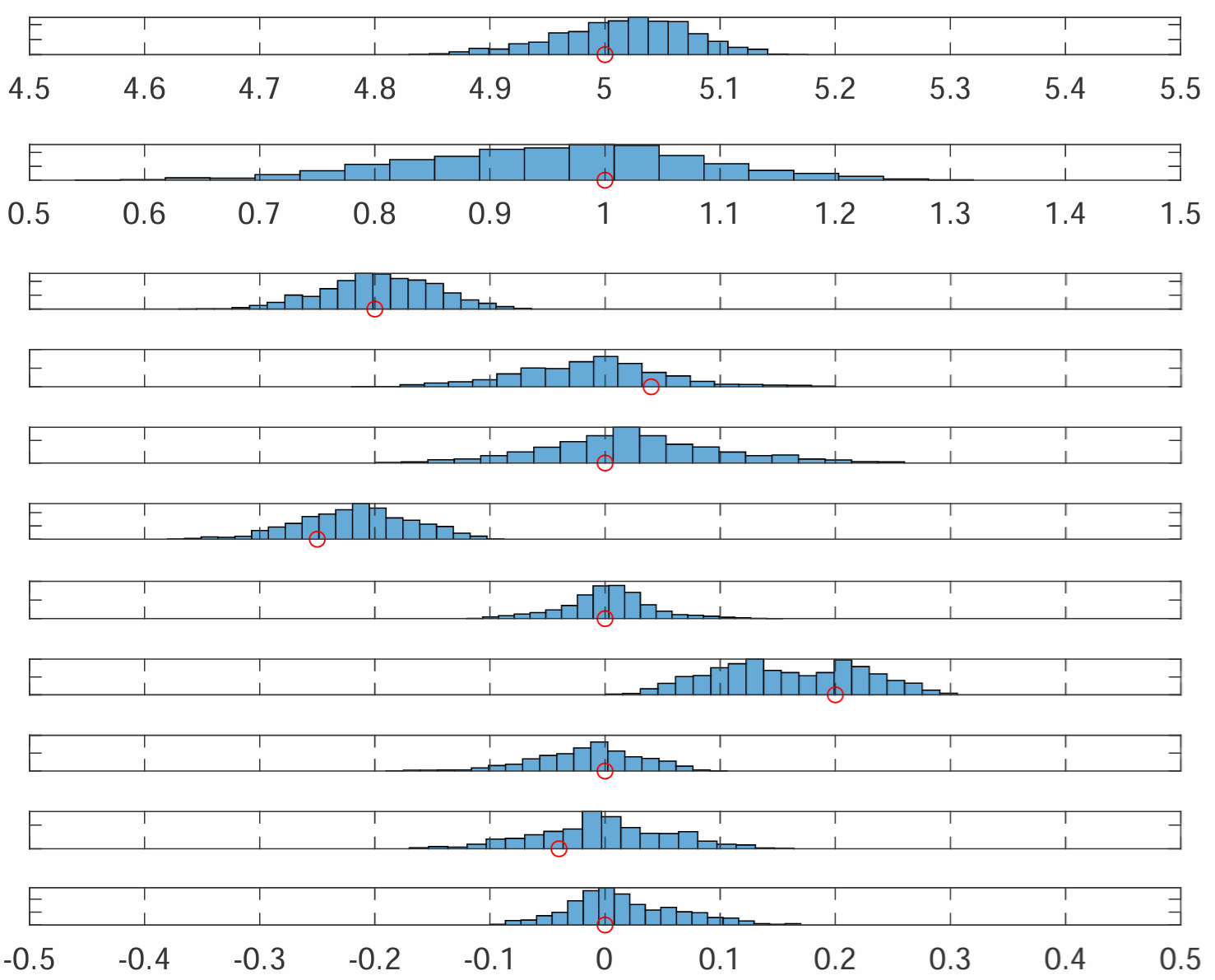}
&\includegraphics[width=0.32\linewidth]{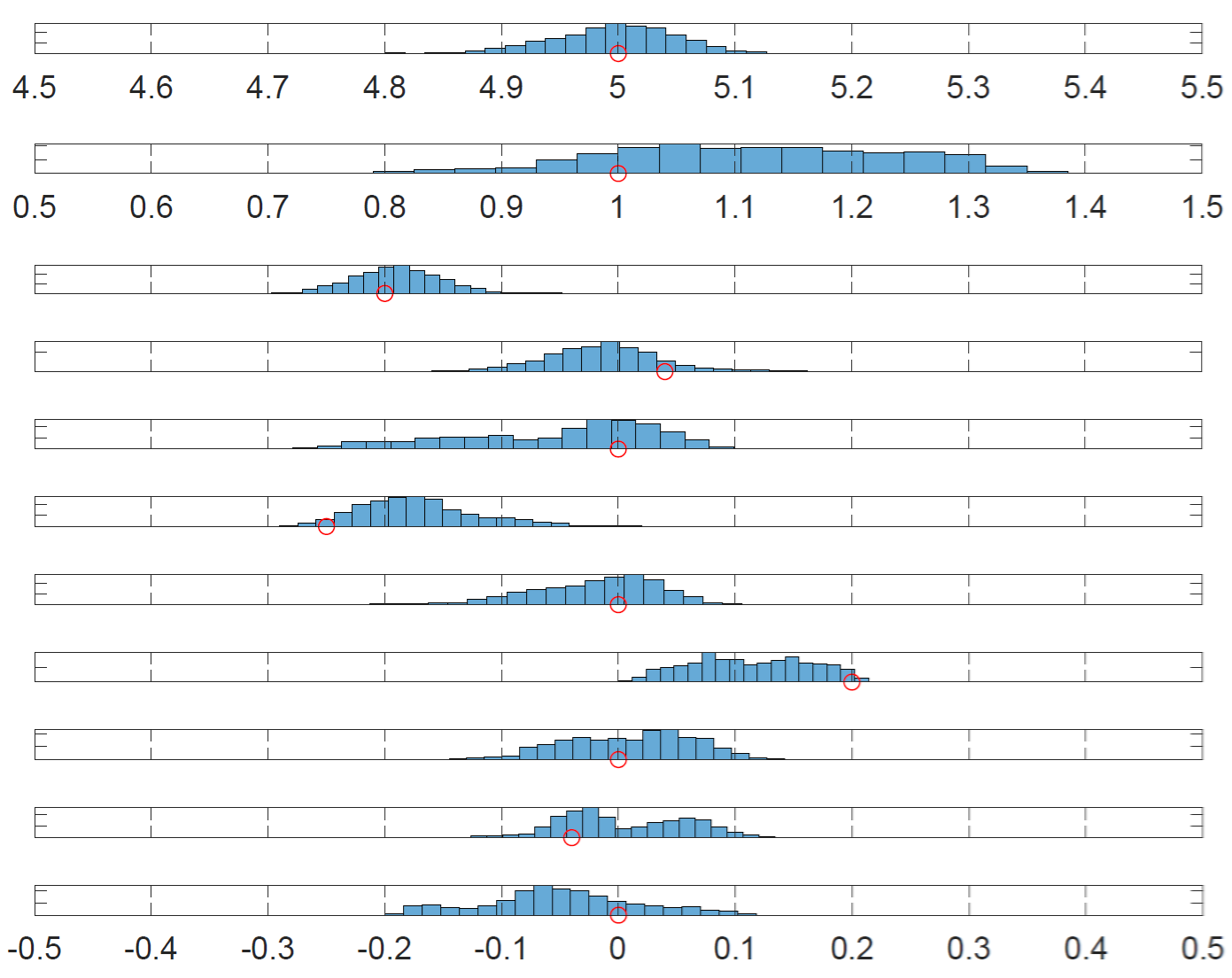}\\
\hline
\includegraphics[width=0.32\linewidth]{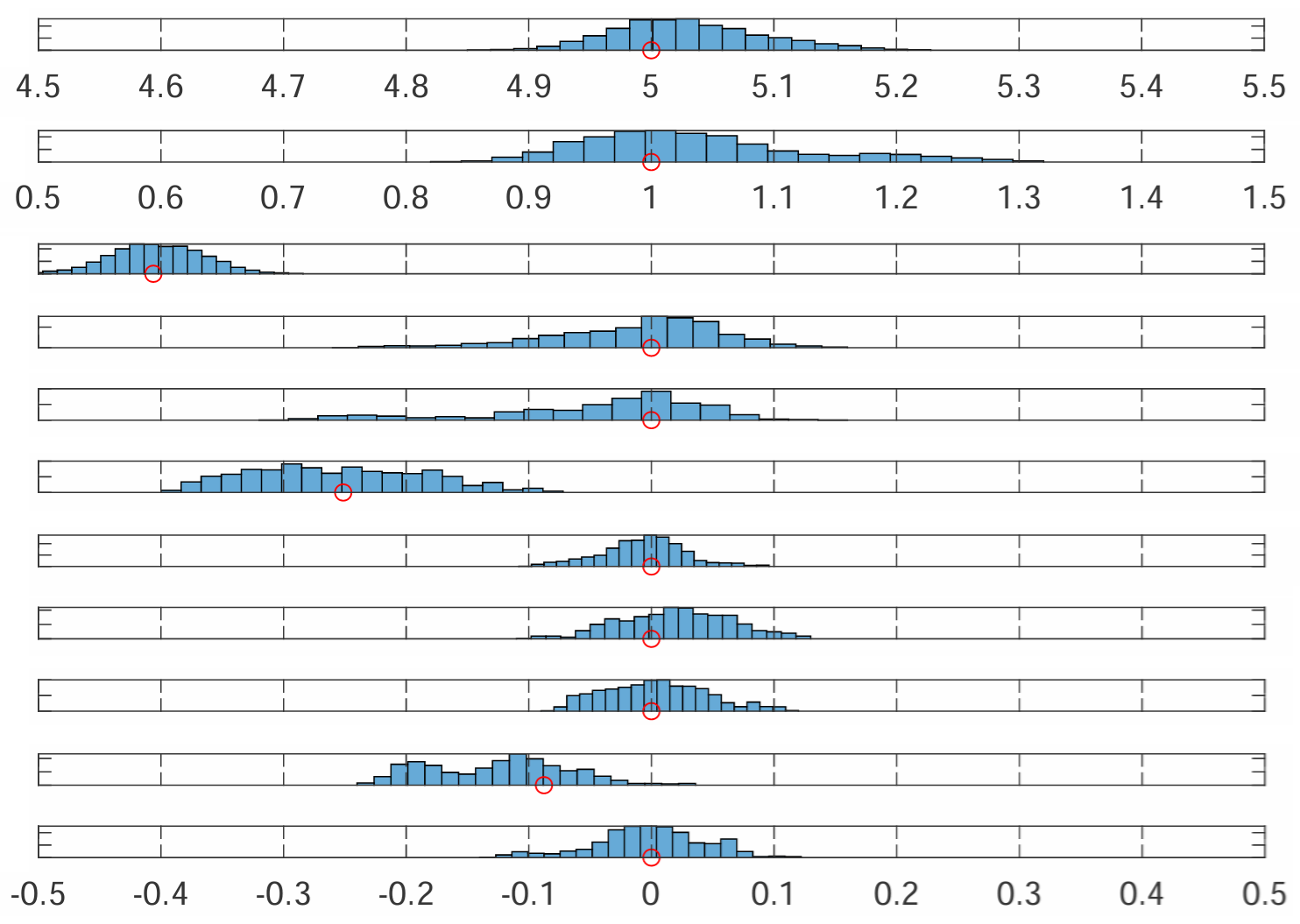}
&\includegraphics[width=0.32\linewidth]{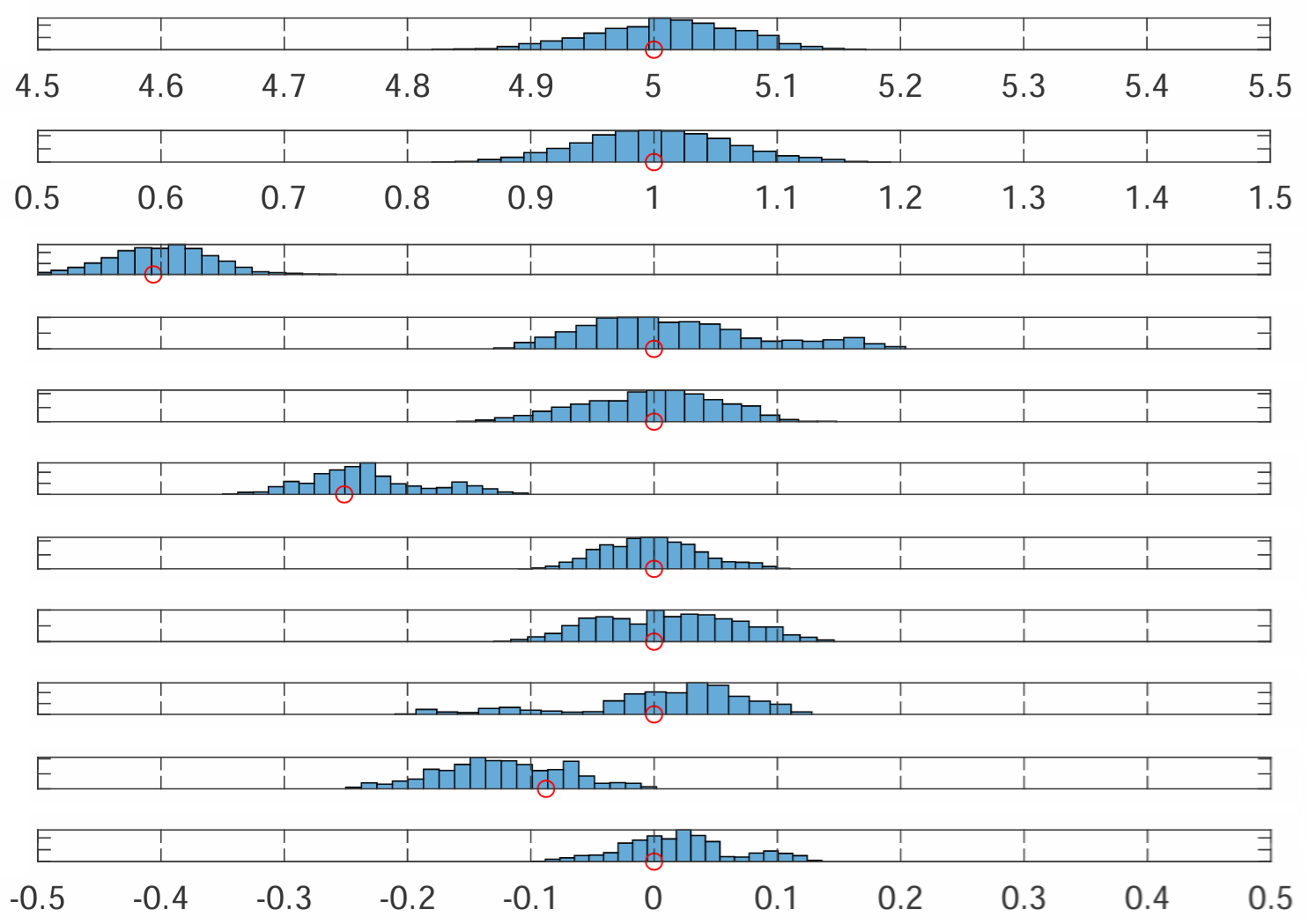}
&\includegraphics[width=0.32\linewidth]{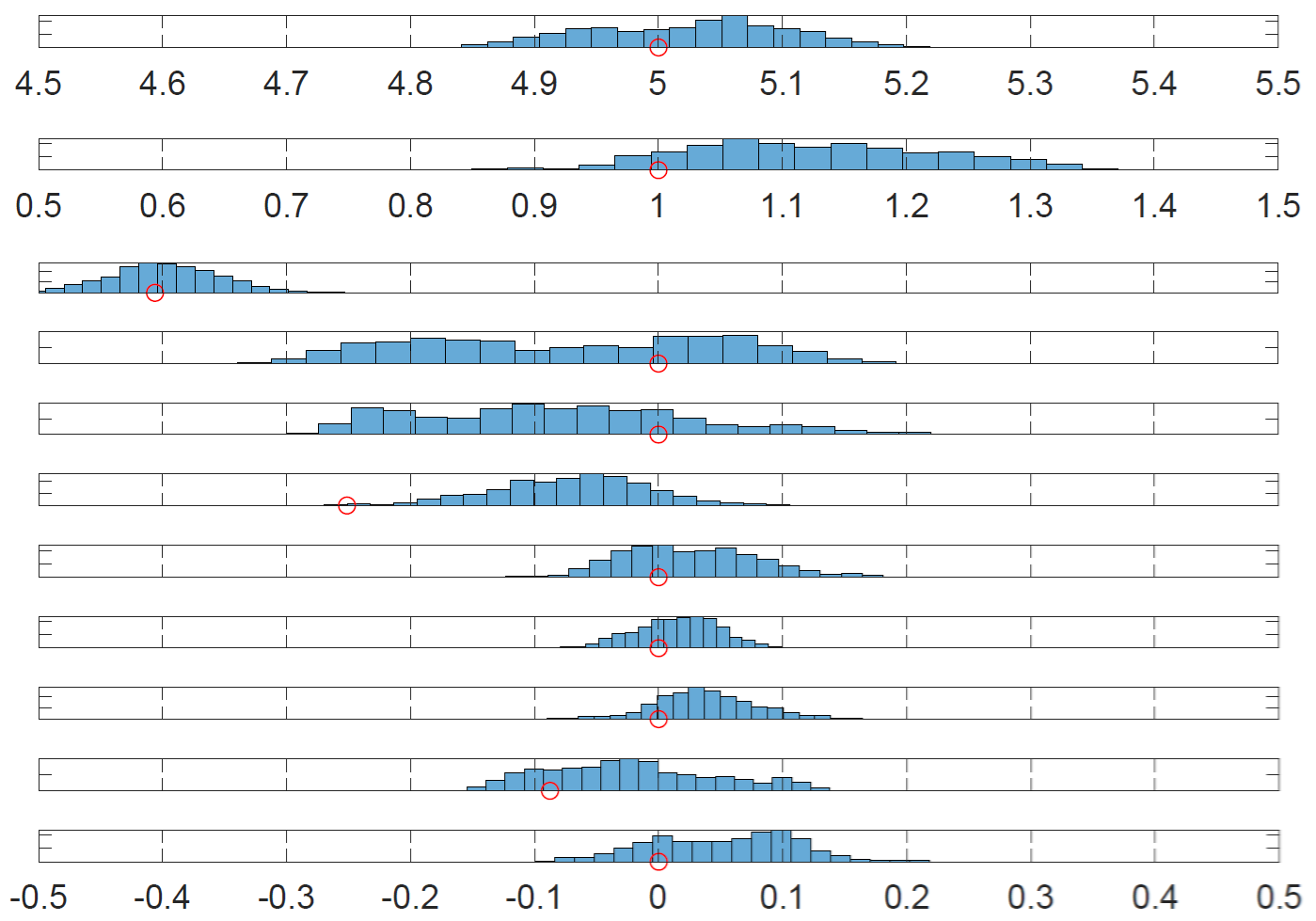}\\
\midrule
DNN & BEM & gPC ($\deg=5$)\\
\toprule
\end{tabular}
\caption{The histograms of the shape parameters $\mathbf{a}^{(i)}$ (first two rows) and $\mathbf{b}^{(i)}$ (last nine rows) for the epochs in $(9*10^4,10^5]$ for Example \ref{ex:2D:inversion} with $5\%$ noise. The red circles indicate the ground truth. From top to bottom for ellipse, pear, kite and rectangle.}
\label{fig:MCMC:2D}
\end{figure}

\begin{figure}[hbt!]
\centering\setlength{\tabcolsep}{0pt}
\begin{tabular}{ccccc}	
\includegraphics[width=0.3\linewidth]{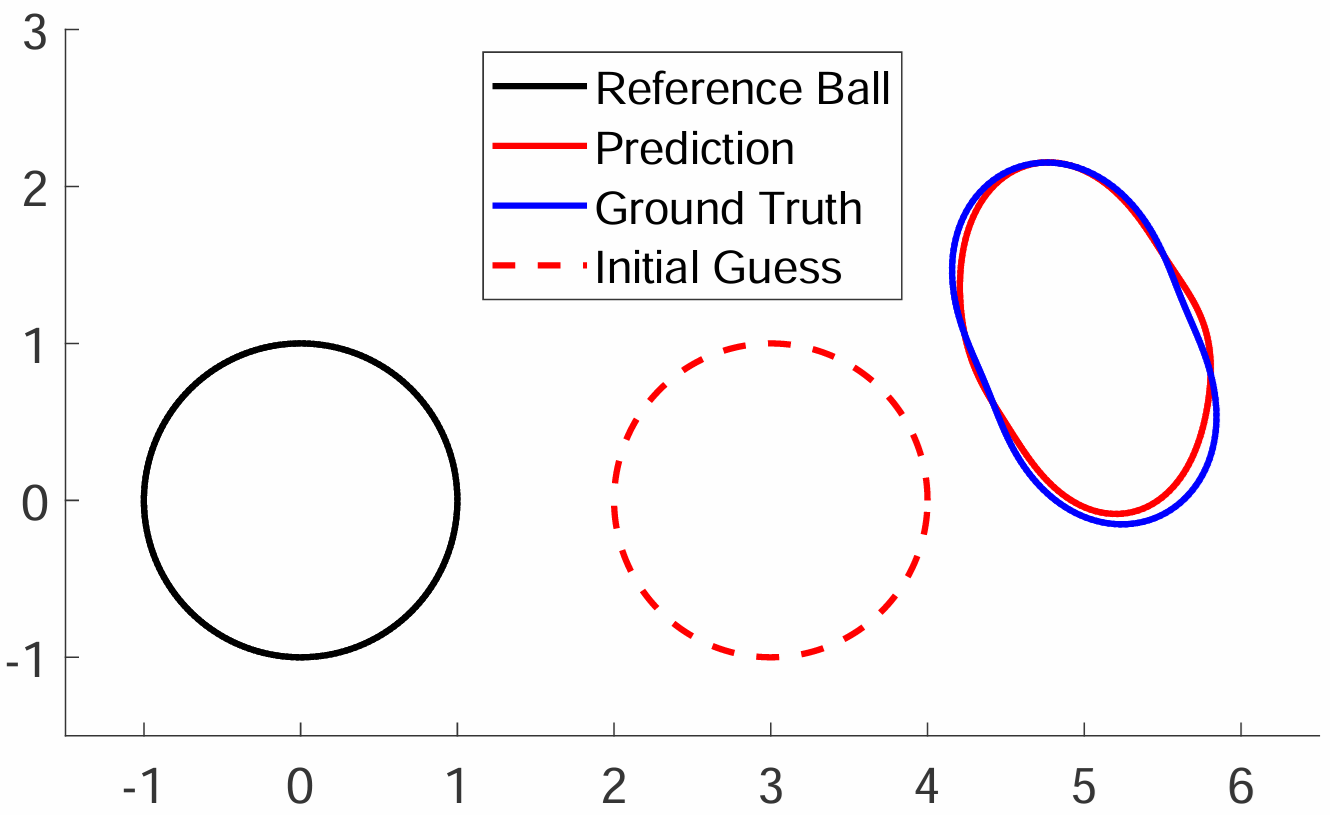}
&&\includegraphics[width=0.3\linewidth]{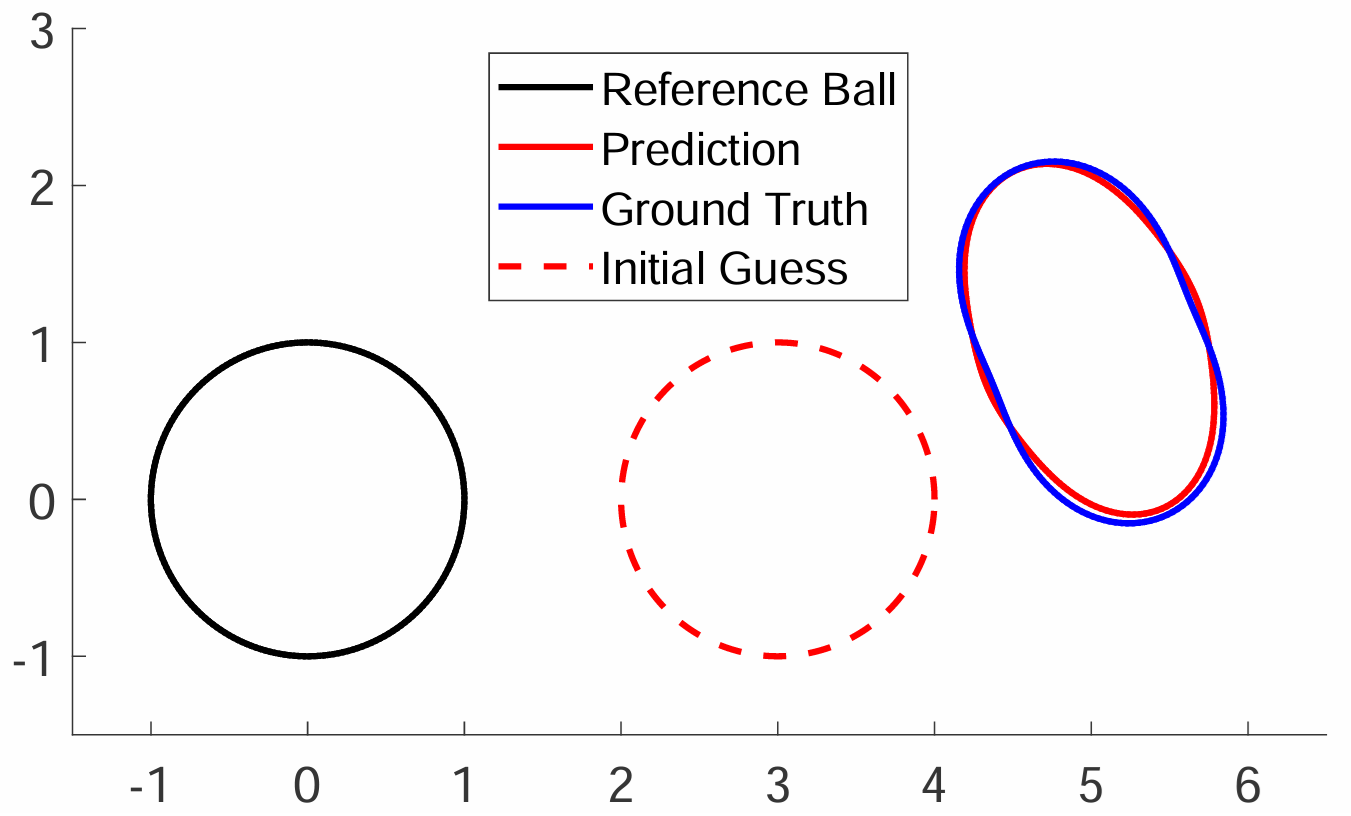}
&&\includegraphics[width=0.3\linewidth]{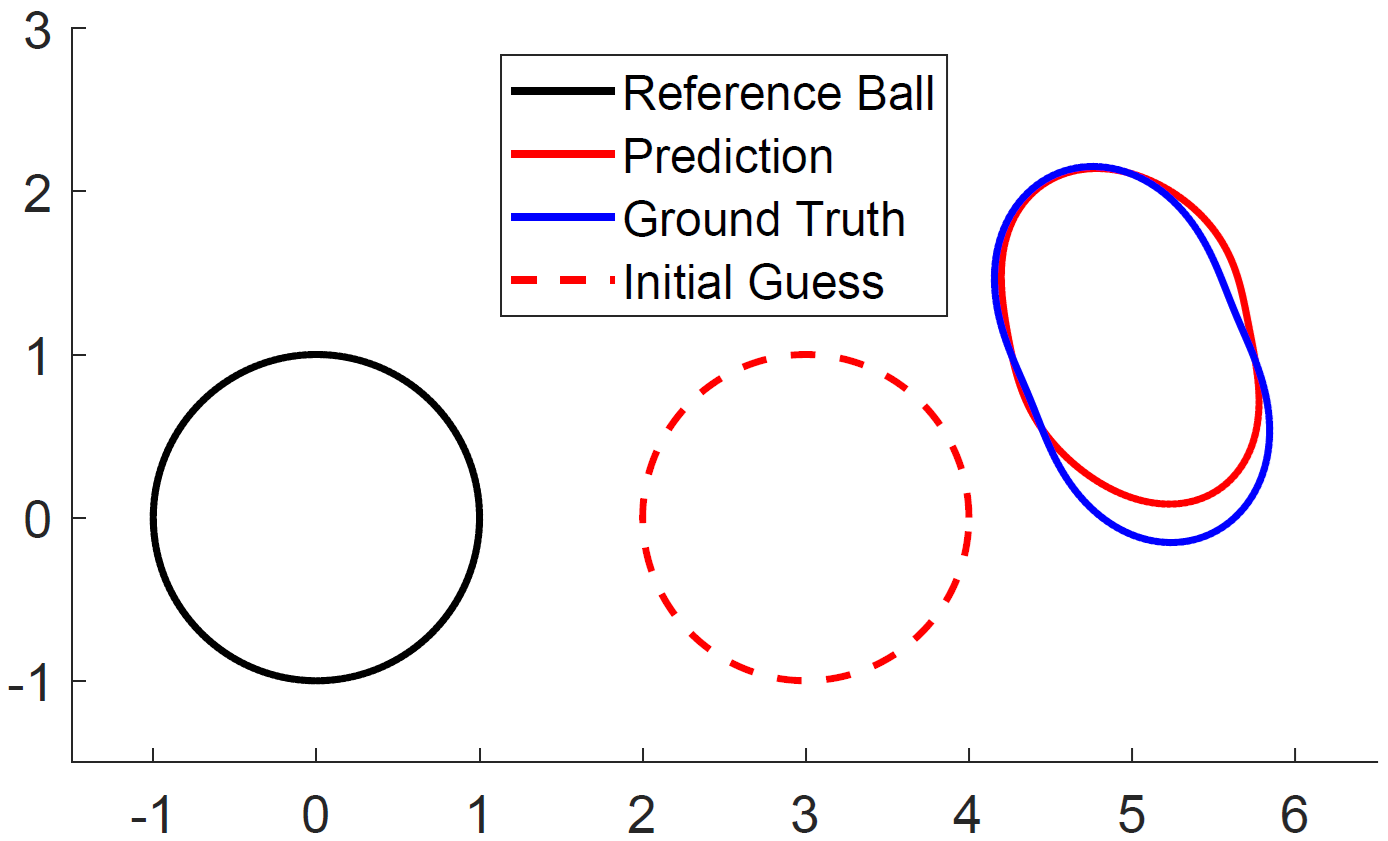}\\
&&&&\\
\includegraphics[width=0.3\linewidth, clip=true, trim = 0mm 0mm 0mm 30mm]{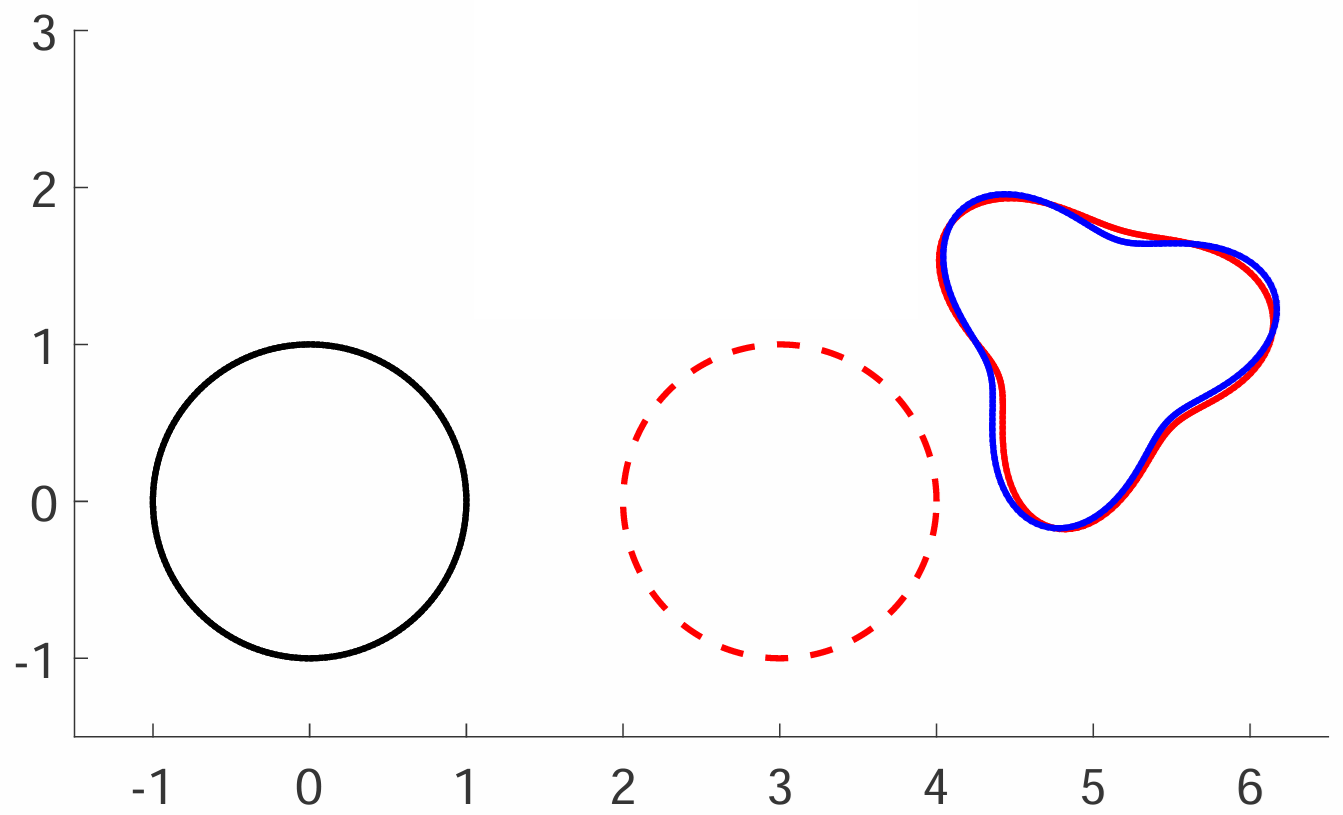}
&&\includegraphics[width=0.3\linewidth, clip=true, trim = 0mm 0mm 0mm 30mm]{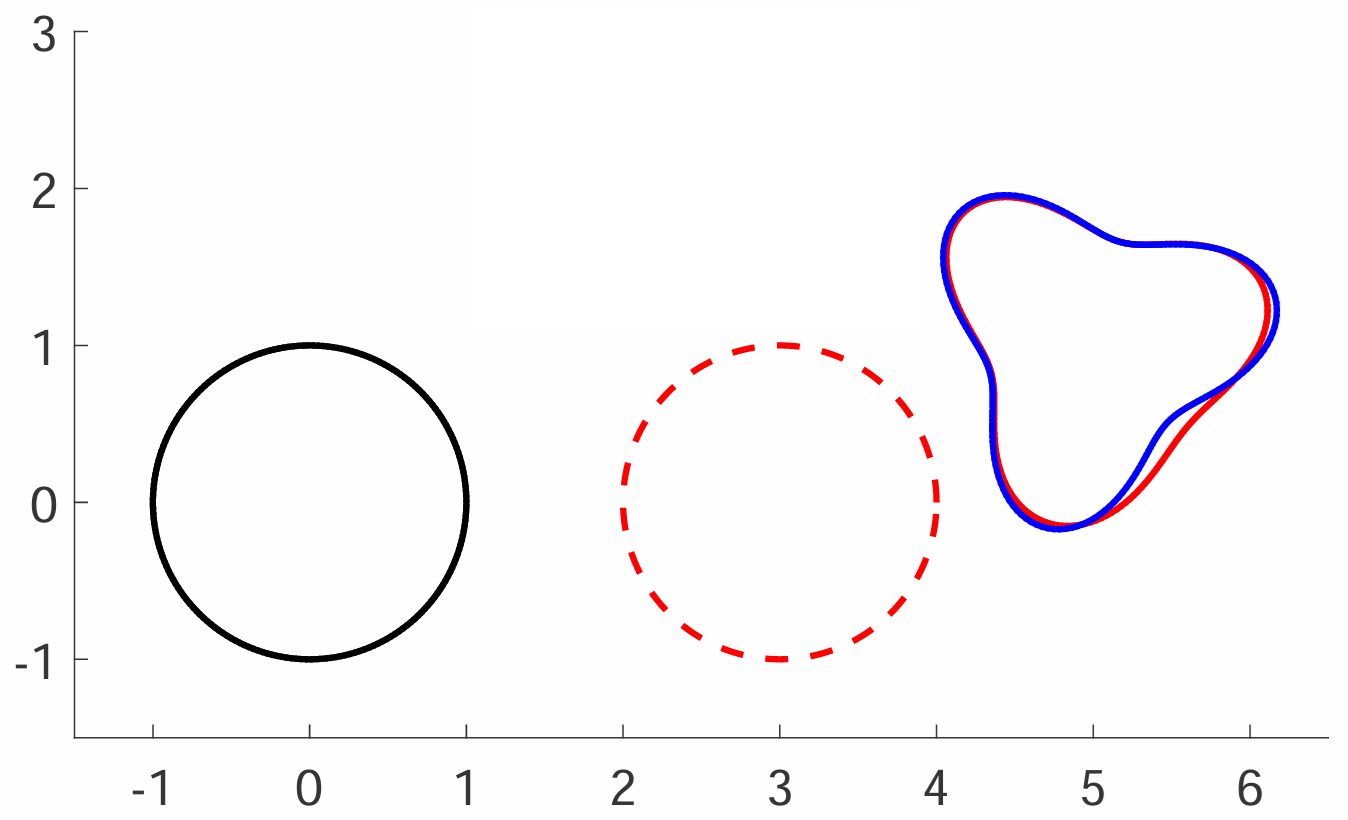}
&&\includegraphics[width=0.3\linewidth, clip=true, trim = 0mm 0mm 0mm 30mm]{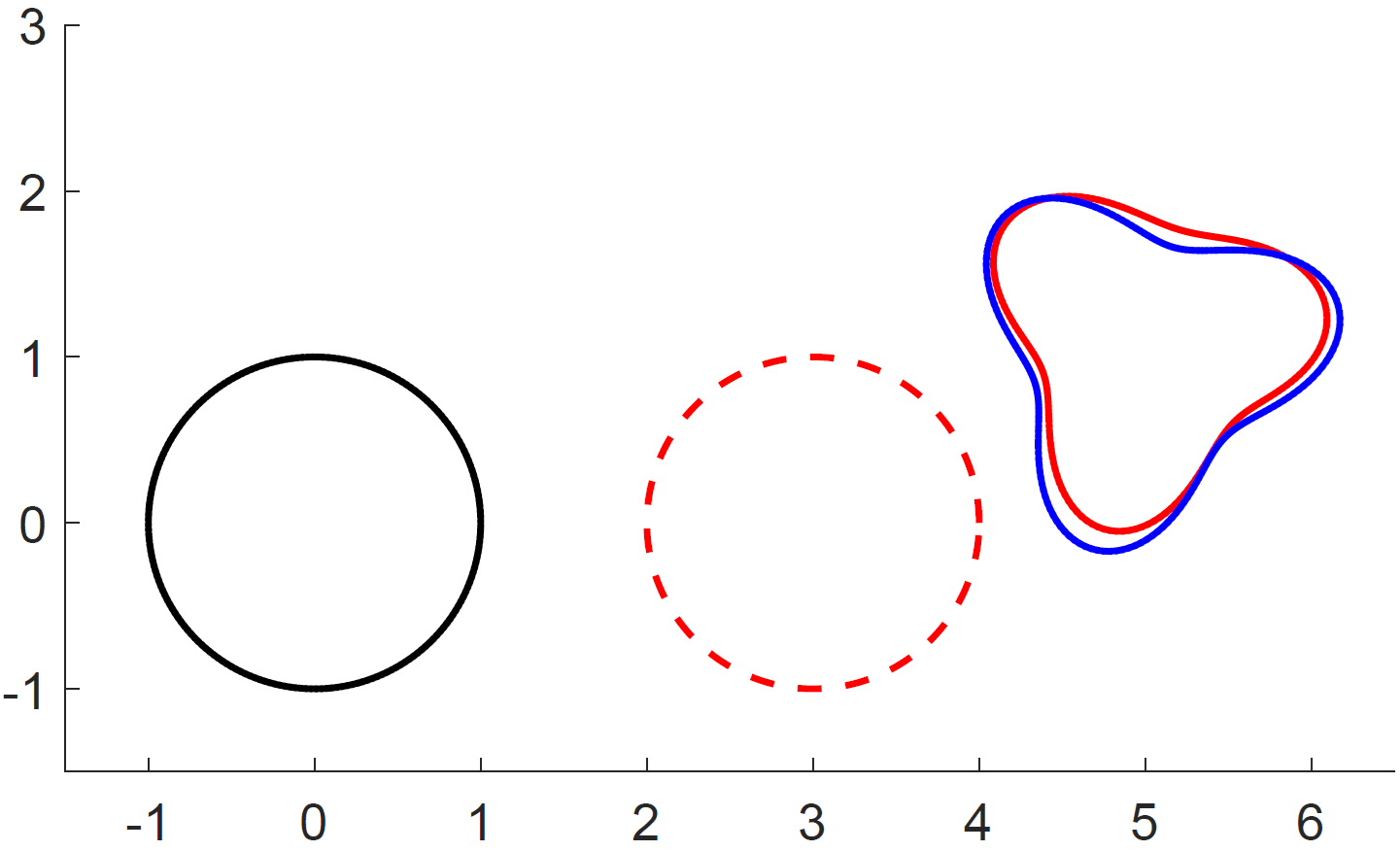}\\
&&&&\\
\raisebox{-0.5\height}{\includegraphics[width=0.3\linewidth, clip=true, trim = 0mm 0mm 0mm 50mm]{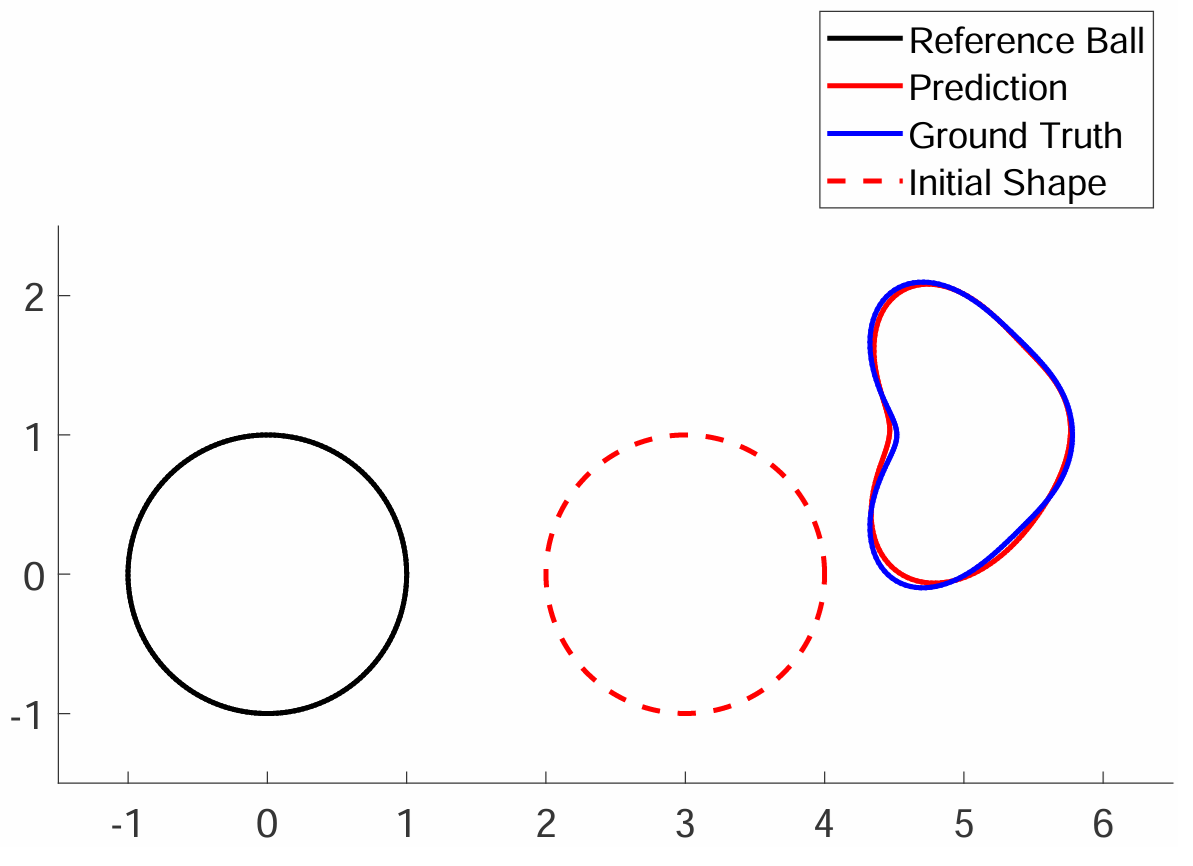}}
&&\raisebox{-0.5\height}{\includegraphics[width=0.3\linewidth, clip=true, trim = 0mm 0mm 0mm 50mm]{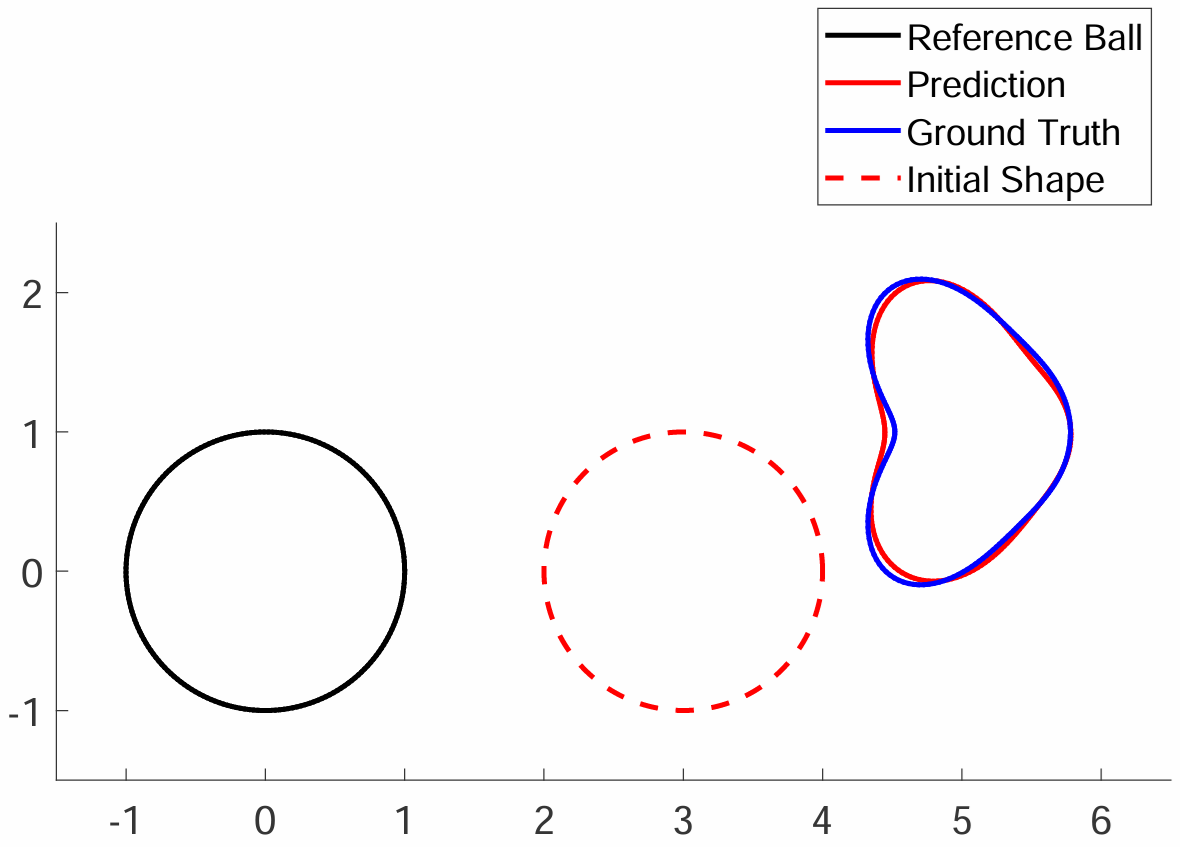}}
&&\raisebox{-0.5\height}{\includegraphics[width=0.3\linewidth, clip=true, trim = 0mm 0mm 0mm 30mm]{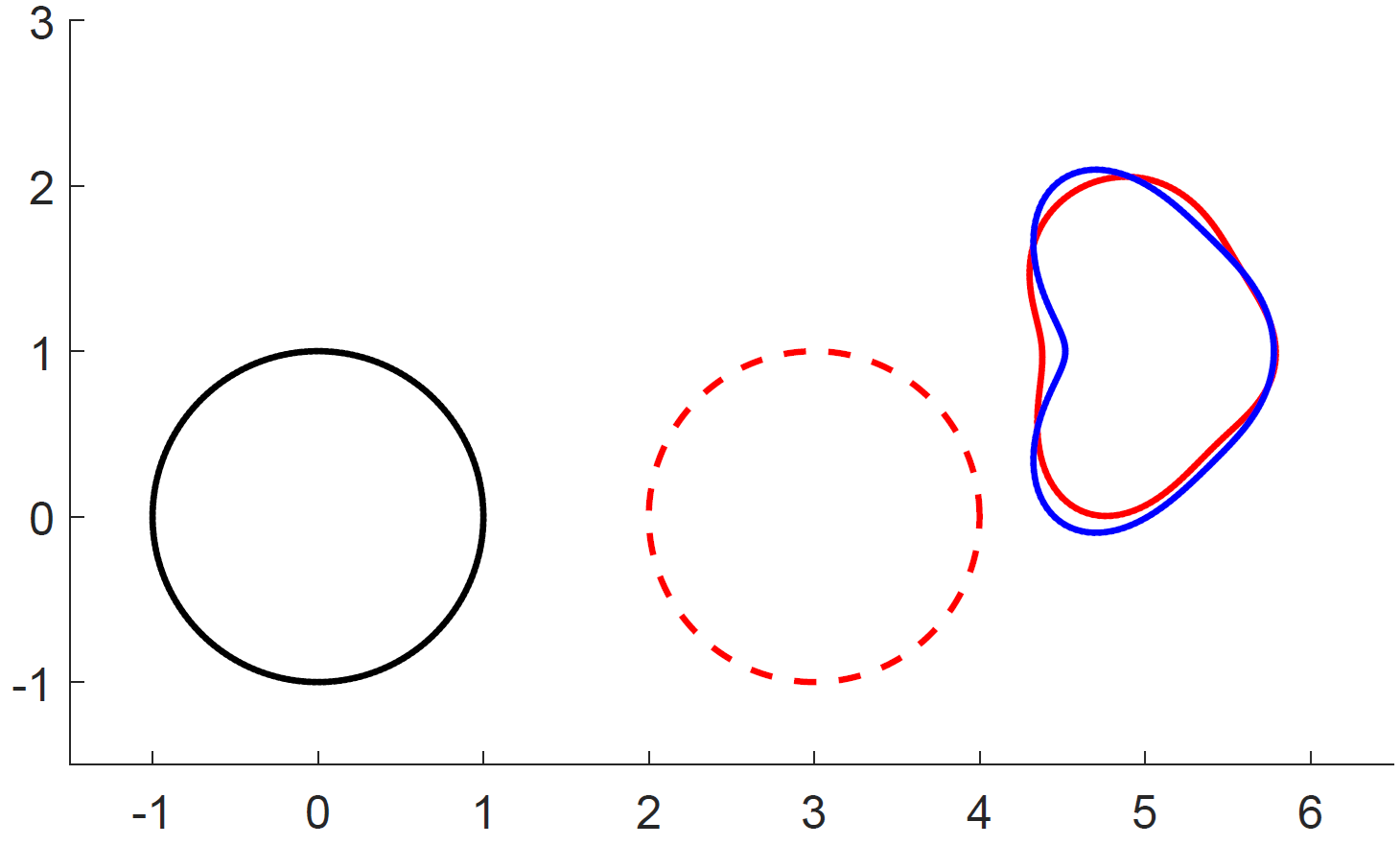}}\\
\raisebox{-0.5\height}{\includegraphics[width=0.3\linewidth, clip=true, trim = 0mm 0mm 0mm 10mm]{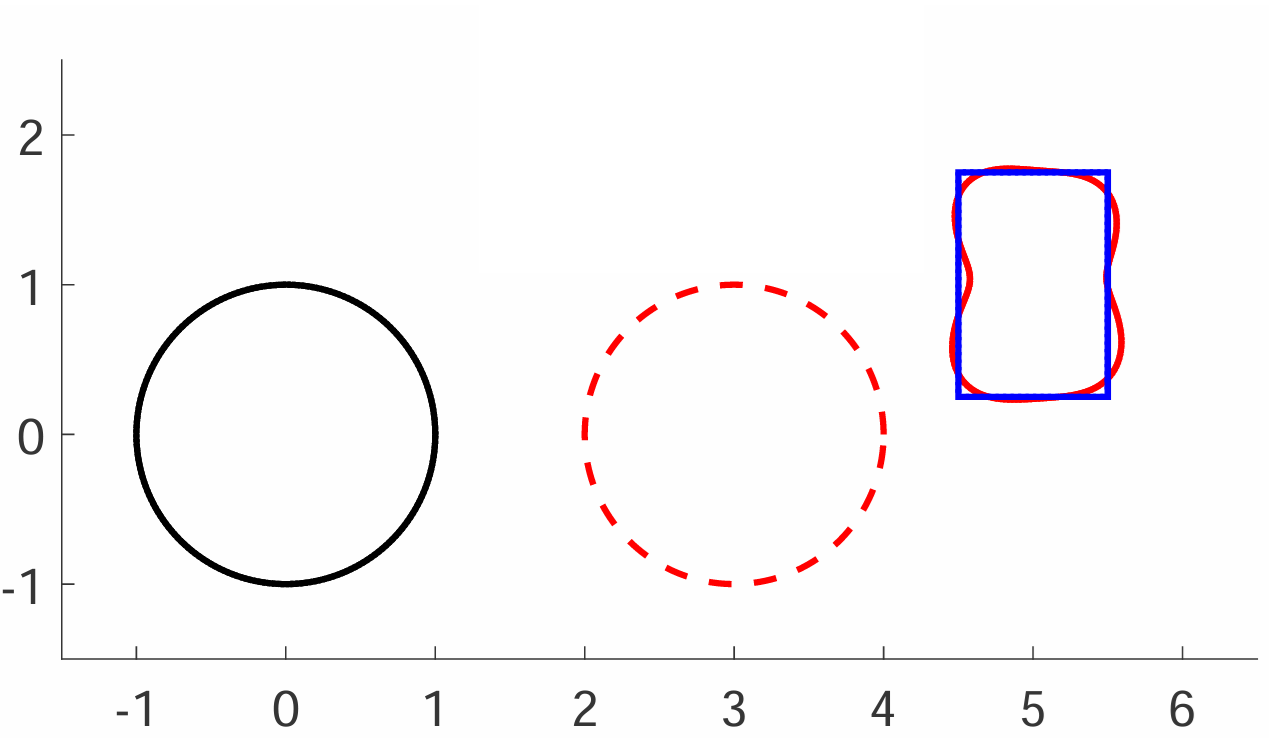}}
&&\raisebox{-0.5\height}{\includegraphics[width=0.3\linewidth, clip=true, trim = 0mm 0mm 0mm 10mm]{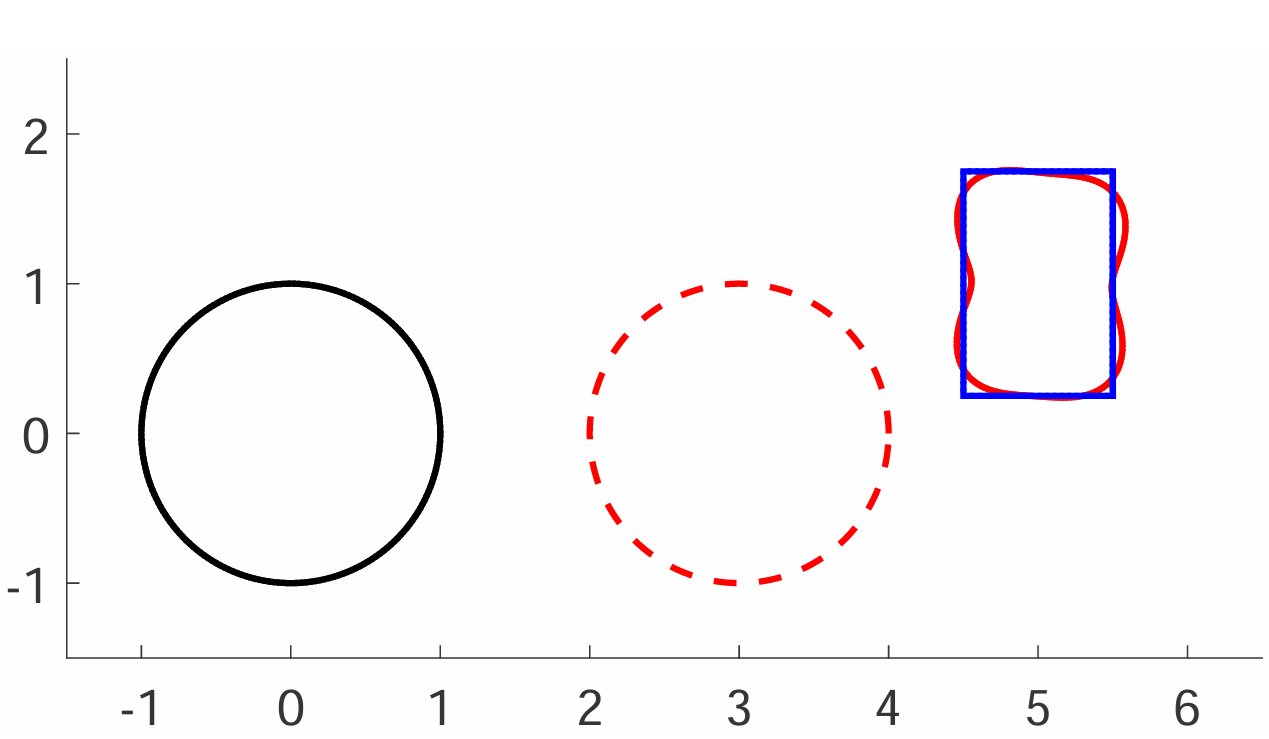}}
&&\raisebox{-0.5\height}{\includegraphics[width=0.3\linewidth, clip=true, trim = 0mm 0mm 0mm 10mm]{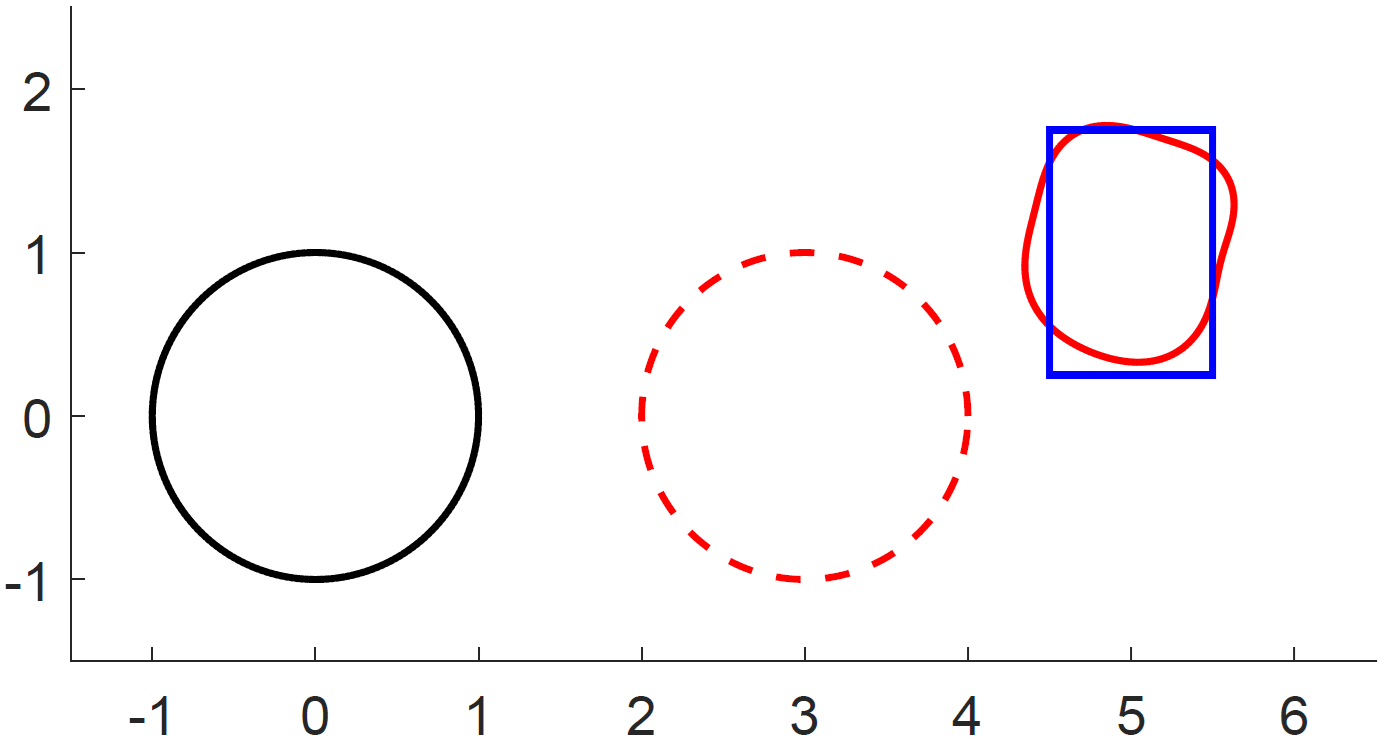}}\\
(a) DNN &\phantom{ab}& (b) BEM &\phantom{ab}& (c) gPC ($\deg=5$)
\end{tabular}
\caption{The reconstructions of the obstacle in 2D by the DNN, BEM and gPC (of degree 5) using the conditional mean estimate of the shape parameters using the epochs in $(9*10^4,10^5]$. From top to bottom for ellipse, pear, kite and rectangle.}\label{fig:2D-reconstruction}
\end{figure}

The initial shape parameters $(\mathbf{a}^{(0)},\mathbf{b}^{(0)})$ for the MCMC algorithm are given by  $\mathbf{a}^{(0)}=(3,0)$, $\rho_{\mathbf{b}^{(0)}}(\cos\theta,\sin\theta)=0$.
For the compatibility with the training data, we choose the penalty $\lambda R$ to be
$\lambda=10^3$ and $R(\Om(\mathbf{a},\mathbf{b}))=\max(0,-\log2+\sum_{m=1}^4\frac{|\widehat{{b}}_m|+|\widehat{{b}}_{-m}|}{m+2}).$ The regularizer $R$ is chosen to enforce the decay property, and the regularization parameter $\lambda$ is determined in a trial-and-error approach (whose rigorous choice is notoriously challenging \cite{ItoJin:2015}).
We set the initial learning rates to be $w=1$ and $w_{m}=1/(2|m|+4)$ for $|m|\le4$. At the end of each iteration $i$, if $\pi^{(i)}>\pi^{(0)}/3$, then we update $\pi^{(0)}\leftarrow\pi^{(0)}/3$ and $(w,(w_{m})_{|m|\le4})\leftarrow(w/2,(w_{m}/2)_{|m|\le4}).$
This is only for the initial adaptation of the step size in the MCMC algorithm. Intuitively, the strategy indicates that when the probability is substantially enhanced at an iteration so that the log probability becomes $1/3$ of the formal milestone, we decrease the step size by $1/2$. In the numerical experiments, the adaptation occurs no more than $10$ times in every MCMC chain (of length $10^5$).
In particular, the adaptation never occurs in the last $10^4$ iterations, and we only use the information in the last $10^4$ epochs to reconstruct the obstacle. We add additive Gaussian random noises to $u_{\rm m}^\infty$ and $v_{\rm m}^\infty$ with $5\%$ of $\|u_{\rm m}^\infty\|_{L^2(\mathbb{S}^1)}$ and $\|v_{\rm m}^\infty\|_{L^2(\mathbb{S}^1)}$, respectively. We compare the numerical results for the MCMC algorithm when using $\widetilde{u}_{\rm nn}^\infty$ and $\widetilde{v}_{\rm nn}^\infty$ in the algorithm in place of the classical BEM solutions. Fig. \ref{fig:MCMC:2D} indicates that the reconstructed obstacle with the DNN surrogates for the data with $5\%$ noise has a comparable accuracy to that with the BEM solver, which agrees well with the accuracy of the DNN surrogate, and the DNN approach has a big advantage in terms of the computational cost. Indeed, the total computing time  for $10^5$ MCMC iterations was less than $10^2$ seconds for neural network surrogates, whereas it is nearly $10^5$ seconds for the BEM approach, achieving a remarkable speedup factor of 1000. Fig. \ref{fig:2D-reconstruction} shows the numerical reconstructions corresponding to the mean shape parameters: the results by the neural network surrogates and the BEM are visually indistinguishable.

To provide further insights into the surrogate approach, we present quantitative results in Table \ref{tab:comparison}, including also the results by the gPC expansion. In the table, to measure the accuracy of the reconstruction obstacle, we employ the Hausdorff distance $d_{\rm H}$ and Jaccard distance $d_{\rm J}$ between two nonempty bounded domains $\Omega_1$ and $\Omega_2$ (in 2D), which are defined respectively by
\begin{align*}
    d_{\rm H} & = d_{\rm H}(\Om_1,\Om_2) = \max\left(\sup\{\operatorname{dist}(\Bx_2,\Om_1)\,:\,\Bx_2\in\Om_2\},\,\sup\{\operatorname{dist}(\Bx_1,\Om_2)\,:\,\Bx_1\in\Om_1\}\right),\\
    d_{\rm J} & = d_{\rm J}(\Om_1,\Om_2) = 1 - \frac{|\Om_1\cap\Om_2|}{|\Om_1\cup\Om_2|}.
\end{align*}
The numerical results indicate that the accuracy of the DNN surrogate is largely comparable with that by the BEM across all the noise levels. However, the DNN approach achieves significant speedup for the MCMC iterations, remarkably by a factor $10^3$ in the two-dimensional case. In contrast, the results by the gPC expansion are less accurate, due to the significant error in approximating the forward maps. Moreover, the gPC method takes longer time than the DNN surrogate for all the examples. These results clearly show the potential of the DNN surrogate approach for inverse obstacle scattering with phaseless far-field data. The results in Fig. \ref{fig:shape:evolution:MCMC} indicate that the DNN surrogate performs better than the standard gPC expansion in terms of the reconstruction accuracy for both in-distribution and out-of-distribution data, which is consistent with Fig. \ref{fig:shape:evolution:Surrogates}.

\begin{example}[MCMC in 3D]
\label{ex:3D:inversion}
The true shape parameters $\mathbf{a}^*$ and $\mathbf{b}^*=(0,\mathbf{b}_1^*,\mathbf{b}_2^*)$, with $\mathbf{b}_\ell=(b_{-\ell},b_{-\ell+1},\dots,b_{\ell})$, are given by
$\mathbf{a}^*=(2,0,2)$, $\mathbf{b}_1^*=(0,-0.3,0.3)$ and $\mathbf{b}_2^*=(0,-0.2,-0.2,-0.2,0.2)$.
\end{example}

The initial values for shape parameters are $\mathbf{a}^{(0)}=(3,0,0)$, $\mathbf{b}_1^{(0)}=(0,0,0)$ and $\mathbf{b}_2^{(0)}=(0,0,0,0,0)$. For the 3D example, the computational expense for the BEM is huge, and thus we do not present the relevant results. For the neural network surrogate, we employ the networks $\widetilde{u}_{\rm nn}^{\infty}$ and $\widetilde{v}_{\rm nn}^{\infty}$ trained in Section \ref{subsection:NN training}.
We add additive Gaussian random noises to $u_{\rm m}^\infty$ and $v_{\rm m}^\infty$ with relative noise levels $5\%$, $10\%$ and $20\%$ of $\|u_{\rm m}^\infty\|_{L^2(\mathbb{S}^2)}$ and $\|v_{\rm m}^\infty\|_{L^2(\mathbb{S}^2)}$, respectively.
We choose the penalty term $\lambda R$ to be
$$\lambda=10^3\quad\mbox{and}\quad R(\Om(\mathbf{a},\mathbf{b}))=\max\left(0,-\log2+\|(\mathbf{b}_1,\mathbf{b}_2)\|_{\ell^2}\right).$$
We set the initial learning rates $w=1$, $w_{1,m}=1/4$ for $|m|\le1$, and $w_{2,m}=1/8$ for $|m|\le2$. At the end of each iteration $i$, if $\pi^{(i)}>\pi^{(0)}/3$, then we update $\pi^{(0)}\leftarrow\pi^{(0)}/3$ and $ (w,(w_{\ell,m})_{\ell\in\{1,2\},|m|\le\ell})\leftarrow(w/2,(w_{\ell,m}/2)_{\ell\in\{1,2\},|m|\le\ell}).$
In the likelihood, we set the  standard deviation $\sigma$ for the phaseless far-field measurement to be dependent on the noise level:  $\sigma^2=10^{-4}$ for the noise level $5\%$ and $\sigma^2=5\times10^{-4}$ for the noise level  $10\%$ and $20\%$.
The numerical results are presented in Figs. \ref{fig:MCMC:3D:history} and \ref{fig:MCMC:3D:reconstruction}.
These results again show that the DNN approach can deliver reasonable numerical reconstructions for up to 20\% noise in the data, and the accuracy of the reconstruction deteriorates as the noise level increases.

\begin{table}[hbt!]
\centering
\setlength{\tabcolsep}{1pt}
\begin{tabular}{rccccccccccccccc}
\toprule
& \multicolumn{3}{c}{Ellipse} & \phantom{a} & \multicolumn{3}{c}{Pear} & \phantom{a} & \multicolumn{3}{c}{Kite} & \phantom{a} & \multicolumn{3}{c}{Rectangle}\\
\cmidrule{2-4} \cmidrule{6-8} \cmidrule{10-12} \cmidrule{14-16}
& DNN & BEM & gPC && DNN & BEM & gPC && DNN & BEM & gPC && DNN & BEM & gPC \\
\midrule
$d_{\rm H}$ \\
$5\%$ & 0.109 & 0.069 & 0.236
 && 0.066 & 0.119 & 0.139
 && 0.061 & 0.072 & 0.140
 && 0.099 & 0.121 & 0.245\\
$10\%$ & 0.075 & 0.043 & 0.210
 && 0.083 & 0.126 & 0.252
 && 0.086 & 0.080 & 0.162
 && 0.179 & 0.183 & 0.310\\
$20\%$ & 0.177 & 0.147 & 0.221
 && 0.142 & 0.133 & 0.376
 && 0.188 & 0.115 & 0.249
 && 0.184 & 0.240 & 0.169\\
$d_{\rm J}$ \\
$5\%$ & 0.091 & 0.077 & 0.156
 && 0.075 & 0.063 & 0.138
 && 0.054 & 0.063 & 0.133
 && 0.097 & 0.097 & 0.229\\
$10\%$ & 0.058 & 0.042 & 0.183
 && 0.064 & 0.089 & 0.175
 && 0.064 & 0.092 & 0.143
 && 0.119 & 0.106 & 0.221\\
$20\%$ & 0.118 & 0.097 & 0.173
 && 0.095 & 0.098 & 0.184
 && 0.144 & 0.125 & 0.163
 && 0.141 & 0.156 & 0.245
\\
Time \\
$10\%$ & 6.8e+1 & 9.9e+4 & 7.6e+2
 && 6.3e+1 & 9.4e+4 & 7.8e+2
 && 6.5e+1 & 9.4e+4 & 7.9e+2
 && 6.1e+1 & 8.9e+4 & 7.8e+2
\\
 \bottomrule
\end{tabular}
\caption{The accuracy and computing time for the MCMC reconstruction for Example \ref{ex:2D:inversion} with various noise levels using different solvers (DNN, BEM and gPC (of degree 5). The accuracy is measured for Hausdorff distance ($d_{\rm H}$) and Jaccard distance ($d_{\rm J}$) between the exact $\Om$ and the prediction by the mean of last $10^4$ accepted shape parameters.} \label{tab:comparison}
\end{table}

\begin{figure}[hbt!]
\centering\setlength{\tabcolsep}{0pt}
    \begin{tabular}{ccc}
        \includegraphics[width=0.48\linewidth]{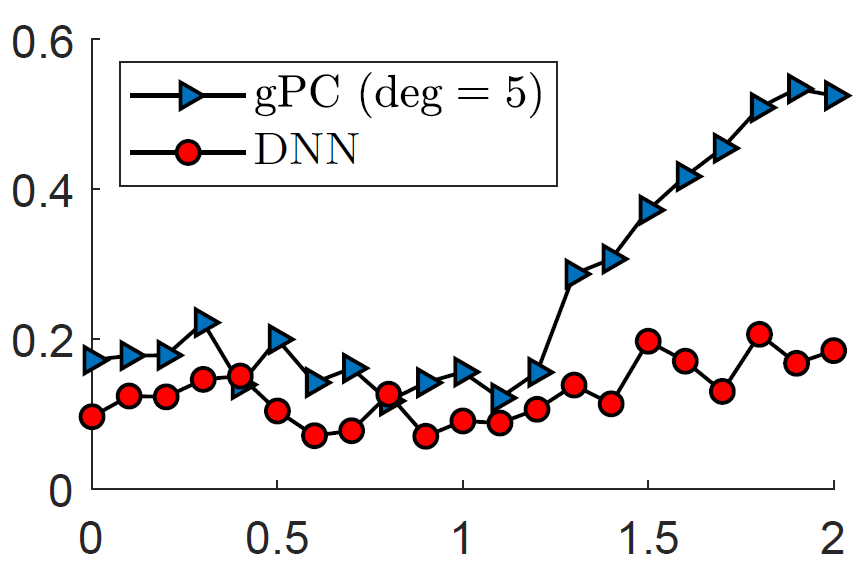}
         &\phantom{a}&
        \includegraphics[width=0.48\linewidth]{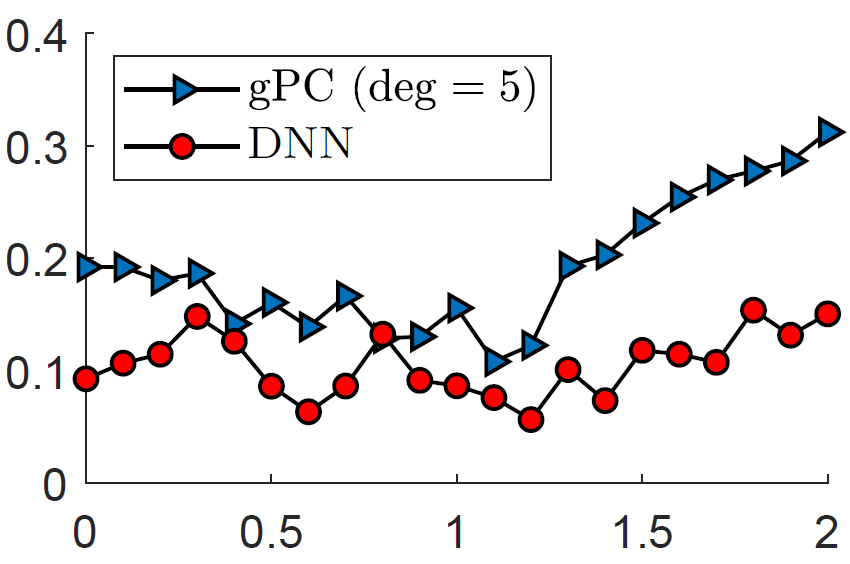}
         \\
         (a) $d_{\rm H}$ && (b) $d_{\rm J}$
    \end{tabular}
    \caption{The  accuracy (in terms of the Hausdorff distance $d_{\rm H}$ and Jaccard distance $d_{\rm J}$) of the MCMC algorithm combined with the DNN and gPC surrogates. The true shape $\Omega$ depends on the parameter $\alpha$ (x-axis), cf. also Fig. \ref{fig:shape:evolution:Surrogates}. The results are for Example \ref{ex:2D:inversion} with the phaseless far-field data with $5\%$ additive white Gaussian noise. The threshold $\log2$ defining $R$ is replaced by $\log3$ to allow more irregularity for the target domain.}
    \label{fig:shape:evolution:MCMC}
\end{figure}

\begin{figure}[hbt!]
\centering\setlength{\tabcolsep}{0pt}
\begin{tabular}{ccc}	
\includegraphics[width=0.33\linewidth]{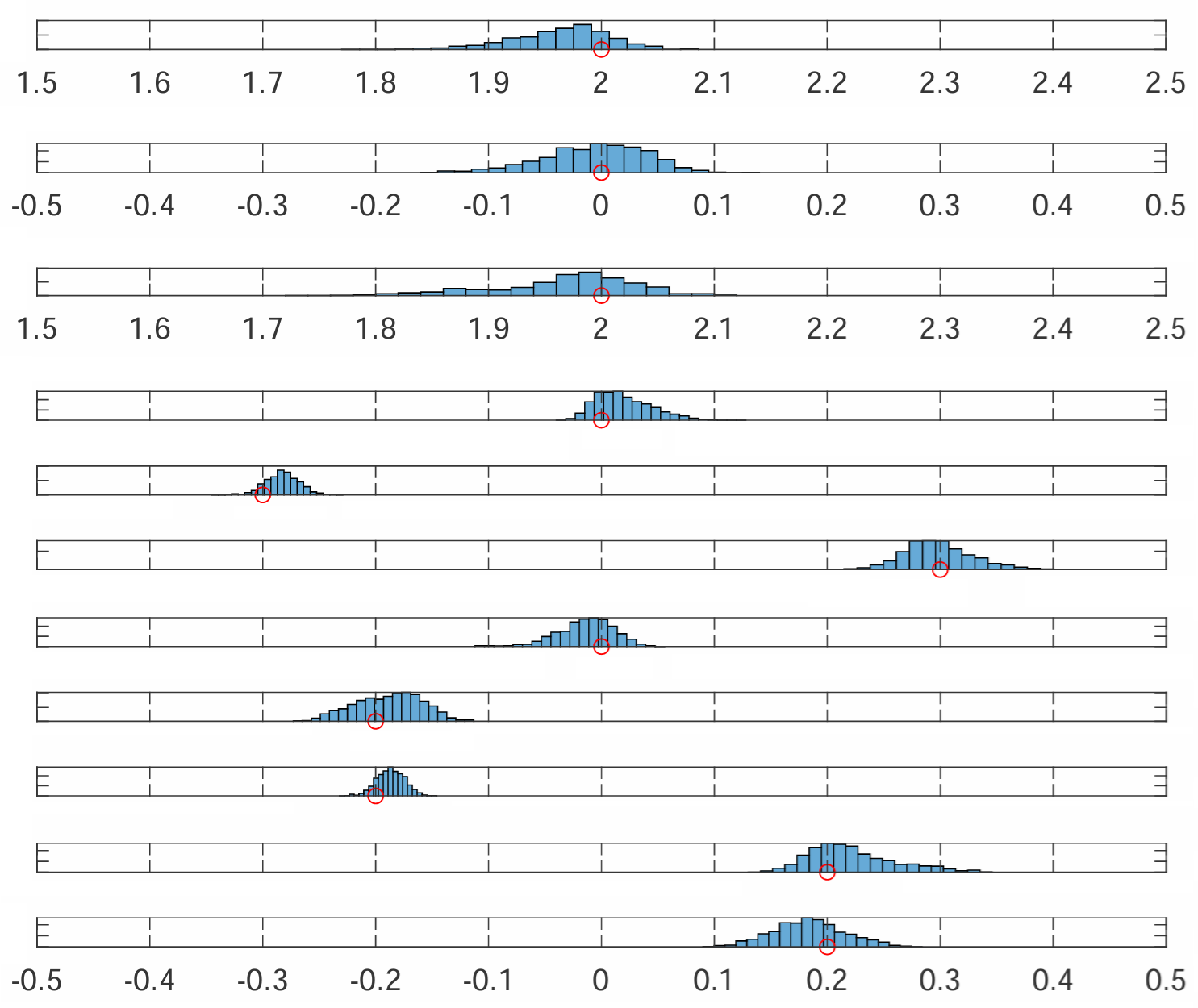} & \includegraphics[width=0.33\linewidth]{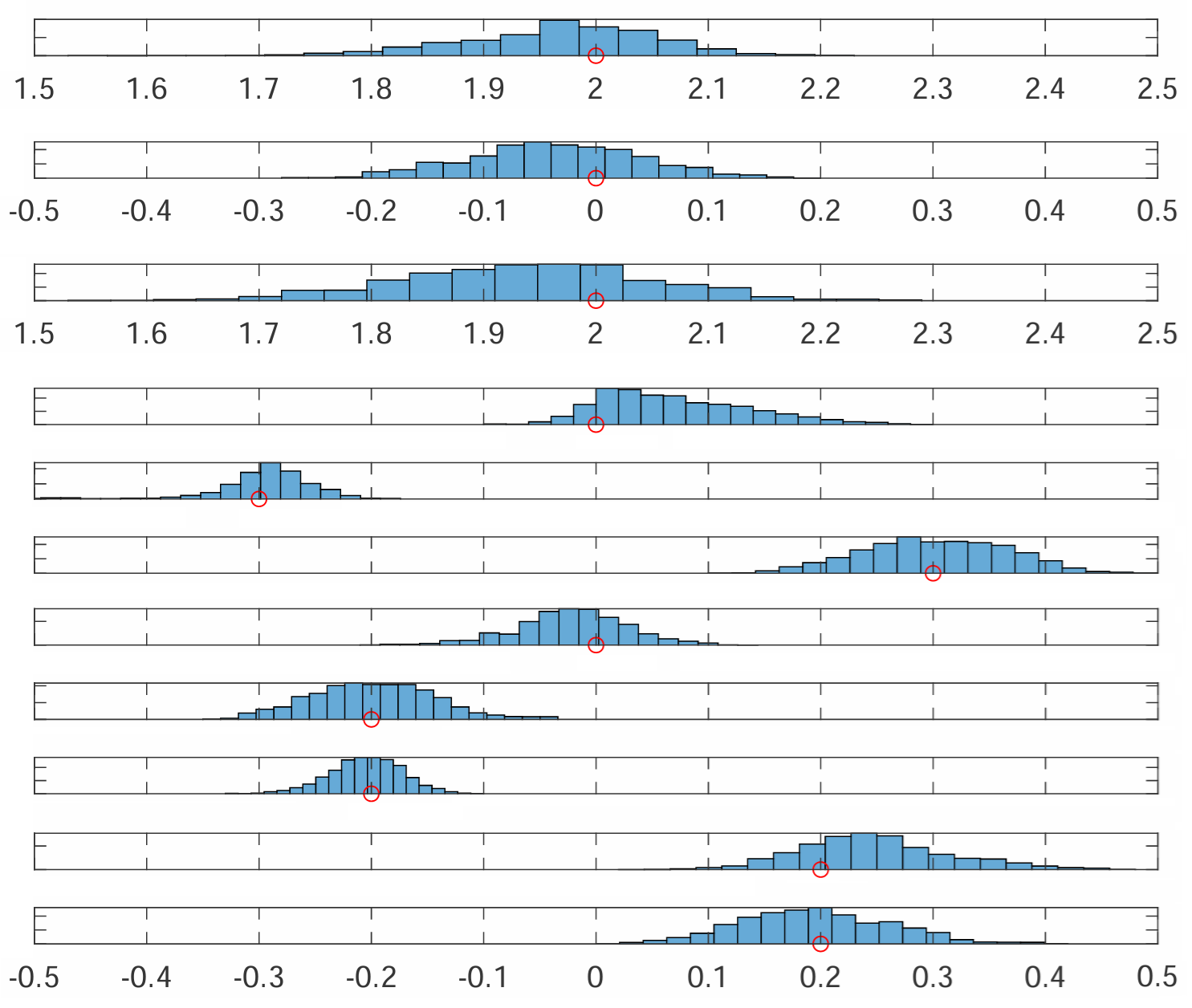}& \includegraphics[width=0.33\linewidth]{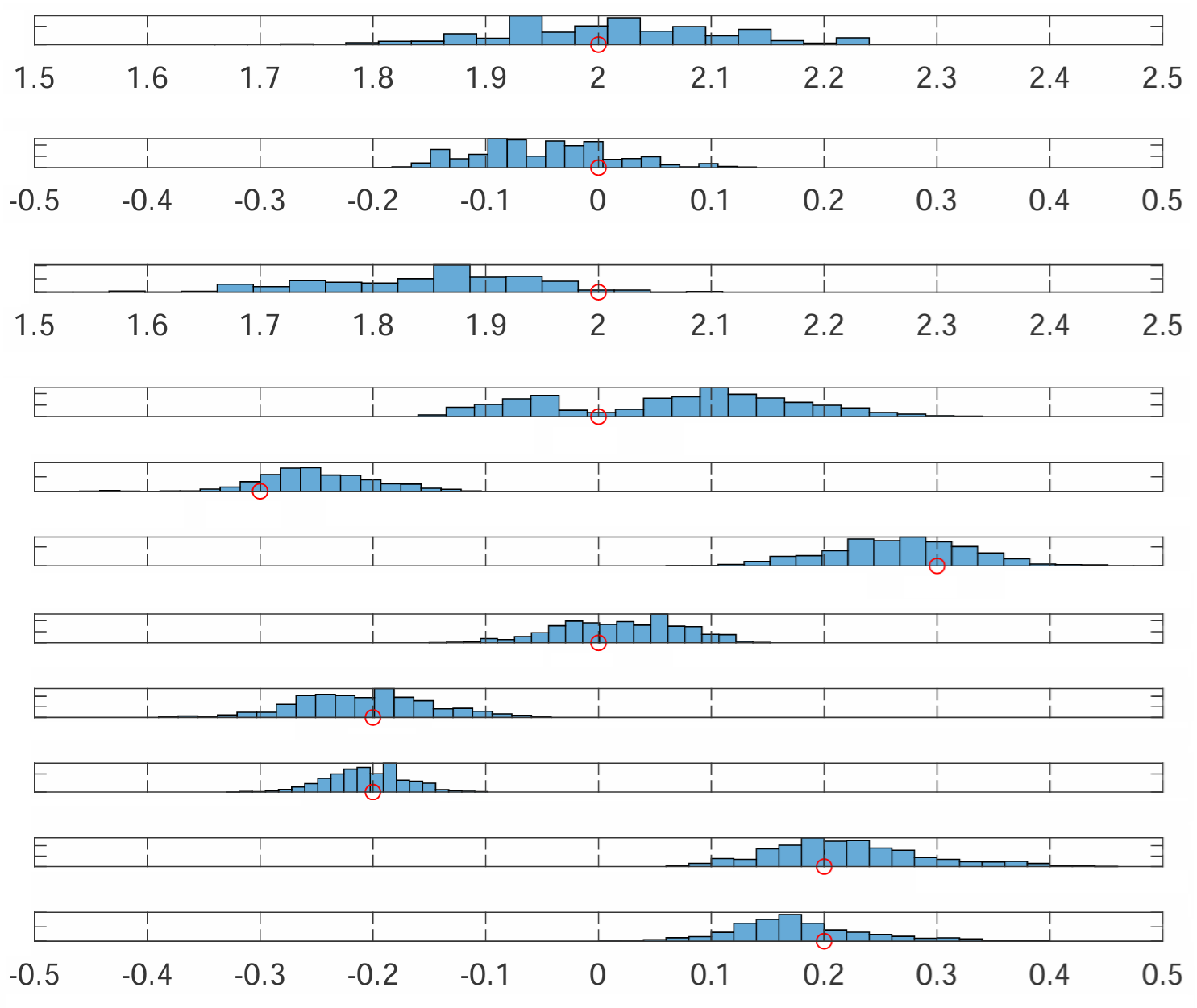}
\\
(a) $5\%$ & (b) $10\%$ & (c) $20\%$
\end{tabular}
\caption{
The histograms of the shape parameters $\mathbf{a}^{(i)}$ (first three rows), $\mathbf{b}_1^{(i)}$ (rows from 4th to 6th) and $\mathbf{b}_2^{(i)}$ (rows from 7th to 11th) for the epochs in $(9*10^4,10^5]$ for Example \ref{ex:3D:inversion} at three noise levels. The red circles indicate the ground truth.}
\label{fig:MCMC:3D:history}
\end{figure}

\begin{figure}[hbt!]
\centering\setlength{\tabcolsep}{0pt}
\begin{tabular}{ccccccc}	
\includegraphics[width=0.28\linewidth]{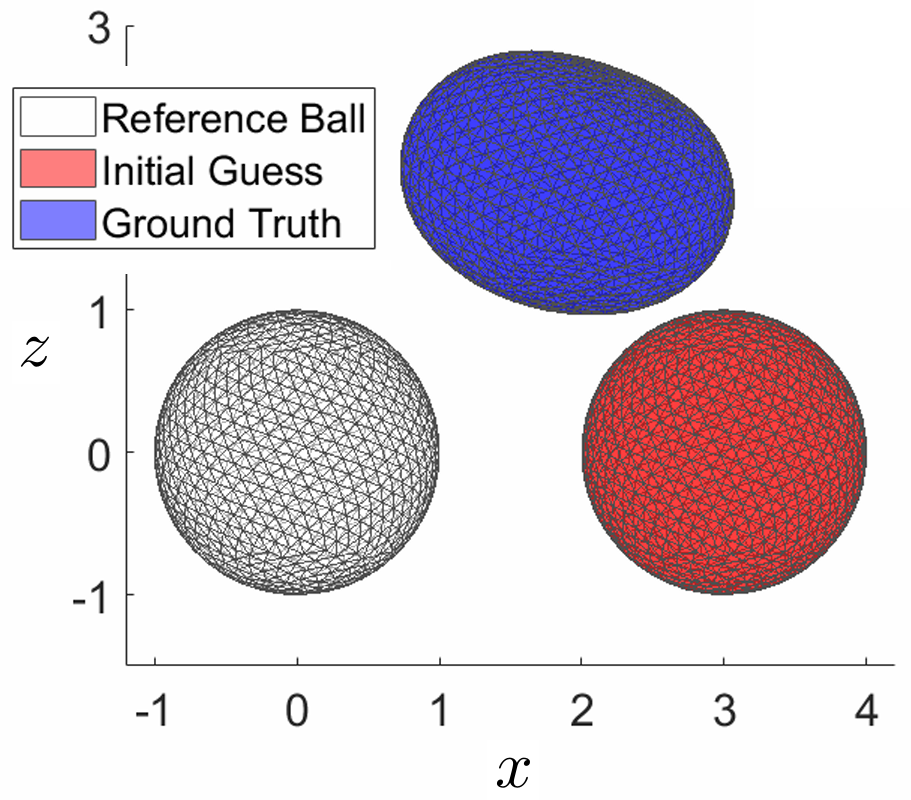}&&\includegraphics[height=3.9cm]{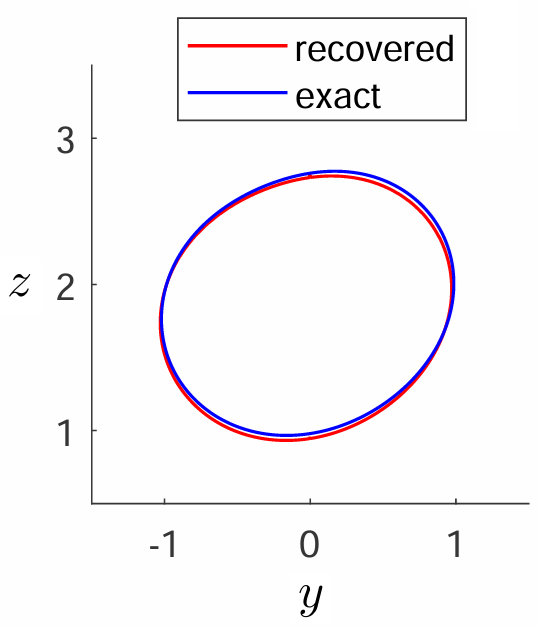}&&
\includegraphics[height=3.9cm]{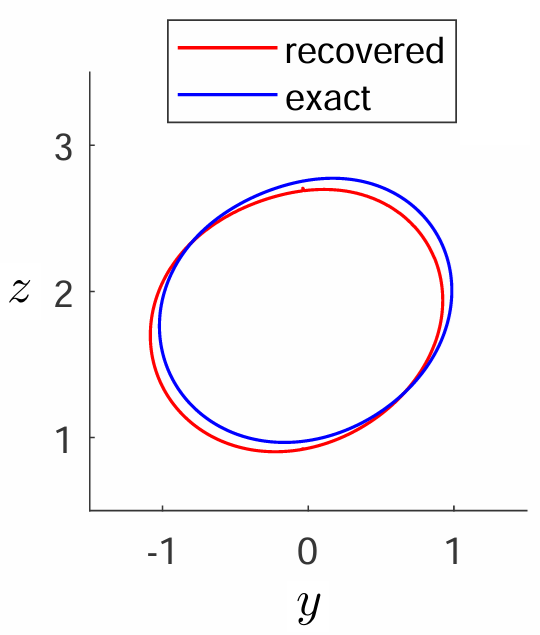}&&
\includegraphics[height=3.9cm]{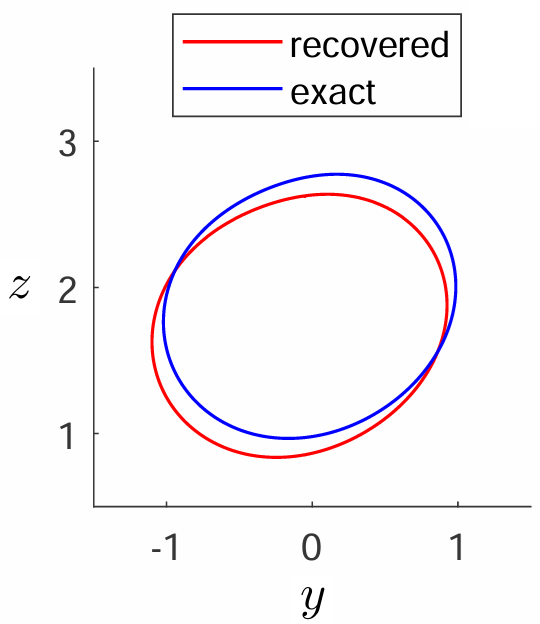}\\
&&\includegraphics[height=3.9cm]{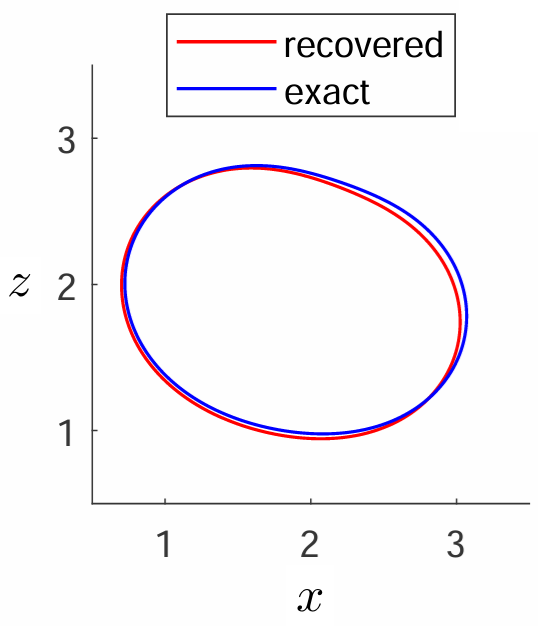}&&
\includegraphics[height=3.9cm]{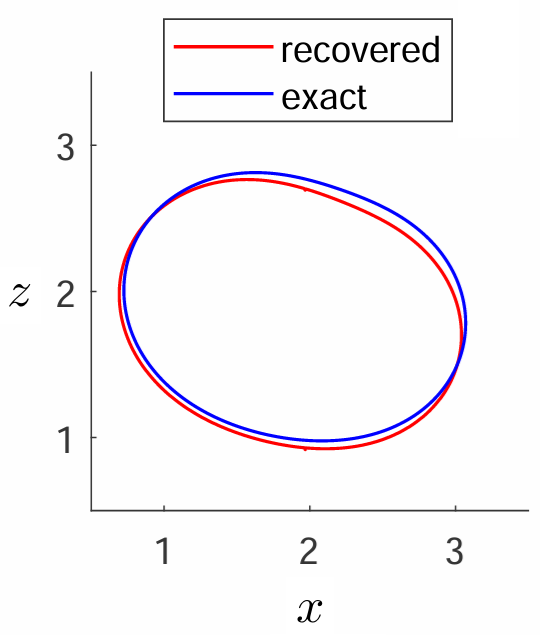}&&
\includegraphics[height=3.9cm]{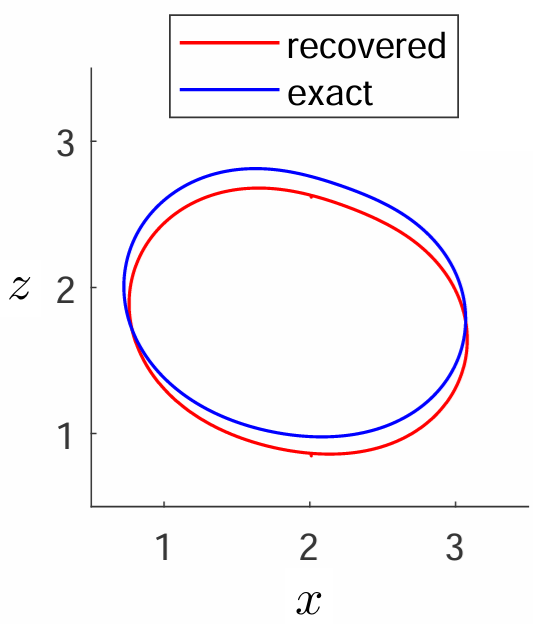}\\
&&\includegraphics[height=3.9cm]{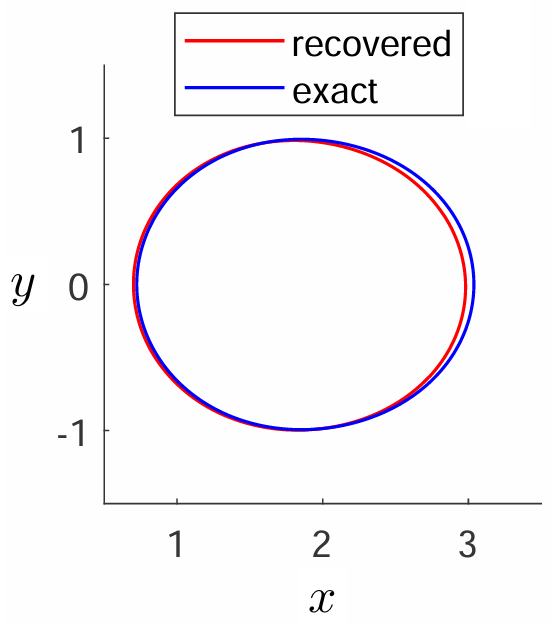}&&
\includegraphics[height=3.9cm]{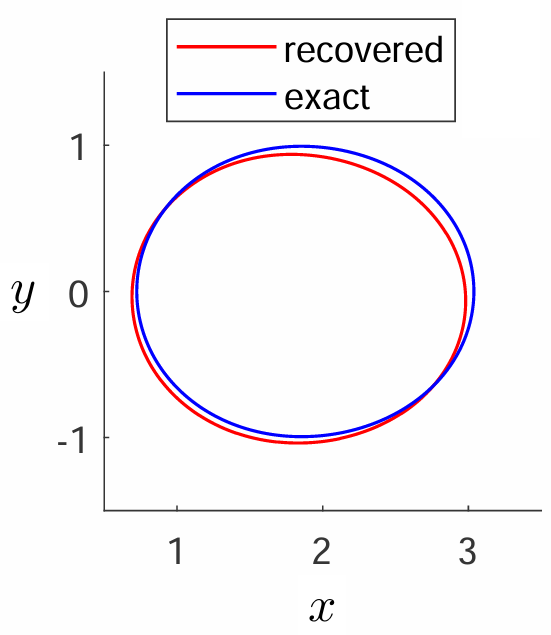}&&
\includegraphics[height=3.9cm]{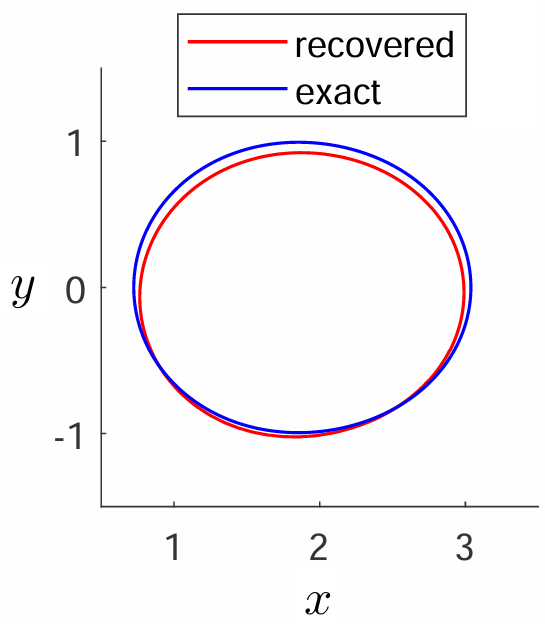}\\
&\phantom{a}& (a) $5\%$ &\phantom{a}& (b) $10\%$ &\phantom{a}& (c) $20\%$
\end{tabular}
\caption{
The cross sections by the planes $x=2$ (top), $y=0$ (middle) and $z=2$ (bottom) of the 3D reconstructions for the mean shape parameters computed using the epochs in $(9*10^4,10^5]$.}
\label{fig:MCMC:3D:reconstruction}
\end{figure}

\section{Conclusion}
In this work we have developed a rigorous numerical approach for reconstructing a sound-soft obstacle from phaseless far field measurements using DNNs. We have rigorously established the feasibility of the approach by providing expression rates of DNNs with the ReLU activation function for approximating the forward maps, via the concept of shape holomorphy. The analysis is based on variational formulations of the direct problems on a bounded domain (involving a nonlocal boundary condition), and can handle the case of piecewise analytic coefficients. The approach can be used directly to accelerate the posterior sampling arising from the Bayesian treatment of the inverse obstacle scattering problem. Numerically we observe significant speed-up in the Markov chain Monte Carlo sampling of the posterior distribution, which shows its significant potential for inverse scattering with phaseless data.

In terms of the practicality of the approach, there are multiple avenues for further research. First, the approach is specifically developed for the class of star-shaped obstacles, which is the primary focus of  existing uniqueness theory for phaseless inverse scattering. Nonetheless, it is of much interest to develop the mathematical theory and algorithms for a more class of obstacles, including nonsmooth (piecewise smooth), non-star-shaped obstacles as well as obstacles with disconnected components. One such extension is a class of parametric shapes by extending the ideas of the disk-to-domain maps  \cite{Arioli:2019:SBP,Juttler:2019:AFPD}. Second, the training of the surrogate model is fully supervised for the specific configuration. It is natural to explore more flexible strategies that can accommodate various challenging variations, e.g., different frequencies, incident directions, or background medium, for which pre-training / test-time adaptation are very promising.
\appendix

\section{Shape holomorphy for point source excitation}

In this part, we show the shape holomorphy of the forward map in the case of point source excitation. Specifically, consider the restriction of $v$ to $E:=B_{\rm o}\backslash(\overline{D}\cup\overline{B}\cup\overline{B_2})$.
Using the Dirichlet-to-Neumann maps $\Lambda^{\rm int}$ and $\Lambda^{\rm ext}$ for the Helmholtz equation on $B_2$ and $\mathbb{R}^d\backslash\overline{B_{\rm o}}$, $v$ satisfies
\begin{equation}\label{eq:bddomain:reformulation:v}\left\{\begin{aligned}
(\nabla\cdot\sigma(\Bx)\nabla + k^2\tau(\Bx)) v &= 0,\quad\mbox{in }E,\\
v &= 0,\quad\mbox{on }\partial D\cup\p B,\\
\partial_\nu (v-v^{\inc}) &= -\Lambda^{\rm int} (v-v^{\inc}),\quad\mbox{on }\partial B_{2},\\
\partial_\nu v &= -\Lambda^{\rm ext} v,\quad\mbox{on }\partial B_{\rm o},
\end{aligned}
\right.
\end{equation}
where the unit normal vector $\nu$ points outside of $B_{\rm o}\backslash \overline{B_2}$. Let $\mathcal{H}:=\{f\in H^{1}(E)\,:\, f|_{\p B\cup \p D}=0\}$.
The weak formulation of problem \eqref{eq:bddomain:reformulation:v} reads: find $v\in H^{1}(E)$ such that
\begin{equation}\label{eq:reform:v}
a(v,\zeta) = b(\zeta),\quad\forall \zeta\in \mathcal{H},
\end{equation}
with the sesquilinear form
$a(v,\zeta):=\int_{\p B_{2}}(\Lambda^{\rm int}v)\overline{\zeta}~\d s+\int_{\p B_{\rm o}}(\Lambda^{\rm ext}v)\overline{\zeta}~\d s + \int_{E}\sigma\nabla v\cdot\nabla \overline{\zeta} - k^2\tau v \overline{\zeta}~\d x$ and the linear form $ b(\zeta):=\int_{\p B_{2}}(\Lambda^{\rm int}v^{\rm i} + \partial_\nu v^{\rm i})\overline{\zeta}~\d s$. Let $A:\mathcal{H}\to \mathcal{H}'$ be the operator induced by the sesquilinear form $a$. Then we have the following well-posedness of problem \eqref{eq:reform:v}.
\begin{lemma}
$A:\mathcal{H}\to\mathcal{H}'$ is an isomorphism.
\end{lemma}
\begin{proof}
By \cite[Corollary 3.1]{Chandler:2008:WNEBTIS}, we have $\Re\left[(\Lambda^{\rm ext} u,{u})_{L^2(\p B_{\rm o})}\right]\ge0$.
Using the series expressions in polar and spherical coordinates for the solution to the Helmholtz equation in the 2D and 3D cases, respectively (see \cite[Chapter 2]{Nedelec:2001:AEE} and \cite{Chandler:2008:WNEBTIS} for circles and spheres, respectively), one can derive $\Re\left[ (\Lambda^{\rm int} u,{u})_{L^2(\p B_{2})}\right]\ge0$.
The rest of the proof is identical with that of Lemma \ref{lemma:iso:pw}.
\end{proof}		

Let $D_T=T(\widehat{D})$ for all $T\in\mathcal{T}$. We denote by $v_T$ the solution to problem \eqref{eq:reform:v} with $D=D_T$. We define, for all $T\in\mathcal{T}$, $\widehat{v}_T:=v_T\circ T$,  and $E_T:=B_{\rm o}\backslash( \overline{D_T}\cup\overline{B}\cup\overline{B_2})$ so that there holds
\begin{align*}
{v}_T\in \mathcal{H}_T&:=\{f\in H^{1}(E_T)\,:\, f|_{\p B\cup \p D_T}=0\},\quad \widehat{v}_T\in\mathcal{H},\,\forall T\in\mathcal{T}.
\end{align*}
Then $\widehat{v}_T$ satisfies
\begin{align}
\label{eq:varform:extended:v}	a_T(\widehat{v}_T,\zeta) &= b(\zeta),\quad\forall \zeta\in \mathcal{H},			\end{align}
with the sesquilinear form $a_T$ and linear form $b_T$ given respectively by
\begin{align*}
	a_T(v,\zeta):=&\int_{\p B_{\rm o}}(\Lambda^{\rm ext}v)\overline{\zeta}~\d s+\int_{\p B_{2}}(\Lambda^{\rm int}v)\overline{\zeta}~\d s \\
     &+ \int_{E}(\sigma\circ T)|J_T| J_T^{-1}J_T^{-\top}\nabla v\cdot\nabla \overline{\zeta} - k^2 (\tau\circ T)|J_T| v \overline{\zeta}~\d \mathbf{x},\\
    b_T(\zeta):=&\int_{\p B_{2}}(\Lambda^{\rm int}v^{\rm i} + \partial_\nu v^{\rm i})\overline{\zeta}~\d s,\quad\forall v,\zeta\in\mathcal{H}.
\end{align*}
Let $A_T:\mathcal{H}\to\mathcal{H}'$ be the operators induced by $a_T$. Similar to Lemma \ref{lem:isomorphism}: if $A_{\rm id}$ is an isomorphism, then so is $A_T$ for every $T\in\mathcal{T}$. Similar to Lemma \ref{lemma:welldefined:complex}: since $\mathcal{T}$ is compact, there exists a $\delta>0$ such that \eqref{eq:varform:extended:v} has the unique solution for all $T\in\mathcal{T}_\delta$.

\begin{theorem}\label{theorem:holomorphy:pointsource}
There exist $\delta>0$ and a holomorphic extension $F:\mathcal{T}_\delta\to H^{1}(E,\mathbb{C})$ of the domain-to-solution map $T\mapsto\widehat{v}_T$.
\end{theorem}
\begin{proof}
Since $\Lambda^{\rm int}$ and $\Lambda^{\rm ext}$ are independent of $T$, the proof of Theorem \ref{theorem:holomorphy:planewave} still applies: the Fr{\'e}chet derivative $F'(T)(H)$ of $F$ at $T\in\mathcal{T}_\delta$ is given by the unique solution $V_{T}(H)\in H^{1}(E,\mathbb{C})$ of the following problem:
${a}_T(V_{T}(H),\zeta)=g_{T,H}(\widehat{v}_T,\zeta)$ for all $\zeta\in\mathcal{H}$, with $g_{T,H}$ given in \eqref{eqn:g-TH}.
\end{proof}

\section{Posteriori approximation}\label{appendix:Posteriori approximation}
In this section, we derive an error bound on the approximate posterior distribution when using the neural network surrogate in place of the exact forward map. Such analysis is well established  \cite{YanZhang:2017,StuartTeckentrup:2018}. For an approximate posterior $\widetilde{\pi}$ to the exact $\pi$, with the shape parameters $\boldsymbol{\eta}$ in the parameter domain $[-1,1]^n$ (with $d_p$ being the parameter dimension), we bound the Kullback–Leibler (KL) divergence $D_{\mathrm{KL}}(\widetilde{\pi} \| \pi)$, defined by
$$
D_{\mathrm{KL}}\left(\widetilde{\pi} \| \pi\right):=\int_{[-1,1]^n} \widetilde{\pi}(\boldsymbol{\eta}) \log \frac{\widetilde{\pi}(\boldsymbol{\eta})}{\pi(\boldsymbol{\eta})} \mathrm{d} \boldsymbol{\eta}.
$$
The negative log-likelihood $\Phi$ and the approximate one $ \widetilde{\Phi}(\boldsymbol{\eta})$ are given respectively by
\begin{align*}
\Phi(\boldsymbol{\eta})&=\sum_{j=1}^3\frac{1}{2\sigma_j^2}\sum_{q=1}^{Q} F_{j,q}(\widetilde{\Om})+ \lambda R(\widetilde{\Omega}),\\
\widetilde{\Phi}(\boldsymbol{\eta})&=\sum_{j=1}^3\frac{1}{2\sigma_j^2}\sum_{q=1}^{Q} \widetilde{F}_{j,q}(\widetilde{\Om})+\lambda R(\widetilde{\Omega}),
\end{align*}
where $\widetilde{F}_{j,q}$ denotes the losses computed using DNN surrogates.
By Theorem \ref{theorem:nn:approx}, the difference between the DNN surrogate at the discrete nodes $\{\hBx_q\}_{q=1}^Q$ and the exact far-field data satisfies
\begin{equation}
    \max_{1 \leq q \leq Q} \left| U(\boldsymbol{\eta})(\hat{\mathbf{x}}_q) - \widetilde{U}_n(\boldsymbol{\eta}, \hat{\mathbf{x}}_q) \right| \leq C' n^{1-1/p}.
    \label{eq:pointwise_error}
\end{equation}

Then we can state the following error bound on the approximate posterior distribution $\widetilde{\pi}$ in terms of the Kullback-Leibler divergence.
\begin{theorem}
Let the functions $U$ and $\widetilde{U}_n$ satisfy Theorem \ref{theorem:nn:approx}. Then for the approximate posterior distribution $\tilde \pi$  and the exact one $\pi$, there esists a constant $C$ independent of $n$ such that
$$
D_{\mathrm{KL}}\left(\widetilde{\pi} \| \pi\right) \leq C n^{1-1/p}.
$$
\end{theorem}
\begin{proof}
Note that the KL divergence can be expressed as:
\begin{equation}
 D_{\mathrm{KL}}\left(\widetilde{\pi} \| \pi\right) = \mathbb{E}_{\pi}[\Phi(\boldsymbol{\eta})-\widetilde{\Phi}(\boldsymbol{\eta})] + \log\left(\frac{Z_{\pi}}{Z_{\widetilde{\pi}}}\right),
\label{eq:kl_decomposition}
\end{equation}
where $Z_{\pi}$ and $Z_{\widetilde{\pi}}$ are the normalization constants.
By Theorem \ref{theorem:nn:approx}, we have
\begin{align*}
|\Phi(\boldsymbol{\eta})-\widetilde{\Phi}(\boldsymbol{\eta})|
&\leq \sum_{j=1}^3 \frac{1}{2\sigma_j^2} \sum_{q=1}^Q |F_{j,q}(\widetilde{\Omega}) - \widetilde{F}_{j,q}(\widetilde{\Omega})|\\
&\leq \sum_{j=1}^3 \frac{QC}{2\sigma_j^2} n^{1-1/p} \leq C n^{1-1/p}.
\end{align*}
Considering the normalizing constant $Z_{\pi}$, we estimate that $ e^{-\|\widetilde{\Phi}(\boldsymbol{\eta})-\Phi(\boldsymbol{\eta})\|_{L^\infty
}}Z_{\widetilde{\pi}}\leq Z_{\pi}\leq e^{\|\widetilde{\Phi}(\boldsymbol{\eta})-\Phi(\boldsymbol{\eta})\|_{L^\infty
}}Z_{\widetilde{\pi}}$.
The second term can be further bounded by
\begin{align*}
\left|\log\left(\frac{Z_{\pi}}{Z_{\widetilde{\pi}}}\right)\right|
&= \left|\log\left(\frac{\int_{[-1,1]^n} e^{-\Phi(\boldsymbol{\eta})}d\boldsymbol{\eta}}{\int_{[-1,1]^n} e^{-\widetilde{\Phi}(\boldsymbol{\eta})}d\boldsymbol{\eta}}\right)\right|\\ &=\left|\log\left(\frac{\int_{[-1,1]^n} e^{-\widetilde{\Phi}(\boldsymbol{\eta})+(\widetilde{\Phi}(\boldsymbol{\eta})-\Phi(\boldsymbol{\eta}))}d\boldsymbol{\eta}}{\int_{[-1,1]^n} e^{-\widetilde{\Phi}(\boldsymbol{\eta})}d\boldsymbol{\eta}}\right)\right|  \\&\leq \sup_{\boldsymbol{\eta}\in  [-1,1]^n}|\widetilde{\Phi}(\boldsymbol{\eta}) - \Phi(\boldsymbol{\eta})|
\leq C n^{1-1/p}.
\end{align*}
Combining these two estimates yields the desired bound.
\end{proof}

\section{MCMC algorithm}\label{appendix:MCMC}

To explore the posterior distribution $\widetilde \pi$, we repeat the following three steps (1)-(3):
\begin{itemize}
\item[(1)]
Set the parameters $N\in\mathbb{N}$, $w>0$, $w_{\ell}>0$, $\lambda\ge0$ and $\sigma>0$, initialize $(\mathbf{a}^{(0)}, \mathbf{b}^{(0)})\in\mathcal{A}_N(B_1,B_2)$ and set $i=1$.
\item[(2a)] The $i$th iteration consists of two rounds: determine $\mathbf{a}^{(i)}$, and then $\mathbf{b}^{(i)}$. Draw independently $x_{j}^{(i)}\sim N(0,w^2)$ and $y_{j}^{(i)}\sim N(0,w_j^2)$ for all $j$, define $\Bx^{(i)}= (x_{j}^{(i)})_{j=1}^n$, $\By^{(i)}=(y_{j}^{(i)})_{j\in\mathbb{N}}$, and let
\begin{align*}
(\widetilde{\mathbf{a}},\widetilde{\mathbf{b}})=\begin{cases}
    (\mathbf{a}^{(i-1)}+\mathbf{x}^{(i)},\mathbf{b}^{(i-1)})&\mbox{for the first round (R1)},\\
    (\mathbf{a}^{(i)},\mathbf{b}^{(i-1)}+\mathbf{y}^{(i)})&\mbox{for the second round (R2)}.
\end{cases}
\end{align*}

\item[(2b)]
Let $\widetilde{u}_{\rm nn}^{\infty}$ and $\widetilde{v}_{\rm nn}^{\infty}$ be the  FNN predictions for $u^\infty[\widetilde{\Om}\cup B]$ and $v^\infty[\widetilde{\Om}\cup B]$, with $\widetilde{\Om} = \Om(\widetilde{\mathbf{a}},\widetilde{\mathbf{b}})$.
Define the acceptance rate $\alpha$ by
\begin{equation*}
\alpha=\begin{cases}
    {\widetilde{\pi}}/{\pi^{(i-1)}}&\mbox{in (R1)},\\
    {\widetilde{\pi}}/{\pi^{*}}&\mbox{in (R2)}
\end{cases}
\end{equation*}
with
\begin{equation*}
\widetilde{\pi}:=\exp\bigg(-\sum_{j=1}^3\frac{1}{2\sigma_j^2}\sum_{q=1}^{Q} F_{j,q}(\widetilde{\Om}) - \lambda R(\widetilde{\Om})\bigg).
\end{equation*}
For the acceptance, randomly draw $\alpha_0\sim \operatorname{Uniform}(0,1)$ and set
\begin{align*}
(\mathbf{a}^{(i)},\pi^{*})&=\begin{cases}
\ds (\mathbf{a}^{(i-1)},\pi^{(i-1)})\quad\mbox{if }\alpha<\alpha_0,\\[1mm]
\ds (\widetilde{\mathbf{a}},\widetilde{\pi})\quad\mbox{if }\alpha\ge \alpha_0\mbox{ in (R1);}
\end{cases}\\ (\mathbf{b}^{(i)},\pi^{(i)})&=\begin{cases}
\ds (\mathbf{b}^{(i-1)},\pi^{*})\quad\mbox{if }\alpha<\alpha_0,\\[1mm]
\ds (\widetilde{\mathbf{b}},\widetilde{\pi})\quad\mbox{if }\alpha\ge \alpha_0\mbox{ in (R2)}.
\end{cases}
\end{align*}

\item[(2c)] (Stopping criterion) Stop the iteration at $i=10^5$ and move to Step (3). Otherwise, increase $i$ by $1$ and move to Step (2a).

\item[(3)] (Reconstruction) Compute the reconstruction using the mean shape parameters of the MCMC trajectories of $(\mathbf{a}^{(i)},\mathbf{b}^{(i)})$.

\end{itemize}

\bibliographystyle{abbrv}
\bibliography{ref_survey}
	
\end{document}